\documentclass{article}
\usepackage{latexsym}
\usepackage{amssymb}
\font\gorditas = msbm8
\def\bbb#1{\hbox {{\gordas #1}}}
\def\errita{\hbox{\gorditas R}}

\font\gordas = msbm10 at 12pt

\def\bbb#1{\hbox {{\gordas #1}}}
\def\erre{{\bbb R}}
\def\ze{{\bbb Z}}
\def\theequation{\thesection.\arabic{equation}}

\newcommand{\R}{{\mathbb{R}}}

\newcommand{\cM}{{\mathcal{M}}}
\newcommand{\cC}{{\mathcal{C}}}

\newcommand{\noi}{{\noindent}}

\begin{document}
{\centerline {\large\bf Occupation Time Fluctuations in Branching
Systems$^1$}}
\begin{center}
\begin{tabular}{llll}
\quad D.A. Dawson$^2$ &\qquad L.G. Gorostiza$^3$&and &A.
Wakolbinger$^4$
\end{tabular}
\end{center}
\vglue.5cm
\noindent
\begin{abstract}
We consider particle systems in locally compact Abelian groups with
particles moving according to a process with symmetric stationary
independent increments and undergoing one and two levels of critical
branching. We obtain long time fluctuation limits for the occupation
time process of the one--and two--level systems.
 We give complete results
 for the case of finite variance branching, where the
fluctuation limits are Gaussian random fields, and partial results
 for an example of infinite variance branching, where the fluctuation
limits are stable random fields. The asymptotics of the occupation
time fluctuations are determined by the Green potential operator $G$
of the individual particle motion and its powers $G^2, G^3$, and by the 
growth
as $t\rightarrow\infty$ of the operator $G_t=\int^t_0T_sds$ and its
powers, where $T_t$ is the semigroup of the motion. The results are
illustrated with two examples
of motions: the symmetric $\alpha$--stable L\'evy process in
$\erre^d$ $(0<\alpha\leq2)$,
 and the so called $c$--hierarchical random walk
in the hierarchical group of order $N$
($0<c<N$).
We show that the two motions have analogous
asymptotics of $G_t$ and its powers that depend on
an order parameter $\gamma$
for their transience/recurrence behavior. This parameter
is $\gamma=d/\alpha-1$ for the
$\alpha$--stable motion, and $\gamma=\log c/\log (N/c)$ for the
$c$--hierarchical random walk.
As a consequence of these analogies, the asymptotics of the occupation time 
fluctuations
of the corresponding branching particle systems are also analogous.
In the case of the $c$--hierarchical random walk, however,
the growth of $G_t$ and its powers is modulated
by oscillations on a logarithmic time scale.
\end{abstract}
\vglue1cm
\noindent
\footnote{\kern-.55cm
$^1$ Research partially supported by a Max Planck award to D.A. Dawson,
 CONACyT grant 27932--E  (Mexico), and a DFG\\[-.2cm]
\hglue.3cm  grant (Germany).\\
$^2$ The Fields Institute, Toronto, Canada, e.mail: 
ddawson@fields.utoronto.ca\\
$^3$ Centro de Investigaci\'on y de Estudios Avanzados, A.P. 14--740,
 Mexico 07000 D.F., Mexico,\\[-.2cm]
\hglue.3cm e.mail: gortega@servidor.unam.mx\\
$^4$ Johann Wolfgang Goethe--Universit\"at, Frankfurt am Main, Germany, 
e.mail: wakolbin@math.uni-frankfurt.de}
\newpage

\noi
{\bf Table of contents}

\bigskip

\noi
{\bf 1.  Intoduction}

\bigskip

\noi
{\bf 2.  General Notions, Results and Comments}

\noi
2.1. The individual motion: powers of the Green potential, strong and weak
 transience
\vglue-.2cm

\noi
2.2. Main results: Occupation time fluctuations for $k$-level critical 
binary
branching particle systems \\[-.2cm]
\hglue.7cm ($k = 0, 1, 2$)
\vglue-.2cm

\noi
2.3. Comments on the assumptions and the results
\vglue-.2cm

\noi
2.4. Order of transience and recurrence, asymptotics of powers of $G_t$, 
and
some special growth functions
\vglue-.2cm

\noi
2.5. Infinite variance branching
\vglue-.2cm

\noi
2.6. Occupation time fluctuations of superprocesses

\bigskip

\noi
{\bf 3.  Two examples of individual motions:
Symmetric $\alpha$-stable
process and $c$-hierarchical}
\vglue-.2cm
{\bf random walk}

\noi
3.1. Transience/recurrence properties of the motions
\vglue-.2cm

\noi
3.2. Occupation time fluctuations limits
\vglue-.2cm

\noi
3.3. Comments on the results

\bigskip

\noi
{\bf 4.  Definitions of constants and functions for the examples}

\noi
4.1. Notation for $\alpha$-stable motion
\vglue-.2cm
\noi
4.2. Notation for $c$-hierarchical random walk

\bigskip

\noi
{\bf 5.  Proofs}

\noi
5.1. Asymptotics of the powers of $G_t$
\vglue-.2cm

\noi
5.2. Main results
\vglue-.2cm
\noi
5.3. Examples
\vglue-.2cm

\noi
5.4. Conditions for the results on infinite variance branching

\bigskip

\noi
{\bf Appendix}

\noi
A.1. Background on 1- and 2-level branching systems
\vglue-.2cm

\noi
A.2. The Palm formula
\vglue-.2cm

\noi
A.3. Tree representation of the Palm measures of $R_t^1$ and
$R_\infty^1$
\vglue-.2cm

\noi
A.4. Second and third moments of $R_\infty^1$

\bigskip

\noi
{\bf Acknowledgements}

\bigskip

\noi
{\bf References}

\newpage
\noindent
{\bf 1. INTRODUCTION}
\vglue.5cm
Consider a particle system described by a random counting measure
$X_t$
on a space of sites $S$, with the same intensity measure
$EX_t$ for all
$t$ which is denoted by $\rho$, and such that $X_t$
converges in distribution as $t \to \infty$ towards an
equilibrium state  which also has the same intensity $\rho$, i.e., the 
system is {\it persistent}.
Then under mild conditions the {\em occupation time fluctuation}
$Y_t
=
\int_0^t (X_s-\rho)
ds$ obeys a law of large numbers, i.e.
$\displaystyle\frac{1}{t}Y_t \to 0$ as $t
\to \infty$ (see e.g. M\'el\'eard and Roelly$^{(32)}$ for a branching 
particle system in $\erre^d$).
The question for which norming $a_{t}$ does a non--trivial limit of
$\displaystyle\frac{1}{a_{t}}Y_t$
exist in distribution as $t\rightarrow\infty$
and what is the limiting random field depends on more specific
properties of
the system. In this paper we  investigate this question for
branching particle
systems
in locally compact Abelian groups,
where the individual particle motion is a process with
symmetric
stationary independent increments,  and for the
so called ``2--level'' branching systems in which not only the
individual particles but also whole families of particles undergo
critical branching. Multilevel branching systems
were introduced by Dawson and Hochberg$^{(9)}$, and they have been studied 
by several authors: Dawson et al$^{(10)}$, Gorostiza$^{(16)}$, Gorostiza et 
al$^{(17)}$, Greven and Hochberg$^{(24)}$, Hochberg$^{(25)}$, Hochberg and 
Wakolbinger$^{(26)}$, Wu$^{(40)}$.

For a transient motion and finite variance branching,
the simple branching  particle system converges as $t \to \infty$ to a
Poisson system of
independently evolving ``clans'', each of which contributes to the
occupation time. The asymptotics of the occupation time
fluctuations should be determined by the
space--time correlations within single clans. Also, the growth of the
occupation time fluctuations  as $t \to \infty$ should depend on
whether there are long time dependencies caused by recurrent visits
of single clans to bounded sets.

Let us first recall some known results. For a critical finite variance
 branching
Brownian
system $X_t$ in $\erre^d$, started off from a Poisson system with
Lebesgue intensity $\lambda$, the right norming  $a_t$ for the
occupation
time
fluctuation is $t^{3/4}$ for $d=3$, $(t \log t)^{1/2}$ for $d=4$,
and $t^{1/2}$ for $d>4$ (see
Cox and Griffeath$^{(5)}$, and also
Iscoe$^{(27)}$, where the corresponding superprocess scenario is treated).
We refer to this as the 1--level branching case.

The same normings appear for the occupation time fluctuations of
Poisson systems
of Brownian particles without branching, which we call 0--level
systems, but
two dimensions lower (Cox and Griffeath$^{(4)}$,
Deuschel and Wang$^{(14)}$), and we shall
see that they also appear in the $2$--level branching case, but now
two dimensions higher than in the $1$--level case.

There is an apparent relation between the critical dimension for
transience
of the motion
and the critical dimension for the classical
$t^{1/2}$--norming of the occupation time fluctuations:
for Brownian particle systems without branching,
2 is the critical dimension
above which the occupation time fluctuations have the classical
norming
and it
is also the dimension above which the particle motion is transient.
In the
1--level branching case, 4 is the critical dimension above which the
occupation time
fluctuations have the classical norming and it is also the dimension
above which
the equilibrium clans are transient (in the sense that they eventually
leave each
bounded region of ${\erre}^d$ forever (St\"ockl and Wakolbinger$^{(38)}$). 
We shall see that
an
analogous result  holds for the $2$--level case.

One of our main objectives is to put these results in a general
context for branching systems in locally compact Abelian groups,
which clarifies the role played by the Green potential operator $G$
of the
particle
motion and
the (operator)
powers  $G^k$ of $G$ in relation with  the various levels
$k=0,1,2,\ldots$ of branching. A key role will be played by the level $k$ 
transience and recurrence properties of the motion, $k=0,1,2,\ldots$ 
defined in Section 2.
In
this paper we will treat only branching levels $k=0,1,2$, but
the results show a pattern which allows one
to guess what the form of the results would be for systems with higher 
levels of
branching.
The analysis of 2--level systems is considerably more difficult than that
of 1--level systems due to the dependencies among
the particles caused by the simultaneous branching of families of 
particles.

Roughly speaking, the bigger the branching level $k$ is, the more
long range dependencies
 are introduced into the system. These dependencies increase the mass
 fluctuations in the system, whereas a strong spreading out of mass
 by the particle motion has a smoothing effect on the mass
 fluctuations.

 In the case of finite variance branching (where we will restrict for
 simplicity to binary branching)
it turns out that finiteness of $G^k$
corresponds to existence of the $k$--level branching equilibrium
clans, and then
$G^{k+1}=\infty$  corresponds to a long time dependence of the visits
of single $k$--level
clans to bounded sets. In the latter case,
a crucial feature of these models is that
the growth of
the
$(k+1)$--st power of the operator
$G_t = \int_{0}^{t}T_{s}ds$ as $t\rightarrow\infty$,
where $T_{t}$ denotes the semigroup
of the individual motion, determines the
right norming $a_t$ for the occupation time fluctuations
of the $k$--level system, whereas for
finite $G^{{k+1}}$ their right norming is the classical $t^{1/2}$.
In the case of finite $G^{k+1}$ the covariance of the limiting
Gaussian field of the occupation time fluctuations
of the $k$--level system contains terms
induced by direct ancestry  and by level $j$--relationship,
$1\le j \le k$.

 For infinite variance
``$(1+\beta)$--branching'' $(0<\beta<1)$, we have so far results only for 
the
 case with $t^{1/(1+\beta)}$--norming.
Then the norming
is determined  by the highest level of branching, and the
occupation time fluctuations converge to stable random fields.

For the superprocess limits of the 1-- and the 2--level branching particle 
systems the corresponding occupation time fluctuation results are
basically analogous to those for the particle systems
and even somewhat simpler.

We will focus on two examples of
particle motions $W_t$: the symmetric $\alpha$--stable
L\'evy process in ${\erre}^d, 0<\alpha\leq2$
(including Brownian motion, $\alpha=2$), and a  ``hierarchical'' random
walk in $\Omega_{N}$, the hierarchical group of order $N$, which is
a direct sum of a countable number of copies of $\ze_{N-1}$, i.e.
$\Omega_N=\{x=(x_1,x_2,\ldots)\,\,\,|x_i\in\{0,1,\ldots,N-1\},x_i
\neq 0$ except for finitely many i$\}$.

 For the symmetric
$\alpha$--stable process in ${\erre}^d$  and $k\geq 0$,
$G_t^{k+1}$ has a power growth in $t$ if
$d/(k+1) < \alpha$, a logarithmic growth if
$\alpha = d/(k+1)$, and $G^{k+1}$  is finite if
$\alpha <d/(k+1)$.

For the  random walk in $\Omega_{N}$
 we consider a probability of jumping a distance $i$
proportional to
$(c/N)^{i-1}$,
where $c$ is a constant such that $0<c<N$.
(Here, the distance between two elements $(x_i)$ and $(y_i)$ of $\Omega_N$ 
is defined to be the highest index $i$ for which $x_i$ and $y_i$ are 
different).
Since the random walk is characterized by this ``mobility parameter''
$c$,
we will call it
the
$c$--{\it hierarchical random walk} in $\Omega_{N}$.
 In this case, for $k\geq 0$,
$G_t^{k+1}$ has a power growth in $t$ if
$c < N^{k/(k+1)}$, a logarithmic growth if
$c = N^{k/(k+1)}$, and $G^{k+1}$ is finite if $c >
N^{k/(k+1)}$. In this example the growths of $G_t$ and its powers have
also  oscillating modulations in a logarithmic time scale.
However, the oscillations vanish in the cases of logarithmic growth.
Oscillatory phenomena have also been observed in
another class of random walks on groups which are direct sums of a 
countable number of copies of a discrete group
(Cartwright$^{(2)}$), and in random walks on the Sierpi\'nski graph
(Barlow and Perkins$^{(1)}$, Grabner and Woess$^{(23)}$
and references therein). A basic reference for random walks on Abelian 
groups is the paper by Kesten and Spitzer$^{(29)}$.

The values of the parameters $\alpha$ (or $d$) and $c$
which  correspond to logarithmic growths of the powers
$G^{k+1}_t$ will sometimes be referred to as ``critical'', since they
 are boundaries between intervals with different regimes of the
longtime behavior of the systems.

We will show that the asymptotics of the occupation time fluctuations
of the branching systems
are analogous for the two examples.
A key role is played by a constant $\gamma >-1$, which is
$\gamma=d/\alpha-1$ for the $\alpha$-stable system in $\erre^d$, and
$\gamma=\log c/ \log (N/c)$
for the $c$--hierarchical system in $\Omega_{N}$.
For $\gamma \in (-1,0)$,  $G_{t}$ grows like $t^{-\gamma}$, for
$\gamma = 0$,  $G_{t}$ grows logarithmically, and  for $\gamma >0$,
$G-G_{t}$ decays like $t^{-\gamma}$. Thus, $\gamma$ is an order
parameter for the transience/recurrence behavior of the motion, and we
will refer to the three above mentioned cases as {\em recurrence of
order --$\gamma$, critical recurrence}, and {\em transience of
order} $\gamma$. It turns out that $G_t^{k+1}$  grows like
$t^{k-\gamma}$ if
$k > \gamma$, grows logarithmically if
$k=\gamma$, and $G^{k+1}$ is finite if $k < \gamma$.

For the symmetric $\alpha$--stable
processes in ${\erre}^d$,  $\gamma$ is restricted to $[(d-2)/2, \infty)$,
whereas for the $c$--hierarchical random walks in $\Omega_N$,
$\gamma$ can range
over the entire interval $(-1, \infty)$.
In this sense
the hierarchical random walks are a richer class of models.
Choosing
$c=N^{1-\alpha/d}$, the $c$--hierarchical random walk in $\Omega_{N}$
has the same
order of transience/recurrence as the $\alpha$--stable process in 
${\erre}^d$,
and this allows to think about non--integer dimensions $d$.

The idea of using hierarchical systems as models of systems with
non--integer dimension was first introduced in the context of
statistical physics. A model of ferromagnetic behavior involving the
case of
\newpage
\noindent
$N=2$ is known as Dyson's hierarchical model and has been
used by Sinai$^{(36)}$ as a framework in which to carry out a
rigorous renormalization group analysis following the ideas of
Wilson$^{(39)}$. In the case of ferromagnetism,
$4$ is the critical dimension and the hierarchical group has been used
to study large scale fluctuations near the critical point in 
$4-\varepsilon$
dimensions. In the case of 1--level branching Brownian motion in
${\erre}^d$, the dimension $2$ is the critical dimension for the
persistence/extinction dichotomy. In Dawson and Greven$^{(7)}$, 1--level
hierarchical
branching random walks (indexed by a sequence $(c_{j})$ rather than
just one parameter $c$) have been analyzed in the so called  mean--field 
limit
 $N\rightarrow\infty$ in the ``nearly 2--dimensional analogue''.
 Since the dimension 4 is the critical dimension for the occupation
 time fluctuations of 1--level branching Brownian motions and for the
 long time behavior of 2--level branching Brownian motions in ${\erre}^d$,
 it is conceivable that the hierarchical mean--field limit can be used to
 carry out a similar analysis of these phenomena ``near dimension 4''.

The outline of the paper is as follows:
Section 2 presents the general results
for the particle systems, and the results for the superprocesses are 
mentioned.
Sections 3 and 4 are devoted to
the two above mentioned examples of  individual motions.
Section 4 contains a list of
 constants and functions that appear in the results of Section 3.
The proofs are
given in section 5. An Appendix contains definitions
and background on
1-- and 2--level branching systems, and some tools.
\vglue1cm
\noindent
{\bf 2. GENERAL NOTIONS, RESULTS AND COMMENTS}
\vglue.5cm
\noindent
{\bf 2.1. The individual motion: powers of the Green potential,
strong and weak transience}

\vglue.5cm
We consider as a space of sites  a locally compact Abelian
group $S$
with Haar measure $\rho$, and as individual particle  motion
a process
 $W_t$  with stationary independent increments.
We assume that for each $s>0, W_s - W_0$ has
 a strictly positive symmetric density with respect to $\rho$,
and that
$W_t$ has
c\`adl\`ag paths.

Let us fix some notation.
The function spaces ${\cal C}_c(S), {\cal C}_\tau(S), {\cal C}^+_c(S), 
{\cal C}^+_\tau(S)$, and the measure spaces ${\cal M}_\tau(S), {\cal 
N}_\tau(S)$, where $\tau$ is a reference function on $S$, are defined in 
the Appendix.
For $\mu\in {\cal M}_\tau(S)$  and $\varphi \in {\cal C}_\tau(S)$, we
write $\langle \mu, \varphi\rangle = \int \varphi d\mu$, and we
denote
$$
(\varphi,\psi)_{\rho} = \int_S \varphi \psi \, d\rho,\qquad
\varphi,\psi\in
{\cal C}_\tau(S).
$$

We designate by $T_t$ the semigroup of $W_t$. Note that $\rho$ is
$T_t$--reversible, i.e.,
$(\varphi, T_t \psi)_{\rho} = (T_t \varphi,
\psi)_{\rho}$ for all $\varphi, \psi \in {\cal C}_\tau(S)$ and $t > 0$, 
which
implies that $\rho$ is $T_t$--invariant for each $t>0$.

The {\em
Green potential operator} $G$ of the motion is defined by
$$
G\varphi =  \int_0^\infty T_t\varphi \, dt, \quad
\quad \varphi \in {\cal C}_c(S).
$$
Together with $T_{t}$, also $G$ is self--adjoint with respect to
$\rho$, so
$$
\langle \rho,
(G\varphi)(G\psi) \rangle = (\varphi, G^2\psi)_{\rho}.
$$
Observe that the (operator) powers of $G$ are given by
\vglue.2cm
\noindent
$\hfill
G^k \varphi = \displaystyle\frac{1}{(k-1)!}
\displaystyle\int_0^\infty t^{k-1}T_t\varphi \, dt,
\quad
k \ge 2.
\hfill (2.1.1)
$\vglue.2cm
\noindent

The quantities $G^{k+1}(x,B) := G^{k+1}1_B(x),  x \in S,  B$ a
measurable subset of $S$, $k = 0,1,...$, have an interpretation in
terms of
occupation times of a mass flow with continuous birth of mass,
which will be helpful later on in the genealogical picture of 1-- and
2--level branching systems. Consider the case
$k=1$ and imagine an initial ``parent'' unit mass at
$x\in S$, which evolves according to the flow $T_t$.
This parent mass
generates
its own amount of ``daughter'' mass at its own site
continuously at rate 1, and this daughter mass is again transported
by the
flow $T_t$. Then $G^2(x,B)$ is simply the total occupation
time of the daughter mass in
$B$; this is immediate from the semigroup property.
To interpret $G^{k+1}(x,B)$ in a similar way, imagine an initial
``k--level''
unit mass at $x$  which evolves according to the flow $T_t$. For
every $j =
k,\dots,1$,   the j--level mass generates its own amount of $(j - 
1)$--level
mass, which is again transported by the flow $T_t$. Then
$G^{k+1}(x,B)$
results as the
total occupation time of 0--level mass in
$B$.

We  use the notation $||\cdot||$ for the supremum norm.
\vglue.5cm
\noindent
{\bf Definition 2.1.1.} {\rm (a) For $k \ge 0$, we say that
$W_t$
is {\em level k transient } if
$$
||G^{k+1}\varphi|| <\infty \,\,\mbox{ for}\, \varphi \in {\cal C}^+_c(S),
$$
and {\em level k recurrent} if
$$
G^{k+1}\varphi \equiv \infty \,\,\hbox{\rm for}\,\, \varphi \in
{\cal C}^+_c(S),\,
\varphi\neq\, 0.
$$
(b) For $k \ge 0$, we say that
$W_t$
is {\em level $k$ strongly transient } if it is level $k+1$
transient,
and  $W_t$ is {\em level k weakly transient} if it is level $k$ transient 
and
level $k+1$ recurrent.}
\vglue.5cm
Note that level $0$ transience and recurrence are (because of the
assumed
irreducibility of $W_t$) just the ordinary
notions of
transience and recurrence,  and note also
that level $0$ strong and weak transience coincide with the
notions
of strong and weak  transience as defined, e.g., in Port and 
Stone$^{(33)}$.
Clearly, level $k$ transience implies level $j$
transience for $j < k$, and
level $k$ recurrence implies level $j$ recurrence for $j > k$.
In terms of the interpretation given above, level $k$ transience
(resp. recurrence)
means that a $k$--level parent unit mass sends up to time infinity a
finite (resp. infinite) amount of 0--level daughter mass to any
bounded
region.
\vglue.5cm
\noindent
{\bf Definition 2.1.2.}  {\rm (a) We define the operator}
$$
G_t\varphi  = \displaystyle\int_0^t T_s\varphi \,ds,\quad \varphi\in
{\cal C}_\tau(S),\quad t>0,
$$
and denote by $G_t^k$  the $k$--th (operator) power of $G_t$, $k \ge
2$.

\noindent
(b)
For transient motion, we define the the bilinear form
$$
R_t(\varphi,\psi) = \left(\varphi,
(G-G_t)\psi
\right)_\rho,\quad \varphi, \psi\in {\cal C}_\tau(S),\quad t>0.
$$
Note that for each $k\geq 1$ and $\varphi\in{\cal C}^+_c(S),\varphi\neq
0$,
$\int^{\infty}_0t^{k-1}R_t(\varphi,\varphi) dt<\infty$
(${\rm resp.} = \infty$) if the motion is level $k$ transient (resp.
level $k$ recurrent).
\vglue.5cm
\noindent
{\bf Definition 2.1.3.} (a)
Let $H$ and, for each $t>0, H_t$ be positive definite bilinear forms on 
${\cal C}_c(S)$, and let $f:[T,\infty)\rightarrow\erre_+$ for some $T>0$. 
We write $H_t\sim f_tH$
if
$$
{1\over f_t} H_t(\varphi,\psi)\rightarrow H(\varphi,\psi)
\quad\mbox{\rm as}\quad t\rightarrow\infty,
\,\,\varphi,\psi\in{\cal C}^+_c(S).$$
If $Q_t$ is a linear operator from ${\cal C}_c(S)$ into ${\cal C}_\tau(S)$, 
the notation $Q_t\sim f_tH$ means that
$H_t\sim f_tH$ holds for $H_t(\varphi,\psi):=(\varphi,Q_t\psi)_\rho$.
\vglue.2cm
\noindent
(b) We call $f_t$ a {\it growth function} if it is
 increasing and
$\lim\limits_{t\rightarrow\infty}f_t=\infty$.
\vglue.5cm
Growth functions will be used to characterize the growth of $G_t$ and its 
powers. The growth functions we shall encounter in the examples are of the 
form $f_t=t^\zeta h_t$, where $\zeta\in(0,1)$ and $h_t$ is either 
identically
equal to 1 or a slowly oscillating function, or $f_t=\log t$.
Part (a) of
 the definition will also be used to characterize the ``order of
 transience'' of $W_t$ in terms of $R_t$ with $f_{t} \to 0.$
\vglue.5cm
\noindent
{\bf 2.2. Main results:
Occupation time fluctuations for $k$--level critical binary branching
particle systems ($k=0,1,2$)}
\vglue.5cm
We consider the following  particle systems in $S$ with the
individual particles moving independently
 according to the process $W_t$:

\begin{itemize}
\item[(i)]{\it {\rm 0}--level system}: The system starts off from a
Poisson
system  with
intensity $\rho$.
\item[(ii)]{\it {\rm 1}--level branching system}:
 The motion is transient, the system starts off from a Poisson
system with intensity $\rho$, and
the particles undergo critical binary branching at rate $V$. We
recall that
this system has a (infinitely divisible) ``Poisson type'' equilibrium
(in the sense of Liemant et al$^{(31)}$, section
2.3), and we denote by $R^1_{\infty}$ its canonical measure, which we
call {\it  equilibrium canonical measure}.
Note that $R^1_{\infty}$ has intensity $\rho$ (see the Appendix, (A.1.11)).
\item[(iii)]{\it {\rm 2}--level branching system}:
 The motion is strongly transient, the system starts off from a Poisson
system of
``2--level particles'' with intensity $R^1_{\infty}$, individual particles
undergo
critical binary branching at rate $V_1$ and clans undergo critical
binary
branching at rate $V_2$.
 Note that $R^1_{\infty}$ is an invariant measure for the
$1$--level dynamics, just as $\rho$ is an invariant measure for the
$0$--level dynamics (however,
$R^1_\infty$ is not reversible for the 1--level dynamics).
\end{itemize}
\vglue.2cm
We refer the reader to Gorostiza$^{(16)}$, and Hochberg and 
Wakolbinger$^{(26)}$
 for a detailed description of the dynamics of 2--level branching systems.
The necessary background for the present paper is given in
the Appendix and should be consulted as the need arises.
\bigskip
\par
For the three particle systems above,
$X_t$ stands for the empirical measure of the locations of all the
particles present
at time
$t$.
Thus, $X_t$ is a random point measure on $S$.
In  the 2--level case $X_t$ corresponds to the {\it aggregated}
system,
i.e., the particles are counted as ``1--level particles''
independently of what
clans  (``2--level particles'') they belong to.
In each one  of  these systems the
intensity is preserved, i.e., $EX_t=\rho$ for all $t> 0$,
as a consequence of  the initial conditions and  the
criticality in the branching cases. We consider the
{\it  occupation time fluctuation}, which is  the random signed
measure
$Y_t$ on $S$ defined by
$$
Y_t=\displaystyle\int^t_0(X_s-\rho)ds,\quad t>0.
$$
\vglue.2cm
The following theorems describe the asymptotic distribution of $Y_t$
as
$t\rightarrow\infty$ for the $k$--level systems, $k=0,1,2,$ described 
above.
All the convergence assertions  are
understood to be in distribution as $t\rightarrow\infty$,  all the
test functions belong to ${\cal C}^+_c(S)$, i.e., all random fields are 
considered over ${\cal C}^+_c(S)$.
The results for the 0--level system are basically
known in special cases, but we include them for completeness and
because they are the initial step in the multilevel ladder.
\vglue.5cm
\noindent
{\bf Theorem 2.2.1.}  Let $X_t$ be the  {\rm 0}--level
system.
\begin{itemize}
\item[(a)]  If the motion $W_t$ is transient,
then  $t^{-1/2}Y_t$
converges
to a Gaussian field with covariance
functional
$2(\varphi, G\psi)_{\rho}$.
\item[ (b)]
If the motion $W_t$ is recurrent with $G_t \sim f_tH$ for some
growth
function
$f_{t}$, then
$(\int_0^t f_s\, ds)^{-1/2}Y_t$ converges
to a Gaussian field with covariance
functional
$2H(\varphi, \psi)$.
\end{itemize}
{\bf  Theorem 2.2.2.}  Let $X_t$ be the {\rm 1}--level branching
system.
\begin{itemize}
\item[ (a)]
If $W_t$ is strongly transient,
then  $t^{-1/2}Y_t$ converges
to a Gaussian field with covariance functional
$$(\varphi,(2G+VG^2)\psi)_{\rho}=2\left(\varphi,\left(I+{V\over
2}G\right)G\psi
\right)_{\rho}.$$
\item[ (b)]  If  $W_t$ is
weakly transient with $G_t^2 \sim f_tH$ for some growth function
$f_{t}$, then
$(\int_0^t f_s\, ds)^{-1/2}Y_t$ converges
to a Gaussian
field with covariance
functional $VH(\varphi, \psi)$.
\end{itemize}
{\bf Theorem 2.2.3.}  Let $X_t$ be the {\rm 2}--level branching
system.
\begin{itemize}
\item[ (a)]
If  $W_t$ is level $1$ strongly
 transient and if there exists $\delta > 5/2$  such
that
\vglue.3cm
\noindent
$\hfill
||T_{t}\varphi||=O(t^{-\delta})\quad
\mbox{\rm as}\quad t\rightarrow\infty,\,\, \varphi\in{\cal C}^+_c(S),
\hfill (2.2.1)$
\vglue.3cm
\noindent
then  $t^{-1/2}Y_t$ converges
to a Gaussian field with covariance
functional
$$\left(\varphi, \left(2G + (V_1+V_2)G^2 + \frac 12 V_1V_2G^3\right)
\psi\right)_{\rho}= 2\left(\varphi,\left(I+{V_2\over
2}G\right)\left(I+
{V_1\over 2}G\right)G\psi\right)_{\rho}.$$
\item[(b)]
If  $W_t$ is level $1$ weakly transient with $G_t^2G\sim
f_tH$
for some growth function
$f_{t}$, then
$(\int_0^t f_s\, ds)^{-1/2}Y_t$ converges to a Gaussian
field with covariance
functional $ \displaystyle{1\over 2} V_1V_2  H(\varphi, \psi).$
\end{itemize}
\vglue.5cm
\noindent
{\bf 2.3. Comments on the assumptions and the results}

\begin{enumerate}
\item In  Subsection 2.4 we will give  conditions on the
motion $W_t$ which imply
the  growth assumptions of  Theorems 2.2.2(b) and 2.2.3(b).
The $\alpha$--stable motion  fits into this
framework.
\item
We do not know if condition (2.2.1) holds in general for
level 1 strongly transient motion. We will show, however, that
it holds for the motions in the examples.
 \item We  give an explanation of the second moment structure
ocurring in Theorem 2.2.2(a).
Let
${\underline R}_{\infty}^1$ be the historical counterpart
of the canonical equilibrium measure $R^1_{\infty}$
(which is a time--shift
invariant measure on the space of clans ranging from time
$ - \infty $ to time
$ + \infty $)
(see Dawson and Perkins$^{(12)}$, Sections 5.4.3, 5.4.4). Firstly, we 
observe
that the second moment measure of $R_{\infty}^1$ is given by (see the
Appendix, (A.1.12))
\vglue.3cm
\noindent
$\hfill
\displaystyle\int
\langle \mu, \varphi \rangle \langle \mu, \psi \rangle R_\infty^1
(d\mu) =
\left(
\displaystyle\varphi, \left(
\displaystyle I+ \displaystyle\frac 12
VG\displaystyle\right) \psi\displaystyle\right)_{\rho}.
\hfill (2.3.1)$
\vglue.3cm
\noindent
The equality (2.3.1) can be
intuitively
understood through the backward tree picture
(Gorostiza and Wakolbinger$^{(22)}$, Section 4, and
references therein, and Appendix, Section 3.A):
The intensity measure of the canonical Palm
distribution
$(R_\infty^1)_x$ (with {\it ego} at site $x$) is
\begin{eqnarray*}\label{1levelpalm}
\qquad\qquad\qquad\qquad
 \int  \langle \mu,  \psi \rangle (R_\infty^1)_x (d\mu)
&=&\psi (x) + E_x\left(\int_0^\infty
\int p_t(W_t,dy) \psi (y) V dt \right)
\nonumber\\
&=&\psi (x) + \frac 12 VG\psi (x),\qquad\qquad
\qquad\qquad\qquad\qquad\,\,  (2.3.2)
\end{eqnarray*}
where $p_t$ denotes the transition probability of the motion.
Then, by the Palm formula (A.2.1)
(Appendix), (2.3.1) follows by integrating
(2.3.2)
 with respect to
$\varphi(x)
\rho(dx)$.
The same reasoning shows that the space--time correlation structure of
$\underline
R_\infty^1$ is given by
$$
E_{\underline R_\infty^1}(\langle X_t, \varphi \rangle
\langle X_{t+s},\psi
\rangle)
= \left(\varphi, \left(I+ \frac 12 VG\right)T_s \psi\right)_{\rho}.
$$
This reveals how the normed second moment measure of the occupation
time behaves:
$$
\frac 1t E_{\underline R_\infty^1}\left(\displaystyle
\int_0^t\langle X_s, \varphi
\rangle ds \displaystyle\int_0^t
\langle X_r,\psi
\rangle dr \right)  \to
2\left(\varphi, \left(\displaystyle I+{V\over
2}G\right)G\psi\right)_{\rho}
  \mbox { as }  t
\to \infty.
$$
\item
Because of the obvious identities
$$\frac 12 (\varphi, G\varphi)_{\rho} = \int_0^\infty (\varphi,
T_{2r}\varphi)_\rho \, dr
= \int_0^\infty  \langle \rho, (T_r \varphi)^2 \rangle
dr,$$
transience of the motion
implies persistence of the $1$--level branching system
(Gorostiza and Wakolbin\-ger$^{(19)}$
 Corollary 2.2).
Therefore the existence of the measure $R^1_\infty$, which is the assumed 
intensity for the Poisson initial condition of the 2--level branching 
system, is
implied by  the strong transience assumption on the system.
\item
We now explain the second moment structure appearing in
Theorem 2.2.3(a). For a
level $1$ strongly transient motion it can be shown along the same
lines
as in Gorostiza et al$^{(17)}$
that the $2$--level branching particle system,
started off from
the Poisson system of $1$--level equilibrium clans,  is persistent.
Using e.g. the argument of
Gorostiza$^{(16)}$ (Lemma 4.6), one derives the following expression
 for
the second moment measure of the canonical measure $Q_\infty$ of the
aggregated
2--level equilibrium with intensity $\rho$
(using the fact that $R_1^\infty$ is the
Poissonization of the canonical measure of the superprocess
counterpart,  Appendix, (A.1.10)):
\vglue.3cm
\noindent
$\hfill
\displaystyle\int
\langle \mu, \varphi \rangle \langle \mu, \psi \rangle Q_\infty
(d\mu) =
\varphi,
\left(
I+ \displaystyle\frac 12
(V_1+V_2)
G +
\displaystyle\frac 14 V_1V_2G^2
\psi\right)_{\rho}.
\hfill (2.3.3)$
\vglue.3cm
\noindent
Then the second moment structure of the occupation time follows as in
the
1--level case.
Equality (2.3.3) can also be understood through
the backward tree picture: The
intensity measure of the Palm distribution $(Q_\infty)_x$ has the
representation (Hochberg and Wakolbinger$^{(26)}$)
\begin{eqnarray*}
\int  \langle \mu,  \psi \rangle (Q_\infty)_x (d\mu)
&=& \psi (x) + E_x\biggl[\biggr.
\int_0^\infty \biggl(\biggr.
\int (V_1+V_2) p_t(W_t,dy) \psi(y) \nonumber\\
&&+
\int_0^t \int V_1p_s(W_t,dz) \int V_2 p_{t-s}(z,dy) \psi(y) \,ds
\biggl.
\biggr) \,dt
\biggl.\biggr] \nonumber \\
&=& \psi (x) + \frac 12 (V_1+V_2)G\psi (x) + \frac 14 V_1V_2G^2\psi
(x).\qquad\qquad\qquad\qquad\qquad  (2.3.4)
\end{eqnarray*}
The summands on the r.h.s. of (2.3.4) have an interpretation in terms
of the genealogy:
$\displaystyle\frac 12 V_1G\psi(x)$ is the contribution of the
1--level relatives,
$\displaystyle\frac 12 V_2G\psi(x)$ is that of the 2--level relatives
breaking off directly
from the individual trunk, and
$\displaystyle\frac 14 V_1V_2G^2\psi(x)$ is that of the
2--level relatives having been generated by 1--level relatives, after
those have
broken off from the individual trunk.
Now (2.3.3) follows by integrating (2.3.4) with respect to
$\varphi(x)
\rho(dx)$.
\item
 The following table subsumes the covariance kernels appearing in the
second moment structures discussed above
(Theorems 2.2.2(a) and 2.2.3(a)). Columns 1, 2 and 3 refer to simple
motion, 1--level branching and 2--level branching, respectively.
$$
\begin{tabular}{l|l|l}
$2G$ \phantom{aa}&\phantom{aa}$2G + VG^2 $&
\phantom{aa}$2G + (V_1+V_2)G^2 + \displaystyle\frac 12
V_1V_2G^3$\\[.5cm]
&\phantom{aa}$=2\left(I+ \displaystyle{V\over 2}G\right)G$\phantom{aa}
&\phantom{aa}$=2\left(I+\displaystyle{V_2\over 2}G\right)\left(I+
\displaystyle{V_1 \over 2}G\right)G$
\end{tabular}
$$

\vglue.2cm
We observe a relationship between the  covariance
kernels for the $1$-- and $2$--level cases: For $V>0$ and an operator
$Q$,
we
define the operator
$$C_V(Q)=\left(I+{V\over 2}G\right)Q.$$
Then the $1$-- and $2$--level covariance kernels are given by
$2C_{V}(G)$ and
$2C_{V_2}(C_{V_1}(G))$, respectively. Thus, the $2$--level covariance
kernel is
like the $1$--level covariance kernel with $V$ replaced by  $V_2$ and
the
operator $G$ replaced by  $C_{V_1}(G)$. Recall that $G$
represents the expected total occupation time of the motion, and note
that
$C_V(G)$ represents the expected total occupation time of the mass
flow of the
motion plus a continuous throwing off of mass with intensity
$\displaystyle{1\over 2}V$, which also
evolves by the same flow.
One can then guess that for a 3--level system the covariance kernel
would be given by
 $$ 2C_{V_{3}}(C_{V_2}(C_{V_1}(G)))=2
\left(I+{V_3\over 2}G\right)\left(I+{V_2\over
2}G\right)\left(I+{V_1\over 2}G
\right)G,$$
and so on for higher levels of branching.
\vglue.5cm
\item
By using an argument of
Dawson and Perkins$^{(12)}$
(Section 5.4.4), one observes that a level $1$
transient motion leads to transient
equilibrium clans.
Indeed, even the expected value $a$  of the total future occupation
time in a bounded region $B$, starting from an equilibrium Palm
cluster, is
finite: Recall from (2.3.2) that the intensity measure  $\nu (dy)$ of
the Palm distribution $(R_\infty^1)_0$ of an equilibrium cluster
(with {\em
ego} at
the origin) is
$$\nu(dy)= \delta_0(dy) +
\frac 12 V G(0,dy),$$
hence $$a = \int\nu(dy) G(y,B) = G(0,B) + \frac 12
VG^2(0,B)
<
\infty.$$
\item
A level $2$
transient motion leads to transient aggregated $2$--level
equilibrium clans. Indeed, by (2.3.4) the intensity measure
$\sigma(dy)$ of the Palm distribution $(Q_\infty)_0$ is
$$
\sigma (dy)  = \delta_0(dy) +\left(
\frac 12 (V_1+V_2)G+ \frac 14 V_1V_2G^2\right)(0,dy),
$$
hence the expected value of the total future occupation
time in a bounded region $B$ is
$$\displaystyle{\int}\sigma(dy) G(y,B) = G(0,B) +
\left(\displaystyle{ \frac 12}
(V_1+V_2)G^2+\displaystyle{\frac 14} V_1V_2G^3\right) (0,B)
<\infty.$$
\item
For the $2$--level  system, if the intensity of the
initial Poisson distribution is $\delta_{\delta_{x}}\rho(dx)$
(instead
of $R_{\infty}^1$), then in the assumption for Theorem 2.2.3(b) the
growth of $G_{t}^2G$ is replaced by the growth
of $G_{t}^{3}$, and an analogous result holds.
\end{enumerate}
\vglue.5cm
\noindent
{\bf 2.4. Order of transience and recurrence,
asymptotics of powers of $ G_{t}$, and some special growth functions}
\vglue.5cm
\noindent
{\bf  Definition 2.4.1.} (a) Let ${\cal H}$ denote the class of 
differentiable functions $h:\erre_+\rightarrow\erre_+$ that are bounded and 
bounded away from $0$ on $[T,\infty)$ for some $T> 0$.
\vglue.5cm
\noindent
(b) For fixed $a\in (0,1)$, let
$$\widetilde{{\cal H}}_a=\{h\in{\cal H}\,\,
|\,\,h_t=h_{at}\quad\mbox{\rm for all}\quad t>0\}$$
(with $T=0$). Note that the elements of $\widetilde{{\cal H}}_a$ are 
periodic in a logarithmic scale.
\vglue.5cm
\noindent
{\bf  Definition 2.4.2. }
(a) For a given $\gamma
> 0$,
we say that $W_t$ is {\it transient
 of order }
$\gamma$
if
$R_t\sim t^{-\gamma}h_t J$ for some
$h\in{\cal H}$ as in Definition 2.4.1 and some
bilinear form $J$ as in Definition 2.1.3.

\noindent
(b) For a given $\gamma
\in (-1,0)$, we say that $W_{t}$ is {\it recurrent of order --}$\gamma$
if $G_{t}\sim  t^{-\gamma}h_t J$ for some $h$ and $J$ as above, and
we say that $W_{t}$ is {\it critically
recurrent} if
$G_{t}\sim  \log t \cdot J$ for some  $J$ as above.
\vglue.5cm
Clearly, for  transience of order $\gamma$ we have
\begin{eqnarray*}
\mbox{\rm level}\,\, k\,\,\mbox{\rm transience if and only if}&
\kern-.5cm& k<\gamma,\\
\mbox{\rm level}\,\, k\,\,\mbox{\rm recurrence if and only if}&
\kern-.5cm& k\geq \gamma,
\end{eqnarray*}
and therefore the following lemma holds.
\vglue.5cm
\noindent
{\bf Lemma 2.4.1.}  If $W_t$ is transient of order $\gamma$, then for each 
$k\geq 0, W_t$ is
\begin{eqnarray*}
 \mbox{\rm level}\,\, k\,\,\mbox{\rm strongly transient if and only if}&
\kern-.5cm& \gamma>k+1,\\
\mbox{\rm level}\,\,\, k\,\,\,\,\mbox{\rm weakly transient if and only if}&
\kern-.5cm&k<\gamma\leq k+1.
\end{eqnarray*}
\vglue.2cm
The next  lemma shows how
the assumptions of Theorems 2.2.2(b) and 2.2.3(b) on the growth of
powers of $G_t$ are implied by  transience of order $\gamma$ with $h\equiv 
1$.
\vglue.5cm
\noindent
{\bf  Lemma 2.4.2.} Let $W_t$ be transient of order $\gamma$
with
$R_t\sim t^{-\gamma}J$.

\noindent
(1) If $0<\gamma<1$, then
$$G_tG\sim t^{1-\gamma}H\quad \mbox{\rm with }\quad H={1\over
1-\gamma}J,$$
and
$$G^2_t\sim t^{1-\gamma}H\quad
\mbox{\rm with }\quad H={2-2^{1-\gamma}\over 1-\gamma}J.$$
(2) If $\gamma=1$, then
$$G_tG\sim G^2_t\sim \log t\cdot H \quad\mbox{\rm with }\quad H=J.$$
(3) If $1<\gamma<2$, then
$$G^2_tG\sim t^{2-\gamma}H\quad
\mbox{\rm with }\quad H={2-2^{2-\gamma}\over(2-\gamma)(\gamma-1)}J,$$
and
$$G^3_t\sim
t^{2-\gamma}H\quad \mbox{\rm with }\quad H=
{3^{2-\gamma}-2^{2-\gamma}-1\over (2-\gamma)(\gamma-1)}J. $$
(4) If $\gamma=2$, then
$$G^2_tG\sim G^3_t\sim\log t\cdot H\quad \mbox{\rm with }\quad H=J.$$
\vglue.5cm
{\rm  More generally, under the assumptions of Lemma 2.4.2 one can show 
that
$$G^\gamma_tG \sim G_t^{\gamma + 1} \sim \log t\cdot J \quad
\mbox{if} \quad \gamma \mbox { is an integer,}$$
and
$$G^{[\gamma]+1}_tG \sim t^{[\gamma]+1-\gamma}c_\gamma J,
\quad G^{[\gamma]+2}_t
 \sim t^{[\gamma]+1-\gamma}c_\gamma'J \quad \mbox{otherwise}$$
(for suitable constants $c_\gamma, c_\gamma'$).
But we will use only the cases $\gamma \le 2$ considered in the
lemma.}
\vglue.5cm
For growth
functions $f_t$ of the form
$f_t=t^\zeta h_t,0<\zeta<1, h\in\widetilde{{\cal H}}_a$,which will occur in 
one of the examples,
the normalizations $(\int^t_0f_sds)^{1/2}$ that appear in
the conclusions of
Theorems 2.2.1(b), 2.2.2(b) and 2.2.3(b) can be replaced by 
$(t^{\zeta+1}\widetilde{h}_t)^{1/2}$, where 
$\widetilde{h}\in\widetilde{{\cal H}}_a$, thanks to the following lemma.

\vskip.5cm

\noindent
{\bf Lemma 2.4.3.}  Let $h\in\widetilde{{\cal H}}_a$. If $\zeta>-1$, then
\vglue.3cm
\noindent
$\hfill
\displaystyle\int^t_1s^\zeta h_sds\sim t^{\zeta+1}
\widetilde{h}_t\,\,\,\,\mbox{\rm as}\, t\rightarrow\infty,\hfill (2.4.1)
$
\vglue.3cm
\noindent
where
$$\widetilde{h}_t=-\log a\cdot \int^\infty_0a^{r(1+\zeta)}h_{a^rt}dr.$$

Note that $\widetilde{h}\in \widetilde{{\cal H}}_a$.
\vglue.5cm
Lemma 2.4.3 applies to the growths of the fluctuations of the 0, 1 and 
2--level hierarchical system, but the l.h.s. of (2.4.1) can be computed
explicitly in this case.
\vglue.5cm
\noindent
{\bf 2.5. Infinite variance branching}
\vglue.5cm

In the case of infinite variance branching we
consider here only the so called
``$(1+\beta)$--branching'',
($0<\beta< 1$), whose offspring generating function
is of the form $s+q(1-s)^{1+\beta}, s\in[0,1]$,
 for some constant $q>0$ (this
law belongs to the domain of normal attraction of a stable law with 
exponent
$1+\beta$).
The picture now is less complete: we have only results
for the
 ``classical'' $t^{1/(1+\beta)}$--norming. In the 1--level case
it can be shown (along the lines of St\"ockl and Wakolbinger$^{(38)}$)
that this regime coincides with
that
of clan transience. We conjecture that an analogous result
 holds also  in the 2--level
case. Note that for $\beta=1$ we have the binary branching system 
considered above. However, we shall see that the results for the finite   
variance case are not special cases of the ones in this subsection.

We consider the following branching systems:

\begin{itemize}
\item[(i)]{\it {\rm 1}--level system}:
The system is as described in Subsection 2.2, except that the particles
undergo $(1+\beta)$--branching at rate $V$.
We assume that is is persistent. Hence
the system (with Poisson initial condition) converges to an
equilibrium  with intensity $\rho$. (Sufficient
conditions
 for this persistence are given in Gorostiza and Wakolbinger$^{(22)}$, 
Theorem 2.1). The equilibrium is then a
 Poisson superposition of ``2--level particles''.
Again we denote the  equilibrium canonical measure by $R^1_\infty$.
\item[(ii)]{\it {\rm 2}--level system}: The system is as described
in Subsection 2.2, except that
individual
particles undergo $(1+\beta_1)$--branching at rate $V_1$, clans undergo
$(1+\beta_2)$--branching at rate $V_2$, and it starts off from a
Poisson system of ``2--level particles'' with intensity $R^1_{\infty}$ (the
 equilibrium canonical measure for the 1--level $(1+\beta_1)$--branching 
system
with rate $V_1$).
\end{itemize}
\vglue.2cm
As above, $Y_t$ denotes the occupation time fluctuation (of the aggregated
system in the 2--level case).

\vskip.5cm

\noindent
{\bf Theorem 2.5.1.}  Let $X_t$ be the 1--level branching system.
 Assume that
\vglue.2cm
\noindent
$\hfill
\langle \rho, (G\varphi)^{1+\beta} \rangle < \infty, \quad \varphi
\in {\cal C}^+_c(S).
\hfill (2.5.1)$
\vglue.2cm
\noindent
 Then
$t^{-1/(1+\beta)}Y_t$ converges
to a random field $Z$ with
Laplace functional
$$
E\exp\{-\langle Z, \varphi \rangle\} = \exp\left\{\frac V{1+\beta}\langle 
\rho,
(G\varphi)^{1+\beta}\rangle\right\},\quad \varphi\in{\cal C}^+_c(S).
$$
\vglue.2cm
Here and in the next theorem
the notation $\langle Z,\varphi\rangle$ means the action of
the random field $Z$ on $\varphi$.
\vglue.5cm
Note that the finite variance
 case $\beta =1$ (Theorem 2.2.2 (a)) is not a special case of
Theorem 2.5.1, since this theorem would provide only the second term of the
covariance functional.
The term $2(\varphi,G\psi)_\rho$ in Theorem 2.2.2(a) does not appear in
Theorem
2.5.1 because the normalization is now
strong enough to kill this contribution to the occupation time fluctuations 
which comes from individuals related in direct line.
\vglue.5cm
\noindent
{\bf Theorem 2.5.2.}  Let $X_t$ be the 2--level branching system.
Assume that $\beta_2 <1$ and
\vglue.2cm
\noindent
$\hfill
\displaystyle\int \langle
\mu, G\varphi
\rangle^{1+\beta_2} R_\infty^1(d\mu) < \infty, \quad \varphi
\in {\cal C}^+_c(S).
\hfill (2.5.2)$
\vglue.2cm
\noindent
Then
$t^{-1/(1+\beta_2)}Y_t$ converges to a
random field
$Z$ with
Laplace functional
$$E\exp\{-\langle Z, \varphi \rangle\} = \exp\left\{\frac 
{V_2}{1+\beta_2}\int
\langle
\mu, G\varphi
\rangle^{1+\beta_2} R_\infty^1(d\mu)\right\},\quad \varphi\in
{\cal C}^+_c(S).$$
\vglue.5cm
 Note that also for the 2--level system the finite variance
case $\beta_1=\beta_2=1$ (Theorem 2.2.3(a)) is not a special case of 
Theorem
 2.5.2.

In the case of $\alpha$--stable motions in $\erre^d$, we shall see
in the next section
that conditions (2.5.1) and (2.5.2) hold for ``high'' dimensions.
\vglue.5cm
\noindent
{\bf 2.6. Occupation time fluctuations of superprocesses}
\vglue.5cm

As stated in the Introduction,
results analogous to the previous ones hold also for
the occupation time fluctuations of
the
1-- and 2-- level superprocesses corresponding to the
branching particle systems.
The results are simpler
because all the mass relationships closer than level $k$ do not contribute 
to the occupation time fluctuation limit of the $k$--level branching 
system.
Theorems 2.2.2(a) and 2.2.3(a) hold with only the terms involving $G^2$ and 
$G^3$ present
 in the
covariances of the limit random fields, respectively.
The proofs are similar to those for
the branching particle system and we shall omit them.
\vglue.5cm
\noindent
{\bf 3. TWO EXAMPLES OF INDIVIDUAL MOTIONS: SYMMETRIC
$\alpha$--STABLE PROCESS AND $c$--HIERARCHICAL RANDOM WALK}
\vglue.5cm
\noindent
{\bf 3.1. Transience/recurrence properties of the motions}
\vglue.5cm
In the first example the particle motion $W_t$ is the spherically
symmetric $\alpha$--stable L\'evy process in $S=\erre^d$,
$(0<\alpha\leq 2)$, and $\rho$ is the Lebesgue measure $\lambda$. In
the second example $W_t$
is a continuous--time random walk in the hierarchical group
$S=\Omega_N$ and $\rho$ is the counting measure $\nu$.  Due to
similarities in the
asymptotic behavior of $G_t$ and its powers for the two examples,
  which we will work
out in  this subsection, the limits of
 the occupation time fluctuations
of the corresponding 1-- and 2--level branching systems will also be
 analogous.

In this subsection we use several constants and functions whose definitions 
are collected in Section 4 for easy reference.

The analogy between the two motions is exhibited by a
constant $\gamma$ which we will define in each case.

For the $\alpha$--stable process we define
\vglue.3cm
\noindent
$\hfill
\gamma=\displaystyle{d\over \alpha}-1. \hfill (3.1.1)$
\vglue.3cm
\noindent

Before defining the $\gamma$ for the hierarchical random walk
 we will give some background.

The {\it hierarchical group} of order $N$ is a countable group
defined by
$$
\Omega_N=\{x=(x_1,x_2,\ldots)\,\,|\,\,x_i\in\{0,1,\ldots,N-1\}, x_i\neq 
0\,\,
\hbox{\rm  except for finitely many}\,\, i\},$$
with addition componentwise mod$(N)$. The {\it hierarchical distance}
$|\cdot|$ on $\Omega_N$ is defined by
$$
|x-y|=\max \{i|x_i\neq y_i\}.$$
A discrete--time {\it hierarchical random walk} in $\Omega_N$ jumps
from $x$ to $y$ such that $|x-y|=i\geq 1$ with probability
$r_i/N^{i-1}(N-1)$, where $r_1,r_2,\ldots$ is a probability
distribution on $\{1,2,\ldots\}$. This type of random walk was
introduced by Spitzer$^{(37)}$ for $N=2$ (the ``light bulb random walk''),
and by Sawyer and Felsenstein$^{(35)}$ for general $N$ in the context of
genetics models.

The continuous--time analogue of the hierarchical random walk with
jump rate $\sigma>0$ has transition density
$$
p_t(0,x)={1\over
N^i}(\delta_{0i}-1)e^{-\sigma(1-f_i)t}+(N-1)\sum^\infty_{j=i+1}{1\over
N^j}e^{-\sigma(1-f_j)t}$$
if $|x|=i\geq 0$, where $f_j$ is given by
$$f_j=r_1+\ldots +r_{j-1}-{r_j\over N-1},\quad j\geq 1,$$
(note that $f_0$ is irrelevant); see e.g. Fleischmann and Greven$^{(15)}$
for additional
information.

Here we will take $r_i$ of the form
$$r_i=\left({c\over N}\right)^{i-1}
\left(1-{c\over N}\right),\quad
i\geq 1,$$
where $c$ is a constant such that $0<c<N$. In this case we have
$f_j=1-a^jb,\quad
 j\geq 1,$ where
\vglue.3cm
\noindent
$\hfill
 a=\displaystyle{c\over N},\quad
 b=\displaystyle{N^2/c-1\over N-1}, \hfill (3.1.2)$
\vglue.3cm
\noindent
(note that $0<a<1, b>1$),  and the transition density becomes
\vglue.3cm
\noindent
$\hfill
p_t(0,x)=
\displaystyle{1\over N^i}(\delta_{0i}-1)e^{-\sigma
a^ibt}+(N-1)\displaystyle\sum\limits^\infty_{j=i+1}
\displaystyle{1\over N^j}e^{-\sigma a^jbt}\hfill (3.1.3)
$
\vglue.3cm
\noindent
if $|x|=i\geq 0$.

Since this random walk is characterized by the constant
$c$ (for fixed $N$), we will call it the {\it $c$--hierarchical
random walk}.
We shall see that $c$ is a mobility constant which plays
a similar role to that of $\alpha$ for the $\alpha$-stable process.

We define
\vglue.3cm
\noindent
$\hfill
\gamma=\displaystyle{\log c\over \log (N/c)},
\quad i.e., \quad c=N^{\gamma/(\gamma+1)},\hfill (3.1.4)$
\vglue.5cm
The following lemma shows that the constant $\gamma$ defined in (3.1.1)
for the $\alpha$--stable
process and in (3.1.4) for the $c$--hierarchical random walk
corresponds to the order of transience/recurrence parameter
(Definition 2.4.2) for the
respective processes, and it also shows the analogies for the semigroups 
$T_t$ and for the growth conditions for
 $G^2_t$ and $G^2_tG$ that appear as assumptions in Theorem 2.2.2(b) and 
Theorem 2.2.3(b) for the two motions.
Each one of the functions $h_t$ appearing in the lemma equals 1 for
 the $\alpha$-stable process, and belongs to $\widetilde{{\cal H}}_a$
for the  c-hierarchical random walk, with $a$ given by (3.1.2), or, 
equivalently, by (4.2.3).
\vglue.5cm
\noindent
{\bf Lemma 3.1.1.} Let $W_t$ be either the $\alpha$--stable process in 
$\erre^d$ or the $c$--hierarchical random walk in $\Omega_N$.

\begin{itemize}
\item[(a)] For all $\gamma>-1$,
$$(\varphi,T_t\psi)_\rho\sim t^{-(\gamma+1)}h_tq_\gamma 
\langle\rho,\varphi\rangle\langle\rho,\psi\rangle.$$
\item[ (b)] For $\gamma>0, W_t$ is transient of order $\gamma$ with
$$R_t(\varphi,\psi)\sim 
t^{-\gamma}h_tq_\gamma\langle\rho,\varphi\rangle\langle\rho,\psi\rangle.$$
\item[(c)] For $-1<\gamma<0,W_t$ is recurrent of order --$\gamma$ with
$$(\varphi,G_t\psi)_\rho\sim 
t^{-\gamma}h_tq_\gamma\langle\rho,\varphi\rangle\langle\rho,\psi\rangle.$$
\item[ (d)] For $\gamma=0, W_t$ is critically recurrent with
$$(\varphi,G_t\psi)_\rho\sim
\log t\cdot q_0\langle\rho,\varphi\rangle\langle\rho,\psi\rangle.$$
\item[ (e)] For $0<\gamma<1, W_t$ is weakly transient with
$$(\varphi,G^2_t\psi)_\rho\sim 
t^{1-\gamma}h_tq_\gamma\langle\rho,\varphi\rangle\langle\rho,\psi\rangle.$$
\item[ (f)] For $\gamma=1, W_t$ is weakly transient with
$$(\varphi,G^2_t\psi)_\rho\sim \log t\cdot 
q_1\langle\rho,\varphi\rangle\langle\rho,\psi\rangle.$$
\item[ (g)]  For $1<\gamma<2, W_t$ is level 1 weakly transient with
$$(\varphi,G^2_tG\psi)_\rho\sim 
t^{2-\gamma}h_tq_\gamma\langle\rho,\varphi\rangle\langle\rho,\psi\rangle.$$
\item[ (h)]  For $\gamma=2, W_t$ is level 1 weakly transient with
$$(\varphi,G^2_tG\psi)_\rho\sim \log t\cdot 
q_2\langle\rho,\varphi\rangle\langle\rho,\psi\rangle.$$
\end{itemize}
\vglue.2cm
The correspondences for the examples are:

 For the $\alpha$--stable process:
$\gamma=\displaystyle{d\over \alpha}-1, \rho=\lambda, h\equiv 1$ in all 
cases,
\vglue.2cm
\noindent
\begin{itemize}
\item[(a)]$q_\gamma=\kappa_{d,\gamma}$,
\item[ (b)]$\alpha<d, 
q_\gamma=\displaystyle{\kappa_{d,\gamma}\over\gamma}$,
\item[ (c)]$\alpha>d, 
q_\gamma=\displaystyle{\kappa_{d,\gamma}\over-\gamma}$,
\item[ (d)]$\alpha=d, q_0=\kappa_{d,0}$,
\item[ (e)]$\displaystyle{d\over 2}<\alpha<d,
q_\gamma=\displaystyle{2-2^{1-\gamma}\over(1-\gamma)\gamma}\kappa_{d,\gamma}$,
\item[(f)]$\alpha=\displaystyle{d\over 2}, q_1=\kappa_{d,1}$,
\item[ (g)]$\displaystyle{d\over 3}<\alpha<
\displaystyle{d\over 2}, q_\gamma=\displaystyle{2-2^{2-\gamma}\over
(2-\gamma)(\gamma-1)\gamma}\kappa_{d,\gamma}$,
\item[ (h)]$\alpha=\displaystyle{d\over 3}, q_2=\kappa_{d,2}$.
\end{itemize}
\vglue.2cm
For the $c$--hierarchical random walk:
$\gamma=\displaystyle{\log c\over\log (N/c)}, \rho=\nu$,
\begin{itemize}
\item[(a)]$q_\gamma=\kappa_{N,\gamma}, h=h^{(1,\gamma+1)}$,
\item[(b)]$c>1, q_\gamma=\kappa_{N,1}, h=h^{(1,\gamma)}$,
\item[(c)]$c<1, q_\gamma=\kappa_{N,\gamma}, h=h^{(2,\gamma)}$,
\item[(d)]$c=1, q_0=\displaystyle{\kappa_{N,0}\over \log N},$
\item[(e)]$1<c<N^{1/2},
q_\gamma=\kappa_{N,\gamma}, h=h^{(3,\gamma-1)}$,
\item[(f)]$c=N^{1/2}, q_1=\displaystyle{2\kappa_{N,1}\over\log N}$,
\item[(g)]$N^{1/2}<c<N^{2/3}, q_\gamma=\kappa_{N,\gamma}, 
h=h^{(3,\gamma-2)}$,
\item[(h)]$c=N^{2/3}, q_2=
\displaystyle{3\kappa_{N,2}\over\log N}$.
\end{itemize}
The constants $\kappa_{d,\alpha}$ and $\kappa_{N,\gamma}$, and the
functions $h^{(\cdot,\cdot)}$ are defined in Section 4: expressions 
(4.1.1), (4.2.1), (4.2.5), (4.2.6) and (4.2.7).
\vglue.5cm
\noindent
{\bf Corollary 3.1.1} (to Lemmas 3.1.1 and 2.4.1).

The $\alpha$--stable process is level $k$ strongly transient if and only if
$$
\alpha<{d\over k+2},$$
and level $k$ weakly transient if and only if
$$
{d\over k+2}\leq\alpha<{d\over k+1}.$$
\vglue.2cm
The $c$--hierarchical random walk is level $k$ strongly transient if and 
only if
$$
c>N^{(k+1)/(k+2)},$$
and level $k$ weakly transient if and only if
$$
N^{k/(k+1)}<c\leq N^{(k+1)/(k+2)}.
$$
\vglue.2cm
{\rm  Corollary 3.1.1
for the $\alpha$--stable process with $k=0$ is well known (Sato$^{(34)}$).}

\bigskip
{\rm The powers of the  Green potential operator of $W_t$ are given in
the next
lemma.}

\vglue.5cm
\noindent
{\bf Lemma 3.1.2.}  Let $\gamma>0$ and $1\leq j<\gamma+1$.

\begin{itemize}
\item[(a)]  For the
$\alpha$--stable process, the integral kernel of $G^j$ is
$G_{d,\gamma,j}(x)$ given by (4.1.2) with $\alpha j<d$.
\item[ (b)]For the $c$--hierarchical random walk, the
integral kernel of $G^j$ is $G_{N,\gamma,j}(x)$ given by (4.2.2)
with $c>N^{(j-1)/j}$.
\end{itemize}
\vglue.2cm
{\rm The next lemma shows that the condition (2.2.1)  in
Theorem 2.2.3(a) is fulfilled in both examples with $\delta = 3$}.
\vglue.5cm
\noindent
{\bf  Lemma 3.1.3.}
Let $W_t$ be either the $\alpha$-stable
process
in $\erre^d$ or the $c$--hierarchical random walk in $\Omega_N$. Then
level $k$ strong transience
implies that $||T_{t}\varphi|| = {\rm
o}(t^{-(k+2)})$,
$\varphi \in {\cal C}^{+}_c(\erre^d)$
(resp. $C^+_c(\Omega_N)$).
\vskip.5cm
We give next upper and lower bounds for the functions $h$
and $\widetilde{h}$  defined in Section 4, which appear in
Lemma 3.1.1 and in Theorem 3.2.1 below for the hierarchical case.
\vglue.5cm
\noindent
{\bf  Proposition 3.1.1.}
$$
\begin{tabular}{ccc}
\qquad\qquad{\rm Function}\qquad\quad&
\qquad\qquad{\rm Lower\,\, bound}\qquad\qquad\qquad
&\qquad\qquad{\rm Upper\,\, bound}\qquad\qquad\qquad\\[.5cm]
$h^{(1,\zeta)}_t,\zeta>0$&
$\displaystyle{\Gamma\over a^{-\zeta}-1}$
&$\displaystyle{a^{-\zeta}\Gamma\over a^{-\zeta}-1}$\\[.5cm]
$h^{(2,\zeta)}_t,-1<\zeta<0$
&$\displaystyle{a^{-\zeta}\Gamma\over 1-a^{-\zeta}}$
&$\displaystyle{\Gamma\over 1-a^{-\zeta}}$\\[.5cm]
$h^{(3,\zeta)}_t,-1<\zeta<0$
&$\displaystyle{a^{-\zeta}(2-2^{-\zeta)})\Gamma\over 1-a^{-\zeta}}$
&$\displaystyle{(2-2^{-\zeta})\Gamma\over 1-a^{-\zeta}}$\\[.5cm]
$\widetilde{h}^{(2,\zeta)}_t,-1<\zeta<0$
&$\displaystyle{a^{-\zeta}\Gamma\over (1-a^{-\zeta})(1-\zeta)}$
&$\displaystyle{\Gamma\over(1- a^{-\zeta})(1-\zeta)}$\\[.5cm]
$\widetilde{h}^{(3,\zeta)}_t,-1<\zeta<0$
&$\displaystyle{a^{-\zeta}(2-2^{-\zeta})\Gamma\over 
(1-a^{-\zeta})(1-\zeta)}$
&$\displaystyle{(2-2^{-\zeta})\Gamma\over (1-a^{-\zeta})(1-\zeta)}$
\end{tabular}
$$
where $\Gamma\equiv\Gamma(\zeta+1)$.
\vglue.5cm
\noindent
{\bf 3.2. Occupation time fluctuations limits}
\vglue.5cm
We give first the occupation
time fluctuation limits
for the two examples  in the case of finite variance
branching. In contrast with the general theorems (Theorems 2.2.1 -- 2.2.3)
we present the results
in a different order,
ending up with the classical
$t^{1/2}$--norming in each case. This is because the emphasis now, in the 
$\alpha$--stable case, is in going
 from the intermediate to the high dimensions $d$.

We denote by ${\cal N}$ a  real--valued centered Gaussian random variable
whose variance is specified in each case.

\vskip.5cm

\noindent
{\bf  Theorem 3.2.1.}  Let $W_t$ be either the $\alpha$--stable process
in $\erre^d$ or the $c$--hierarchical random walk in $\Omega_N$.
\vglue.3cm
\noindent
 0--level:

\begin{itemize}
\item[ (i)] For 
--$1<\gamma<0,(t^{1-\gamma}\widetilde{h}_t^{(2,\gamma)})^{-1/2}Y_t$
converges to ${\cal N}\rho$, where ${\cal N}$ has variance
$$2q^{(0)}_\gamma .$$
\item[ (ii)] For $\gamma=0, (t\log t)^{-1/2}Y_t$ converges
to ${\cal N}\rho$, where ${\cal N}$ has variance
$$2q^{(0)}_0.$$
\item[ (iii)] For $\gamma >0, t^{-1/2}Y_t$ converges to a
Gaussian field with covariance kernel
$$2Q_{\gamma,1}(x).$$
\end{itemize}
1--level:
\begin{itemize}
\item[ (i)]For $0<\gamma<1, (t^{2-\gamma }\widetilde{h}^{(3,\gamma 
-1)}_t)^{-1/2}Y_t$ converges
to ${\cal N}\rho$, where ${\cal N}$ has variance
$$V q^{(1)}_\gamma.$$
\item[(ii)] For $\gamma=1, (t\log t)^{-1/2}Y_t$ converges
to ${\cal N}\rho$, where ${\cal N}$ has variance
$$Vq^{(1)}_1.$$
\item[(iii)] For $\gamma >1, t^{-1/2}Y_t$ converges to a
Gaussian field with covariance kernel
$$2Q_{\gamma,1}(x)+VQ_{\gamma,2}(x).$$
\end{itemize}
2--level:
\begin{itemize}
\item[ (i)] For $1<\gamma<2, (t^{3-\gamma}\widetilde{h}^{(3,\gamma 
-2)}_t)^{-1/2}Y_t$
converges to ${\cal N}\rho$, where ${\cal N}$ has variance
$${V_1V_2\over 2} q^{(2)}_\gamma.$$
\item[ (ii)] For $\gamma=2, (t\log t)^{-1/2}Y_t$ converges
to ${\cal N}\rho$, where ${\cal N}$ has variance
$${V_1V_2\over 2}q^{(2)}_2.$$
\item[ (iii)]For $\gamma >2, t^{-1/2}Y_t$ converges to a
Gaussian field with covariance kernel
$$2Q_{\gamma,1}(x)+(V_1+V_2)Q_{\gamma,2}(x)+{V_1V_2\over 
2}Q_{\gamma,3}(x).$$
\end{itemize}
\vglue.5cm
The correspondences for the examples are as follows (recall that $a_t$ 
denotes the normalization):

For the
$\alpha$  stable process: $\gamma=\displaystyle{d\over \alpha}-1$, in this 
case the functions $\widetilde{h}^{(\cdot , \cdot)}$ are constant, and the
constants are included in the $q_\gamma$'s.
\vglue.5cm
\noindent
0--level:

\begin{itemize}
\item[ (i)]
$\alpha>d$,
$a_t=t^{(1-d/2\alpha)}, q^{(0)}_\gamma=
\displaystyle{\kappa_{d, \gamma}\over (1-\gamma)(-\gamma)}
=\displaystyle{\kappa'_{d, \alpha}\over (2-d/\alpha)(1-d/\alpha)}$.
\item[(ii)] $\alpha=d, a_t=(t\log t)^{1/2}, q^{(0)}_0=\kappa_{d,0}=
\kappa'_{d,d}$.
\item[(iii)]
$\alpha<d, a_t=t^{1/2},
Q_{\gamma,1}(x)=G_{d,\gamma,1}(x)=G'_{d,\alpha,1}(x)$.
\end{itemize}
1--level:

\begin{itemize}
\item[ (i)] $\displaystyle{d\over 2}<\alpha<d$,
 $a_t=t^{(3/2-d/2\alpha)},
q^{(1)}_\gamma=
\displaystyle{\kappa_{d,\gamma}(2-2^{1-\gamma})\over (2-\gamma )(1-\gamma 
)\gamma}=
\displaystyle{\kappa'_{d,\alpha}(2-2^{2-d/\alpha})\over
(3-d/\alpha)(2-d/\alpha)(d/ \alpha-1)}$.
\item[(ii)] $\alpha=\displaystyle{d\over 2}, a_t=(t\log t)^{1/2},
q^{(1)}_1=\kappa_{d,1}=\kappa'_{d,d/2}$.
\item[ (iii)] $\alpha<
\displaystyle{d\over 2}, a_t=t^{1/2},
Q_{\gamma,j}(x)=G_{d,\gamma,j}(x)=G'_{d,\alpha,j}(x), j=1,2$.
\end{itemize}
2--level:

\begin{itemize}
\item[(i)]
$\displaystyle{d\over 3}<\alpha<
\displaystyle{d\over 2}$,
$a_t=t^{(2-d/2\alpha)},
q^{(2)}_\gamma=
\displaystyle{\kappa_{d,\gamma}(2-2^{2-\gamma})\over (3-\gamma)(2-\gamma) 
(\gamma -1)\gamma}=
\displaystyle{\kappa'_{d,\alpha}(2-2^{3-d/\alpha})\over
(4-d/\alpha)(3-d/\alpha)(d/ \alpha -2)(d/ \alpha-1)}$.
\item[(ii)]
$\alpha=\displaystyle{d\over 3}, a_t=(t\log t)^{1/2},
q^{(2)}_2=\kappa_{d,2}=\kappa'_{d,d/3}$.
\item[{\it (iii)}] $\alpha<
\displaystyle{d\over 3},a_t=t^{1/2},
Q_{\gamma,j}(x)=G_{d,\gamma,j}(x)=G'_{d,\alpha,j}(x), j=1,2,3$.
\end{itemize}
\vglue.3cm
For the
$c$--hierarchical random walk:
$\gamma=\displaystyle{\log c\over \log (N/c)},
c=N^{\gamma/(\gamma+1)}$.
\vglue.3cm
\noindent
0--level:

\begin{itemize}
\item[(i)] $c<1$,
$a_t=t^{\log (N^{1/2}/c)/\log (N/c)}(\widetilde{h}^{(2,\gamma)}_t)^{1/2},
q^{(0)}_\gamma=\kappa_{N,\gamma}=\kappa'_{N,c}$.
\item[ (ii)] $c=1,a_t=(t\log t)^{1/2}, q^{(0)}_0=
\displaystyle{\kappa_{N,0}\over -\log a}= \displaystyle{(N-1)^2\over \sigma 
(N^2-1)\log N}$.
\item[ (iii)] $c>1,a_t=t^{1/2},
Q_{\gamma,1}(x)=G_{N,\gamma,1}(x)=G'_{N,c,1}(x)$.
\end{itemize}
1--level:

\begin{itemize}
\item[(i)] $1<c<N^{1/2}$,
$a_t=t^{\log (N/c^{3/2})/\log (N/c)}(\widetilde{h}^{(3,\gamma 
-1)}_t)^{1/2},
q^{(1)}_\gamma=\kappa_{N,\gamma}=\kappa'_{N,c}$.
\item[(ii)] $c=N^{1/2}, a_t=(t\log t)^{1/2},q^{(1)}_1=
 \displaystyle{\kappa_{N,1}\over \log a^{-1/2}}=\displaystyle{2(N-1)^3\over 
(\sigma (N^{3/2}-1))^2\log N}$.
\item[ (iii)] $c>N^{1/2},
Q_{\gamma,j}(x)=G_{N,\gamma,j}(x)=G'_{N,c,j}(x), j=1,2$.
\end{itemize}
2--level:

\begin{itemize}
\item[ (i)] $N^{1/2}<c<N^{2/3}$,
$a_t=t^{\log (N^{3/2}/c^{2})/\log (N/c)}(\widetilde{h}^{(3,\gamma 
-2)}_t)^{1/2},
q^{(2)}_\gamma=\kappa_{N,\gamma}=\kappa'_{N,c}$.
\item[ (ii)] $c=N^{2/3},a_t=(t\log t)^{1/2}, 
q^{(2)}_2=\displaystyle{\kappa_{N,2}\over \log 
a^{-1/3}}=\displaystyle{2(N-1)^4\over (\sigma (N^{4/3}-1))^3\log N}$.
\item[(iii)] $c > N^{2/3}, a_t=t^{1/2},Q_{\gamma,j}(x)=
G_{N,\gamma,j}(x)=G'_{N,c,j}(x), j=1,2,3$.
\end{itemize}

We turn now to the case of infinite variance branching with
 $\alpha$-stable motion in $\erre^{d}$. The question is when do the 
assumptions for Theorems 2.5.1 and 2.5.2 hold.
\vglue.5cm
\noindent
{\bf  Proposition 3.2.1.}
For the symmetric $\alpha$--stable motion in $\erre^d$, condition
(2.5.1) holds if and only if
$d >\alpha\left(1+\displaystyle{\frac1\beta}\right).$
\vglue.5cm
\noindent
{\bf  Proposition 3.2.2.}
For the
symmetric $\alpha$--stable motion in $\erre^d$,
condition
(2.5.2) holds if and only if $\beta_2 < \beta_1$ and $d >
\alpha\left(1+
{\displaystyle\frac1\beta_2}\left(1+\displaystyle{\frac 1\beta_1}
\right)\right).$
\vglue.5cm
\noindent
{\bf 3.3. Comments on the results}

\begin{enumerate}
\item Note that for the 2--level  branching Brownian
system $(\alpha=2)$ the normings $t^{3/4},(t\log t)^{1/2}$ and $t^{1/2}$ 
mentioned in the Introduction correspond to dimensions $d=5, d=6$ and $d> 
6$, respectively, i.e., two dimensions higher than for the $1$--level 
system. In the Brownian case  the critical dimensions for
the
1-- and 2--level branching systems are $d=4$ and $d=6$, corresponding
to
$\gamma=1$ and $\gamma=2$, respectively. For 1--level Brownian systems,
  occupation time large deviation results
have been obtained by Deuschel and Rosen$^{(13)}$ (and references therein).
\item For the $\alpha$--stable process in $\erre^d$ the order of
transience/recurrence  parameter $\gamma$ defined in (3.1.1) can take 
values
only in the interval $[(d-2)/2, \infty)$. On the other hand,
for the $c$-hierarchical
random walk in $\Omega_N$
 the possible values of the parameter $\gamma$ defined in (3.1.4)
range over the
whole interval $(-1, \infty)$, which means that in $\Omega_N$
 there is a rich structure of naturally ordered
random walks.
\item Note that for
both  the $\alpha$--stable motion and the $c$--hierarchical random walk, in 
all cases (i) and (ii) of Theorem 3.2.1 the occupation time
fluctuation limits in different regions of the space $S$ are perfectly 
correlated. An intuitive explanation for this might come from the recurrent 
visits of each $k$--level equilibrium clan, $k=0,1,2$,
to all bounded regions $B\subset S$.
\item Equating the parameters $\gamma$ for the $\alpha$--stable
process (3.1.1) and the $c$--hierarchical random walk (3.1.4) we
obtain
$$c=N^{1-\alpha/d}.$$
For this value of $c$, by Lemma 3.1.1 the $c$--hierarchical random
walk in $\Omega_N$ and the $\alpha$--stable process in $\erre^d$
have  the same order of  transience/recurrence. Consequently, by Theorems
2.2.1, 2.2.2 and 2.2.3
the asymptotics of the occupation time fluctuations
are analogous for the corresponding $k$--level branching systems,
$k=0,1,2$.
 The only differences are in the constants and the kernels of the
powers of the Green operators which appear in the fluctuation limits in 
Theorem 3.2.1.
The same observation holds for branching systems of ``$\alpha $--stable''
random walks on the lattice $\ze^{d}$. In passing we note that
$\alpha$--stable motions with $\alpha < 2$ do not have finite moments
of order $ \le \alpha$, but this plays no role in the  asymptotics of the
occupation times. The corresponding $c$--hierarchical random walks have
finite moments of all orders.

\item For the $c$--hierarchical random walk with
$c=c_N=\eta N^{k/(k+1)}, \,\,k\geq 0$,  $\eta >1$,
the powers of the Green potential operator take a simple form in the
limit $N \to \infty$ :
all the powers of order $1 \le j \le k$  vanish
as $N \to \infty$,  and
for the $(k+1)$--st power we observe from (4.2.2) and Lemma 3.1.2(b) that
$$
\lim_{N\rightarrow\infty}G_{N,c_{_{N}},k+1}(x)=
\frac{\eta^{k+1}}{\sigma^{k+1}(\eta^{k+1}-1)}
(\eta^{k+1})^{-|x|}.
$$
In particular, $\lim_{N\rightarrow\infty}G_{N,\eta,1}$ and
$\lim_{N\rightarrow\infty}G_{N,\eta N^{1/2},2}$ have the same
spatial asymptotics. This indicates that ``near $\eta =1$'' a
similar analysis as  was carried out by Dawson and Greven$^{(7)}$ for 1-  
-level branching
hierarchical random walks ($k=0$) might also be possible for
2--level  branching hierarchical random walks ($k=1$).

\item
The results for the $\alpha$--stable motion with the $t^{1/2}$--norming can 
be extended to test functions in ${\cal C}^+_\tau(\erre^d)$.
 For example, for the 2--level system
(with $d>3\alpha$)
we can take $\tau (x)=(1+|x|^2)^{-q}$ with $d/2< q< (d+\alpha)/2$ (Dawson 
and Gorostiza$^{(6)}$).
Moreover, the results can be extended to convergence of ${\cal
S}'(\erre^d)$--valued random fields, where ${\cal S}'(\erre^d)$ is
the space of tempered distributions on $\erre^d$, using an argument of
Iscoe$^{(27)}$.

\item Similarly to the previous comment, for the $c$--hierarchical random 
walk the results with the $t^{1/2}$--norming can be extended to test 
functions in
${\cal C}^+_\tau(\Omega_N)$ with an appropriate function
$\tau$. For example, for the 2--level system, $\tau$ should be a function 
in
$L^1(\Omega_N,\nu)$ such that the function
$x\mapsto\sum\limits_y\tau(y)(N^2/c^3)^{|x-y|}$ is bounded.
\end{enumerate}
\vglue.5cm
\noindent
{\bf 4. DEFINITIONS OF CONSTANTS AND FUNCTIONS FOR THE EXAMPLES}
\vglue.5cm
\noindent
{\bf 4.1.
Notation for $\alpha$--stable motion:}
\vglue.3cm
\noindent
$\hfill
\kappa_{d,\gamma}=
\displaystyle{1\over (2\pi)^d}
\displaystyle\int_{\errita^d}e^{-|x|^{d/(\gamma+1)}}
dx=\kappa'_{d,\alpha}=\displaystyle{1\over
(2\pi)^d}
\displaystyle\int_{\errita^d}e^{-|x|^{\alpha}}dx,
\hfill (4.1.1)
$\vglue.3cm
\noindent
$\hfill
G_{d,\gamma,j}(x)=C_{d,\gamma,j}|x|^{-d(1-j/(\gamma+1))}=
G'_{d,\alpha,j}(x)=C'_{d,\alpha,j}|x|^{-(d-j\alpha)},
\hfill (4.1.2)
$\vglue.3cm
\noindent
where
$$
C_{d,\gamma,j}=
{
\Gamma
\left(
\displaystyle{d(\gamma+1-j)\over 2(\gamma+1)}\right)
\over
2^{jd/(\gamma+1)}\pi^{d/2}
\Gamma\left(\displaystyle{dj\over 2(\gamma+1)}\right)}
=C'_{d,\alpha,j}={\Gamma\displaystyle{\left({d-j\alpha\over 2}\right)}
\over 2^{j\alpha}\pi^{d/2}\Gamma\displaystyle{\left({j\alpha\over
2}\right)}},
$$
and $j$ is a positive integer such that $j<\gamma+1$, i.e. $\alpha
j<d$.
\vglue.5cm
\noindent
{\bf 4.2. Notation for $c$--hierarchical random walk:}
\vglue.3cm
\noindent
$\hfill
\kappa_{N,\gamma}=
(N-1)^{\gamma+2}(\sigma
(N^{(\gamma+2)/(\gamma+1)}-1))^{-(\gamma+1)}\hfill (4.2.1)
$
\vglue.3cm
\noindent
$$
=\kappa'_{N,c}
=(N-1)^{\log(N^2/c)/\log(N/c)}(\sigma(N^2/c-1))^{-\log N/\log(N/c)},
$$
\vglue.3cm
\noindent
$\hfill
G_{N,\gamma,j}(x)=C_{N,\gamma,j}N^{-|x|(1-j/(\gamma+1))}\,\,=\,\,
G'_{N,c,j}(x)=C'_{N,c,j}\displaystyle\left({N^{j-1}\over c^j}\right)^{|x|},
\hfill (4.2.2)
$
\vglue.3cm
\noindent
where
\begin{eqnarray*}
\lefteqn{C_{N,\gamma,j}
=\left({N-1\over\sigma(N^{(\gamma+2)/(\gamma+1)}-1)}\right)^j
\left[
\delta_{0,|x|}-1+{(N-1)N^{j-1}\over N^{j\gamma/(\gamma+1)}-N^{j-1}}\right] 
}\\
&=&C'_{N,\gamma,j}=
\left(\displaystyle{N-1\over\sigma(N^2/c-1)}\right)^j
\left[
\delta_{0,|x|}-1+\displaystyle{(N-1)N^{j-1}\over
c^j-N^{j-1}}\right],
\end{eqnarray*}
and $j$ is a positive
integer such that $j<\gamma+1$, i.e. $c>N^{(j-1)/j}$.

Note that $a$ and $b$ defined in (3.1.2) are also expressed as
\vskip.3cm
\noindent
$\hfill
a=N^{-1/(\gamma+1 )},\quad
b=\displaystyle{N^{(\gamma+2)/(\gamma+1)}-1\over N-1}.
\hfill (4.2.3)$
\vglue.3cm
To simplify notation we write
\vskip.3cm
\noindent
$\hfill
\theta=\sigma \displaystyle{N^{(\gamma+2)/(\gamma+1)}-1\over N-1}.\hfill 
(4.2.4)$
\vskip.3cm
\noindent
For $\zeta>0$, let
\vglue.3cm
\noindent
$\hfill
h^{(1,\zeta)}_t=
\displaystyle\sum\limits_{j=-\infty}^\infty (\theta a^jt)^\zeta e^{-\theta 
a^jt},\qquad t>0,\hfill (4.2.5)
$
\vglue.3cm
\noindent
and for $-1<\zeta<0$, let
\vglue.3cm
\noindent
$\hfill
h^{(2,\zeta)}_t=
\displaystyle\sum\limits^\infty _{j=-\infty }(\theta a^jt)^\zeta
(1-e^{-\theta a^jt}),\qquad t>0,\hfill (4.2.6)
$
\vglue.3cm
\noindent
$\hfill
h^{(3,\zeta)}_t=
\displaystyle\sum\limits^\infty _{j=-\infty }
(\theta a^jt)^\zeta (1-e^{-\theta a^jt})^2,\qquad t>0,\hfill (4.2.7)
$
\vglue.3cm
\noindent
$\hfill
\widetilde{h}^{(2,\zeta)}_t=
\displaystyle\sum\limits^\infty _{j=-\infty}
(\theta a^jt)^{\zeta-1} (e^{-\theta a^jt}-1+\theta a^jt),\qquad t>0,
\hfill (4.2.8)
$
\vglue.3cm
\noindent
$\hfill
\widetilde{h}^{(3,\zeta)}_t=
\displaystyle\sum\limits^\infty _{j=-\infty}
(\theta a^jt)^{\zeta-1}
\left(2e^{-\theta a^jt}-
\displaystyle{1\over 2}e^{-2\theta a^jt}+
\theta a^jt-\displaystyle{3\over 2}\right),\qquad t>0, \hfill (4.2.9)
$
\vskip.5cm
\noindent
with $a$ and $\theta$ given by (4.2.3) and (4.2.4). Note that
all the functions defined in (4.2.5)--(4.2.9) belong to
$\widetilde{{\cal H}}_a$ for $a$ given by (4.2.3). The functions 
$\widetilde{h}^{(2, \zeta)}$ and $\widetilde{h}^{(3, \zeta)}$ correspond
(asymptotically)
to $h^{(2, \zeta)}$ and $h^{(3, \zeta)}$, respectively, by Lemma 2.4.3,
but in this case they are obtained by explicit calculation of
the l.h.s. of (2.4.1).
\par
\noindent
{\bf 5. PROOFS}
\setcounter{section}{5}
\setcounter{equation}{0}
\vglue.5cm
\noindent
{\bf 5.1. Asymptotics of the powers of $G_t$}
\vglue.5cmm
\noindent
{\it Proof of Lemma 2.4.2:}
\par
Note that it suffices to do the proofs with $\varphi=\psi$.
Fix $\varphi\in{\cal C}^+_c(S), \varphi\neq 0$, and denote
 $J=J(\varphi,\varphi)$ and
$R_t=R_t(\varphi,\varphi)$. By assumption,  $R_t\sim t^{-\gamma}J$.
\vglue.5cm
\noindent
(1)
\vglue-1.3cm
$$
(\varphi,G_tG\varphi)_\rho=\left(
\varphi,\int^t_0\int^{\infty}_0T_{s+u}\varphi
dsdu\right)_\rho=\int^t_0R_udu\sim J{1\over 1-\gamma}t^{1-\gamma},
$$
and
\begin{eqnarray*}
\lefteqn{(\varphi,(G_tG-
G^2_t)\varphi)_\rho=
\left(\varphi,G_t\int^{\infty}_tT_s\varphi ds\right)_\rho
=\left(\varphi,\int^t_0\int^{\infty}_tT_{s+u}\varphi dsdu\right)_\rho
=\int^t_0R_{s+t}ds}\\
\kern-.5cm&&\sim J\int^t_0(s+t)^{-\gamma}ds\,\,=\,\,
J{1\over 1-\gamma}((2t)^{1-\gamma}-t^{1-\gamma})\,\,
=\,\,{J\over 1-\gamma}(2^{1-\gamma}-1)t^{1-\gamma},
\end{eqnarray*}
hence
$$(\varphi,G^2_t\varphi)_\rho=(\varphi,(G^2_t-
G_tG)\varphi)_\rho+(\varphi,G_tG\varphi)_\rho\sim
{J\over 1-\gamma}(2-2^{1-\gamma})t^{1-
\gamma}.$$
\vglue.5cm
\noindent
$(2)\hfill
(\varphi,G_tG\varphi)_\rho=\displaystyle\int^t_0R_udu\sim J
\displaystyle\int^t_1u^{-1}du\sim c \log t,\hfill$%
$$(\varphi,(G_tG-G^2_t)\varphi)_\rho=\displaystyle\int^t_0R_{s+t}
ds\hfill$$%
$$\qquad\qquad\quad \sim J\displaystyle\int^t_0(s+t)^{-1}ds
\sim J(\log(2t)-\log t)=o(\log t),$$
hence the assertion for $G^2_t$ follows.
\vglue.5cm
\noindent
$(3)\hfill
(\varphi, G^2_tG\varphi)_\rho=\displaystyle\int^t_0
\displaystyle\int^t_0R_{s+u}drdu\hfill\qquad\qquad\qquad\qquad$%
$$\sim
J\displaystyle\int^t_1
\displaystyle\int^t_1 (s+u)^{-\gamma}dsdu \sim
\displaystyle{J\over (2-\gamma)(\gamma-1)}
(2-2^{2-\gamma})t^{2-\gamma},$$%
$$(\varphi,(G^2_tG-G^3_t)\varphi)_\rho=\displaystyle\int^t_0
\displaystyle\int^t_0 R_{s+u+t}duds\qquad\qquad\qquad\qquad$$%
$$\sim
J\displaystyle\int^t_0\displaystyle\int^t_0(s+u+t)^{-\gamma}duds
\sim \displaystyle{J\over (2-\gamma)(\gamma-1)}(-3^{2-\gamma}+2\cdot
 2^{2-\gamma}-1)t^{2-
\gamma},$$
\newpage
\noindent
hence the assertion for $G^3_t$ follows.
\vglue.2cm
\noindent
$(4)\hfill
(\varphi,G^2_tG\varphi)_\rho=\displaystyle\int^t_0\displaystyle\int^t_0
R_{s+u}drdu\hfill$%
$$\qquad \sim
J\displaystyle\int^t_1\displaystyle\int^t_0(s+u)^{-2}dsdu\sim c\log t,$$%
$$(\varphi,(G^2_tG-G^3_t)\varphi)_\rho=\displaystyle\int^t_0
\displaystyle\int^t_0 R_{s+u+t}duds$$%
$$\sim J
\displaystyle\int^t_0\displaystyle\int^t_0
(s+u+t)^{-2}dsdu=o(\log t),$$
hence the assertion for $G^3_t$ follows. $\hfill\Box$
\vglue.5cm
\noindent
{\it Proof of Lemma 2.4.3:}
$$\int^t_1s^\zeta h_sds=-\log a\cdot 
\int^\tau_0a^{-r(1+\zeta)}h_{a^{-r}}dr,$$
where $\tau=-\displaystyle{\log t\over \log a}$.
Hence, since $a<1$ and $h$ is bounded, we have
\begin{eqnarray*}
t^{-(1+\zeta)}\int^t_1s^\zeta h_s ds&=&-\log a\cdot 
\int^\tau_0a^{(\tau-r)(1+\zeta)}h_{a^{-r}}dr\\
&=&-\log a\cdot\int^\tau_0 a^{r(1+\zeta)}h_{a^{-(\tau-r)}}dr\\
&\sim&-\log a\cdot \int^\infty_0 a^{r(1+\zeta)}h_{a^rt}dr.
\end{eqnarray*}
\vskip-1cm
$\hfill\Box$
\vglue.5cm
\noindent
{\bf 5.2. Main results}

\par
We will not include the proofs for the $0$--level
system (Theorem 2.2.1)
 because they are simpler versions of the proofs for the branching systems.
The proofs for the 1-- and 2--level
systems follow the idea of the method employed
by Iscoe$^{(27)}$ for  (1--level) superprocesses in $\erre^d$.
The particle systems are somewhat harder to deal with than the
superprocesses, but the main point is
it has been necessary to modify the method
in order to deal with the new
technical difficulties that arise from  the second level branching.
We will use the modified approach
  also for the 1--level system.
\vglue.5cm
\noindent
{\it Proof of Theorems 2.2.2  and 2.5.1 (1--level branching system):}

In order to simplify notation we write
\vglue.3cm
\noindent
$\hfill
\varphi_t =F_t^{-1/(1+\beta)}\varphi\,\,\,\,\hbox{\rm for}\,\,
\varphi\in{\cal C}^+_c(S), \varphi\neq 0\,\,\hbox{\rm and}\,\,\beta\leq 
1,\,\,
\mbox{\rm where}\,\, F_t=\displaystyle\int^t_0f_sds,\hfill (5.2.1)
$

\vskip.3cm

\noindent
and $f_s$ is a growth function.
For Theorem 2.2.2(a) and Theorem 2.5.1, $f_t$ is interpreted as
$f_t\equiv 1$.
\par
 We have,
by  (A.1.1) and (A.1.2) (Appendix),
$$E\exp\left\{-
F_t^{-1/(1+\beta)}\biggl\langle \int^t_0X_sds,\varphi\biggr\rangle\right\}=
\exp\{-\langle\rho,u_{\varphi_t}(t)\rangle\},$$
where $u_{\varphi_t}(s,x)$ with values in $ [0,1]$
is the unique  solution of

\vskip.3cm

\noindent
$\hfill
u_{\varphi_t}(s)=-\displaystyle{V\over1+\beta}\displaystyle\int^s_0
T_{s-r}(u_{\varphi_t} (r)^{1+\beta})dr+
\displaystyle\int^s_0 T_{s-r}(\varphi_t(1-u_{\varphi_t} (r)))dr
\hfill (5.2.2)
$

\vskip.3cm

\noindent
Hence, by $T_t$--invariance of $\rho$ and
$E\langle\int^t_0X_sds,\varphi\rangle=t\langle\rho,\varphi\rangle$,
for the occupation time fluctuation $Y_t$ we have

\vskip.3cm

\noindent
$\hfill
E \exp\{-F_t^{-1/(1+\beta)}\langle Y_t,\varphi\rangle\}=
\exp\displaystyle\left\{{V\over 1+\beta}(I_1(t)+I_2(t))+I_3(t))
\displaystyle\right\},
\hfill (5.2.3)
$
\vskip.3cm
\noindent
where
\vskip.3cm
\noindent
$\hfill
I_1(t) = \displaystyle\int_0^t\langle \rho, 
w_{\varphi_t}(s)^{1+\beta}\rangle
ds,\hfill (5.2.4)
$\vskip.3cm
\noindent
$\hfill
I_2(t) = \displaystyle\int_0^t\langle \rho, u_{\varphi_t}(s)^{1+\beta} -
w_{\varphi_t}(s)^{1+\beta}\rangle ds, \hfill (5.2.5)
$\vskip.3cm
\noindent
$\hfill
I_3(t) = \displaystyle\int_0^t\langle
\rho, \varphi_t u_{\varphi_t}(s)\rangle ds,\hfill (5.2.6)
$\vskip.3cm
\noindent
with
\vskip.3cm
\noindent
$\hfill
w_\varphi (s)(x):= \displaystyle\int^s_0 T_r\varphi(x)dr=G_s\varphi(x),\,\,
x\in S, \varphi\in{\cal C}^+_c(S). \hfill (5.2.7)
$
\vskip.5cm
We will prove the following limits as $t\rightarrow\infty$:
\vskip.4cm
\noindent
For Theorem 2.2.2(a):

\vskip.3cm
\noindent
$\hfill
I_1(t) \to (\varphi,
G^2\varphi)_\rho . \hfill (5.2.8)
$
\vskip.3cm
\noindent
For Theorem 2.2.2(b):
\vskip.3cm
\noindent
$\hfill
I_1(t) \to  H(\varphi, \varphi). \hfill (5.2.9)
$
\vskip.3cm
\noindent
For Theorem 2.5.1:
\vskip.3cm
\noindent
$\hfill
I_1(t)\to \langle\rho,(G\varphi)^{1+\beta}\rangle.
\hfill (5.2.10)
$
\vskip.3cm
\noindent
For  $\beta\leq 1$:
\vskip.3cm
\noindent
$\hfill
I_2(t)\to 0. \hfill (5.2.11)
$
\vskip.3cm
\noindent
For $\beta <1$:
\vskip.3cm
\noindent
$\hfill
I_3(t)\to 0. \hfill (5.2.12)
$\vskip.3cm
\noindent
For Theorem 2.2.2(a):
\vskip.3cm
\noindent
$\hfill
I_3(t)  \to (\varphi, G\varphi)_\rho. \hfill (5.2.13)
$\vskip.3cm
\noindent
For Theorem 2.2.2(b):
\vskip.3cm
\noindent
$\hfill
I_3(t)  \to 0. \hfill (5.2.14)
$\vskip.3cm
\noindent
These limits  will yield the  conclusions of the theorems.

\vskip.2cm

\noindent
{\it Proof of (5.2.8) and (5.2.9):} From (5.2.4) and (5.2.7) we have
$$I_1(t)=F_t^{-1}\int^t_0\langle\rho,(G_s\varphi)^2\rangle ds.$$
By L'H\^opital's rule, for (5.2.8) we have
$$
I_1(t)=t^{-1}\int^t_0\langle\rho,(G_s\varphi)^2\rangle ds
\sim(\varphi,G^2_t\varphi)_\rho\rightarrow(\varphi,G^2\varphi)_\rho,
$$
and for (5.2.9),
$$
I_1(t)=F_t^{-1}
\int^t_0
\langle\rho,(G_s\varphi)^2\rangle ds
\sim{1\over f_t}(\varphi,G^2_t\varphi)_\rho\rightarrow
H(\varphi,\varphi).
$$
{\it Proof of (5.2.10):} The same as (5.2.8).
\vskip.5cm
\noindent
{\it Proof of (5.2.11):} We rewrite (5.2.5) as
$$-I_2(t)=\int^t_0\langle\rho,w_{\varphi_t}(s)^{1+\beta}
[1-(u_{\varphi_t}(s)/w_{\varphi_t}(s))^{1+\beta}]\rangle ds.
$$
We have from (5.2.2) and (5.2.7)
\vskip.3cm
\noindent
$\hfill
u_{\varphi_t}(s)-w_{\varphi_t}(s)=-\displaystyle{V\over 1+\beta}
\displaystyle \int^s_0 T_{s-r}
(u_{\varphi_t}(r)^{1+\beta})ds-
\displaystyle\int^s_0 T_{s-r}(\varphi_tu_{\varphi_t}(r))dr
\leq 0, \hfill (5.2.15)
$\vskip.3cm
\noindent
hence
$$
0\leq1-(u_{\varphi_t}(s)/w_{\varphi_t}(s))^{1+\beta}\leq 1.
$$
Therefore, since the convergence of $I_1 (t)$ implies uniform
integrability, it suffices to show that
$$\lim_{t\to\infty}{u_{\varphi_t}(s)\over w_{\varphi_t}(s)}=1\,\,
\,\,\hbox{\rm for all}\,\, s,\rho-a.e.$$
We have from (5.2.7) and (5.2.15)
$$
{u_{\varphi_t}(s)\over 
w_{\varphi_t}(s)}=1-(G_s\varphi)^{-1}(J_t(s)+K_t(s)),$$
where
$$
J_t(s)={V\over 1+\beta}F_t^{1/(1+\beta)}\int^s_0 T_{s-r}(
u_{\varphi_t}(r)^{1+\beta})dr \geq 0,$$
and
$$
K_t(s)=\int^s_0 T_{s-r}(\varphi u_{\varphi_t}(r))dr\geq 0.$$
By (5.2.15),
\begin{eqnarray*}
J_t(s)&\leq&
{V\over 1+\beta}F_t^{1/(1+\beta)}\int^s_0 T_{s-r}(w_{\varphi_t}
(r)^{1+\beta})dr\\
&=&{V\over 1+\beta}F_t^{-\beta/(1+\beta)}\int^s_0 T_{r-s}((G_r\varphi)
^{1+\beta})dr\to 0.
\end{eqnarray*}
Similarly, $K_t(s)\rightarrow 0$, so (5.2.11) is proved.
\vskip.3cm
\noindent
{\it Proof of (5.2.12):} We have from (5.2.6) and (5.2.15)
$$
I_3(t)\leq t^{-2/(1+\beta)}
\int^t_0(\varphi,G_s\varphi)_\rho ds
=t^{(\beta-1)/(1+\beta)}
t^{-1}\int^t_0 (\varphi,G_s\varphi)_\rho ds,
$$
and the result follows
since $(\varphi,G_s\varphi)_\rho\rightarrow(\varphi,G\varphi)_\rho<\infty$
 as $s\to\infty$.
\vskip.3cm
\noindent
{\it Proof of (5.2.13):} We rewrite (5.2.6) as
$$I_3(t)= t^{-1}\int^t_0\langle \rho,\varphi\widetilde{u}_{\varphi_t}(s)
\rangle ds,$$
where $\widetilde{u}_{\varphi_t}(s):= t^{1/2}u_{\varphi_t}(s)$. Since
$$
t^{-1}\int^t_0\langle\rho,\varphi G_s\varphi\rangle ds\to
(\varphi, G\varphi)_\rho,
$$
it suffices to prove that
$$A_t:=t^{-1}\int^t_0\langle\rho,\varphi\widetilde{u}_{\varphi_t}(s)\rangle 
ds
-t^{-1}\int^t_0\langle\rho,\varphi G_s\varphi\rangle ds\to 0
$$
We have, from (5.2.2)
$$
\widetilde{u}_{\varphi_t}(s)=-{V\over 2}t^{-1/2}\int^s_0 T_{s-r}
(\widetilde{u}_{\varphi_t}(r)^2)dr + G_s\varphi-t^{-1/2}\int^s_0
T_{s-r}(\varphi\widetilde{u}_{\varphi_t}(r))dr,
$$
hence
$$
|A_t|\leq t^{-3/2}\biggl({V\over 2}\int^t_0\int^s_0\langle\rho,
\varphi T_{s-r}(\widetilde{u}_{\varphi_t}(r)^2)\rangle drds\biggr.
+\int^t_0\int^s_0\langle\rho,\varphi T_{s-
r}(\varphi\widetilde{u}_{\varphi_t}(r))\rangle drds \biggl. \biggr),
$$
but, from (5.2.2), $\widetilde{u}_{\varphi_t}(r)\leq G_r\varphi$ (since 
$f_t\equiv 1$), so
\begin{eqnarray*}
|A_t|&\leq& t^{-3/2}\biggl({V\over 2}
\int^t_0\int^s_0\langle\rho,\varphi
T_{s-r}(G_r\varphi)^2\rangle drds
+\int^t_0\int^s_0\langle\rho,\varphi T_{s-r}(\varphi G_r\varphi)\rangle 
drds
\biggr)\nonumber
\end{eqnarray*}
$\hfill
= t^{-3/2}\displaystyle\biggl(\displaystyle{V\over 2}
\displaystyle\int^t_0\langle\rho,\varphi G_{t-
r}(G_r\varphi)^2\rangle dr
+\displaystyle\int^t_0\langle\rho,\varphi G_{t-r}(\varphi
G_r\varphi)\rangle dr\displaystyle\biggl. \displaystyle\biggr),\hfill 
(5.2.16)
$\vskip.3cm
\noindent
therefore
$$|A_t|\leq t^{-1/2}\biggl\langle\rho,{V\over 2}\varphi 
G(G\varphi)^2+\varphi
G(\varphi G\varphi)\biggr\rangle.$$
Now, by strong transience,
$$\langle\rho,\varphi G(G\varphi)^2\rangle\,\,\,\leq\,\,\,
||G\varphi||\,\,
||G^2\varphi||\langle\rho,\varphi\rangle <\infty,$$
and by transience,
$$\langle\rho,\varphi G(\varphi G\varphi)\rangle \,\,\leq\,\,\parallel\!\!
G\varphi\!\!\parallel^2
\langle\rho,\varphi\rangle<\infty,$$
hence the result follows.
\vskip.3cm
\noindent
{\it Proof of (5.2.14):}
We can follow the same argument used for (5.2.13)  replacing $\varphi$ by
$\varphi F^{-1/2}_t$.
Both terms on the r.h.s. of (5.2.16) can be shown to converge to $0$
 by L'H\^opital's rule and the assumptions.

It follows from (5.2.3)--(5.2.6) and the limits (5.2.8)--(5.2.14) that
\begin{eqnarray*}
&&E {\rm exp} \{-(F_t)^{-1/(1+\beta)}\langle Y_t,\varphi\rangle\}\\[.5cm]
&&\rightarrow\left\{\begin{array}{ll}
{\rm exp}
\left\{\displaystyle{V\over 2}(\varphi,G^2\varphi)_\rho+(\varphi, 
G\varphi)_\rho\right\}&\mbox{\rm for Theorem 2.2.2(a),}\\[.5cm]
{\rm exp}\left\{\displaystyle{V\over 2}H(\varphi,\varphi)\right\}
&\mbox{\rm for Theorem 2.2.2(b),}\\[.5cm]
{\rm exp}
\left\{\displaystyle{V\over 1+\beta}\langle \rho,(G\varphi)^{1+\beta}
\rangle\right\}&\mbox{\rm for Theorem 2.5.1}
\end{array}\right.
\end{eqnarray*}
as $t\rightarrow\infty$.

Finally, the convergence of the bilateral Laplace functional implies the 
weak convergence of $Y_t$ as $t\rightarrow\infty$ (Iscoe$^{(27)}$, pp. 
106--107 and 112). $\hfill\Box$
\vskip.5cm
\noindent
{\it Proof of Theorems 2.2.3 and 2.5.3 (2--level branching system):}
\vskip.5cm

We will follow the same steps for the proof of the 1--level case, but now 
some
of them are harder.
The  problem is that, while the test functions $\varphi\in{\cal C}^+_c(S)$
and the measure $\rho$ for the 1--level system are not so difficult
to work with, for the 2--level system the test
 functions $\mu\mapsto\langle\mu,\varphi\rangle$
and the measure $R^1_\infty$
(which now plays the role of $\rho$) raise new technical questions that are 
not easy
 to deal with,
in particular involving the third moments of $R^1_\infty$.
The background on the
2--level system in the Appendix should be consulted at this point.

We continue to use the notation $\varphi_t$ introduced in  (5.2.1),
but now with
$\beta=\beta_2$, and we put $f_t\equiv 1$ for
Theorem 2.2.3(a) and Theorem 2.5.2.

For $\varphi\in{\cal C}^+_c(S)$ we have,
by (A.1.13) and (A.1.14) (Appendix),
$$
E\exp\left\{-F_t^{-1/(1+\beta_2)}
\biggl\langle\int^t_0X_sds,\varphi
\biggr\rangle\right\}=
\exp\{-
\langle\!\langle R^1_\infty , {\mathbf u}_{\varphi_t}(t)
\rangle\!\rangle\},$$
where
\vskip.3cm
\noindent
$\hfill
{\mathbf u}_{\varphi_t}(s)=-
\displaystyle{V_2\over 1+\beta_2}\displaystyle\int^s_0U_{s-
r}({\mathbf u}_{\varphi_t} (r)^{1+\beta_2})dr+
\displaystyle\int^s_0 U_{s-
r}(\langle\cdot,\varphi_t\rangle(1-{\mathbf u}_{\varphi_t}(r)(\cdot)))dr.
\hfill (5.2.17)
$\vskip.3cm
\noindent
Hence, by $U_t$--invariance of $R^1_\infty$  and
$E\langle\int^t_0X_sds,\varphi\rangle=t\langle\rho,\varphi\rangle$,
\vskip.3cm
\noindent
$\hfill
E\exp\{- F_t^{-1/(1+\beta_2)}\langle Y_t, \varphi \rangle\}  =
\exp
\displaystyle\left\{\displaystyle\frac {V_2}{1+\beta_2}
(I_1(t)+I_2(t))+I_3(t))\displaystyle\right\}, \hfill (5.2.18)
$
\vskip.3cm
\noindent
where
\vskip.3cm
\noindent
$\hfill
I_1(t) =\displaystyle\int_0^t\langle\!\langle R_\infty^1,
{\mathbf w}_{\varphi_t}(s)^{1+\beta_2}\rangle\!\rangle ds, \hfill (5.2.19)
$\vskip.3cm
\noindent
$\hfill
I_2(t) =\displaystyle\int_0^t\langle\!\langle R_\infty^1,
{\mathbf u}_{\varphi_t}(s)^{1+\beta_2} - {\mathbf w}_{\varphi_t}
(s)^{1+\beta_2}\rangle\!\rangle ds, \hfill (5.2.20)
$\vskip.3cm
\noindent
$\hfill
I_3(t) = \displaystyle\int_0^t\langle\!\langle
R_\infty^1, \langle \cdot,\varphi_t\rangle
{\mathbf u}_{\varphi_t}(s)(\cdot)\rangle\!\rangle ds, \hfill (5.2.21)
$\vskip.3cm
\noindent
with, by (A.1.4),
\vskip.3cm
\noindent
$\hfill
{\mathbf w}_\varphi(s)(\mu) :=
\displaystyle\int^s_0U_r(\langle\cdot,\varphi\rangle)(\mu)dr=
\displaystyle\int_0^s \langle \mu , T_r \varphi \rangle
dr = \langle \mu , G_s \varphi \rangle,\quad \mu\in{\cal M}_\tau(S).
\hfill (5.2.22)
$\vskip.5cm
We will prove the following limits as $t\rightarrow\infty$:

\noindent
For Theorem 2.2.3(a):
\vskip.3cm
\noindent
$\hfill
I_1(t) \to (\varphi,
G^2\varphi)_\rho+{V_1\over 2}(\varphi, G^3\varphi)_\rho.\hfill (5.2.23)
$\vskip.3cm
\noindent
For Theorem 2.2.3(b):
\vskip.3cm
\noindent
$\hfill
I_1(t) \to
\displaystyle{V_1\over 2} H(\varphi, \varphi).\hfill (5.2.24)
$\vskip.3cm
\noindent
For Theorem 2.5.2:
\vskip.3cm
\noindent
$\hfill
I_1(t) \to \langle\!\langle R^1_\infty,
\langle\cdot,G\varphi\rangle^{1+\beta_2}\rangle\!\rangle.\hfill (5.2.25)
$\vskip.3cm
\noindent
For $\beta_2\leq 1$:
\vskip.3cm
\noindent
$\hfill
I_2(t)\to 0. \hfill (5.2.26)
$\vskip.3cm
\noindent
For $\beta_2<1$:
\vskip.3cm
\noindent
$\hfill
I_3(t)\to 0.\hfill (5.2.27)
$\vskip.3cm
\noindent
For Theorem 2.2.3(a):
\vskip.3cm
\noindent
$\hfill
I_3(t)\to (\varphi,G\varphi)_\rho+\displaystyle{V_1\over 2}(\varphi,
G^2\varphi)_\rho.\hfill (5.2.28)
$\vskip.3cm
\noindent
For Theorem 2.2.3(b):
\vskip.3cm
\noindent
$\hfill
I_3(t)\to 0.\hfill (5.2.29)
$\vskip.3cm
\noindent
{\it Proof of (5.2.23) and (5.2.24):} We have from (5.2.19) and (5.2.22)
$$I_1(t)=F_t^{-1}\int^t_0\langle\!\langle
R^1_\infty,\langle\cdot,G_s\varphi\rangle^2\rangle\!\rangle ds.$$
By L'H\^opital's rule and using (A.1.12) we have, for (5.2.23),
\begin{eqnarray*}
I_1(t)&=&t^{-1}\int^t_0
\langle\!\langle R^1_\infty,\langle\cdot, G_s\varphi\rangle^2
\rangle\!\rangle ds
\sim \langle\!\langle R^1_\infty,\langle\cdot, G_t\varphi\rangle^2
\rangle\!\rangle \\
&=&(\varphi,G^2_t\varphi)_\rho+{V_1\over 2}(\varphi,G^2_tG\varphi)_\rho
\rightarrow (\varphi,G^2\varphi)_\rho+{V_1\over 
2}(\varphi,G^3\varphi)_\rho,
\end{eqnarray*}
and for (5.2.24),
\begin{eqnarray*}
I_1(t)&=&F_t^{-1}\int^t_0\langle\!\langle R^1_\infty,
\langle\cdot,G_s\varphi\rangle^2\rangle\!\rangle ds\\
&\sim& {1\over f_t}\langle\!\langle
R^1_\infty,\langle\cdot,G_t\varphi\rangle^2\rangle\!\rangle\\
&=& {1\over f_t}\left(
(\varphi,G^2_t\varphi)_\rho+{V_1\over 
2}(\varphi,G^2_tG\varphi)_\rho\right)\\
&&\rightarrow {V_1\over 2}H(\varphi,\varphi).
\end{eqnarray*}
{\it Proof of (5.2.25):} The same as (5.2.23).
\vskip.5cm
\noindent
{\it Proof of (5.2.26):} We rewrite (5.2.20) as
\begin{eqnarray*}
-I_2(t)&=&\int^t_0\langle\!\langle R^1_\infty,
{\mathbf w}_{\varphi_t}(s)^{1+\beta_2}-
{\mathbf u}_{\varphi_t}(s)^{1+\beta_2}\rangle\!\rangle ds\\
&=&\int^t_0\langle\!\langle
R^1_\infty,{\mathbf w}_{\varphi_t}(s)^{1+\beta}[1-({\mathbf
u}_{\varphi_t}(s)/{\mathbf w}_{\varphi_t}(s))^{1+\beta_2}]
\rangle\!\rangle ds.
\end{eqnarray*}
Since $0\leq 1-({\mathbf u}_{\varphi_t}(s)/{\mathbf
w}_{\varphi_t}(s))^{1+\beta_2}\leq 1$ by
(5.2.2) and  (A.1.16), and since $I_1(t)$
converges,
 it suffices to prove that
$$\lim_{t\rightarrow\infty}{{\mathbf u}_{\varphi_t}(s)\over{\mathbf
w}_{\varphi_t}(s)}=1\,\,\mbox{\rm for all}\,\, s,R^1_\infty-a.e.$$
We have from (5.2.22) and  (A.1.16)
$${{\mathbf u}_{\varphi_t}(s)(\mu)\over
{\mathbf w}_{\varphi_t}(s)
(\mu)}=1-\langle\mu,G_s\varphi\rangle^{-1}(J_t(s)+K_t(s))(\mu)\geq 0,$$
where
$$J_t(s)(\mu)=
{V_2\over 1+\beta_2}F_t^{1/(1+\beta_2)}\int^s_0U_{s-r}({\mathbf
u}_{\varphi_t}(r)^{1+\beta_2})(\mu)dr\geq 0$$
and
$$K_t(s)(\mu)=\int^s_0U_{s-r}(\langle\cdot,\varphi\rangle{\mathbf
u}_{\varphi_t}(r)(\cdot))(\mu)dr\geq 0.$$
By (A.1.17) and (5.2.22)
\begin{eqnarray*}
J_t(s)(\mu)&\leq&
{V_2\over 1+\beta_2}F_t^{1/(1+\beta_2)}
\int^s_0U_{s-r}({\mathbf w}_{\varphi_t}(r)^{1+\beta_2})(\mu)dr\\
&=&{V_2\over
1+\beta_2}F_t^{-\beta_2/(1+\beta_2)}\int^s_0U_{s-r}(\langle\cdot,G_r\varphi
\rangle^{1+\beta_2})(\mu)dr\rightarrow 0.
\end{eqnarray*}
Similarly, $K_t(s)(\mu)\rightarrow 0$, so the result follows.
\vglue.3cm
\noindent
{\it Proof of (5.2.27):} We have from (5.2.21), (5.2.22) and (A.1.12)
\begin{eqnarray*}
\lefteqn{I_3(t)\leq
t^{-2/(1+\beta_2)}\int^t_0
\langle\!\langle
R^1_\infty,\langle\cdot,\varphi\rangle \langle\cdot,G_s
\varphi\rangle\rangle\!\rangle
ds}\\
&&=t^{(\beta_2-1)/(1+\beta_2)}t^{-1}\int^t_0
\left(
(\varphi,G_s\varphi)_\rho+{V_1\over 2}(\varphi,G_sG\varphi)_\rho)\right)ds.
\end{eqnarray*}
Since
$$(\varphi,G_s\varphi)_\rho+{V_1\over 2}(\varphi,G_sG\varphi)\rightarrow
(\varphi,G\varphi)_\rho+{V_1\over 2}(\varphi,G^2\varphi)_\rho
\quad\mbox{\rm as}\,\,s\rightarrow\infty,$$
the result follows.
\vglue.5cm
\noindent
{\it Proof of (5.2.28):} We rewrite (5.2.21) as
$$I_3(t)=t^{-1}\int^t_0\langle\!\langle
R^1_\infty,\langle\cdot,\varphi\rangle
\widetilde{{\mathbf u}}_{\varphi_t}(s)(\cdot)\rangle\!\rangle ds,$$
where $\widetilde{{\mathbf u}}_{\varphi_t}(s)=t^{1/2}{\mathbf
u}_{\varphi_t}(s)$. Since
$$
t^{-1}\int^t_0\langle\!\langle R^1_\infty,\langle\cdot,\varphi\rangle
\langle\cdot, G_s\varphi\rangle\rangle\!\rangle
\,\,\rightarrow\,\, (\varphi,G\varphi)_\rho+{V_1\over
2}(\varphi,G^2\varphi)_\rho,$$
by (A.1.12), it suffices to prove that
$$A_t:=t^{-1}\int^t_0\langle\!\langle R^1_\infty
,\langle\cdot,\varphi\rangle\widetilde{{\mathbf
u}}_{\varphi_t}(s)(\cdot)\rangle\!\rangle
ds-t^{-1}\int^t_0\langle\!\langle
R^1_\infty,\langle\cdot,\varphi\rangle \langle\cdot,G_s\varphi
\rangle\rangle\!\rangle ds\rightarrow 0.$$
\vglue.2cm
We have from (A.1.16)
$$
\widetilde{{\mathbf u}}_{\varphi_t}(s)(\mu)=-{V_2\over 2}t^{-1/2}\int^s_0
U_{s-r}(\widetilde{{\mathbf u}}_{\varphi_t}(r)^2)(\mu)dr
+\langle\mu,G_s\varphi\rangle-t^{-1/2}\int^s_0U_{s-r}
(\langle\cdot,\varphi\rangle\widetilde{{\mathbf 
u}}_{\varphi_t}(r)(\cdot))(\mu)dr,
$$
hence
\begin{eqnarray*}
|A_t|&\leq&t^{-3/2}\left(
{V_2\over 2}\int^t_0\int^s_0\int
R^1_\infty(d\mu)\langle\mu,\varphi\rangle U_{s-r}(\widetilde{{\mathbf
u}}_{\varphi_t}(r)^2)(\mu) drds\right.\\
&&\qquad\quad \left.+\int^t_0\int^s_0\int R^1_\infty(d\mu)
\langle\mu,\varphi\rangle U_{s-r}(
\langle\cdot,\varphi\rangle\widetilde{{\mathbf
u}}_{\varphi_t}(r))(\mu) drds\right).
\end{eqnarray*}
Now
$\widetilde{{\mathbf
u}}_{\varphi_t}(r)(\mu)\leq\langle\mu,G_r\varphi\rangle$ (since $f_t\equiv 
1$).
Note that this estimate is not too rough because $\widetilde{{\mathbf
u}}_{\varphi_t}(r)(\mu)\rightarrow\langle\mu,G_r\varphi\rangle$ as 
$t\rightarrow\infty$.
Hence
$$|A_t|\leq const.(H_1(t)+H_2(t)),$$
where
\begin{eqnarray*}
H_1(t)&=&t^{-3/2}\int^t_0\int^s_0\int
R^1_\infty(d\mu)\langle\mu,\varphi\rangle
U_{s-r}(\langle\cdot,G_r\varphi\rangle^2)(\mu) drds,\\
H_2(t)&=&t^{-3/2}\int^t_0\int^s_0\int
R^1_\infty(d\mu)\langle\mu,\varphi\rangle
U_{s-r}(\langle\cdot,\varphi\rangle\langle\cdot,G_r\varphi\rangle)
(\mu)
drds.
\end{eqnarray*}
We will show that $H_1(t)\rightarrow 0$.
The proof that $H_2(t)\rightarrow 0$ is similar.

Using (A.1.8) we have
$$H_1(t)\leq \, const. \sum^3_{j=1}J_j(t),$$
where
\begin{eqnarray*}
J_1(t)&=&t^{-3/2}\int^t_0\int^s_0\int
R^1_\infty(d\mu)\langle\mu,\varphi\rangle\langle\mu,T_{s-r}G_r
\varphi\rangle^2drds,\\
J_2(t)&=&t^{-3/2}\int^t_0\int^s_0\int
R^1_\infty(d\mu)\langle\mu,\varphi\rangle\langle\mu,T_{s-r}(G_r\varphi)^2
\rangle drds,\\
J_3(t)&=&t^{-3/2}\int^t_0\int^s_0\int
R^1_\infty(\mu)\langle\mu,\varphi\rangle
\int^{s-r}_0\langle\mu,T_{u}(T_{s-r-u}G_r\varphi)^2\rangle
du drds.
\end{eqnarray*}
By (A.1.12) and (A.1.13) we obtain
$$J_1(t)\leq const. \sum^5_{j=1}K_{1,j}(t),\,\,\,\,\,J_2
(t)\leq const. \sum^2_{j=1}K_{2,j}(t),\,\,\,\,\,
J_3(t)\leq const. \sum^2_{j=1}K_{3,j}(t),$$
where
\begin{eqnarray*}
K_{1,1}(t)&=&t^{-3/2}\int^t_0\int^s_0\langle\rho,\varphi(T_{s-r}
G_r\varphi)^2\rangle
drds,\\
K_{1,2}(t)&=&t^{-3/2}\int^t_0\int^s_0\langle\rho,\varphi T_{s-r}
G_r\varphi\cdot T_{s-r}GG_r\varphi\rangle drds,\\
K_{1,3}(t)&=&t^{-3/2}\int^t_0\int^s_0\langle\rho,\varphi
G(T_{s-r}G_r\varphi)^2\rangle drds,\\
K_{1,4}(t)&=&t^{-3/2}\int^t_0\int^s_0\int^\infty_0\left\langle\rho,\varphi 
G
 T_u(T_uT_{s-r}G_r\varphi)^2 \right\rangle dudrds,\\
K_{1,5}(t)&=&t^{-3/2}\int^t_0\int^s_0\int^\infty_0\langle\rho,
T_{s-r}G_r\varphi\cdot GT_u(T_u\varphi\cdot T_uT_{s-r}G_r\varphi)\rangle
dudrds,\\
K_{2,1}(t)&=&t^{-3/2}\int^t_0\int^s_0\langle\rho,\varphi
T_{s-r}(G_r\varphi)^2\rangle drds,\\
K_{2,2}(t)&=&t^{-3/2}\int^t_0\int^s_0\langle\rho,\varphi
GT_{s-r}(G_r\varphi)^2\rangle drds,\\
K_{3,1}(t)&=&t^{-3/2}\int^t_0\int^s_0\int^{s-r}_0\left\langle\rho,\varphi
T_u(T_{s-r-u}G_r\varphi)^2\right\rangle dudrds,\\
K_{3,2}(t)&=&t^{-3/2}\int^t_0\int^s_0\int^{s-r}_0\langle\rho,\varphi
G T_u(T_{s-r-u}G_r\varphi)^2\rangle dudrds.
\end{eqnarray*}
We will show that each of these terms converges to $0$ as 
$t\rightarrow\infty$. Recall that
$||G^j\varphi||<\infty, j=1,2,3$, for $\varphi\in{\cal C}^+_c(S)$.
$$K_{1,1}(t)\leq||G\varphi|| t^{-3/2}\int^t_0
\langle\rho,\varphi G_sG\varphi\rangle ds\leq 
||G\varphi||\,\,||G^2\varphi||\langle\rho,\varphi\rangle 
t^{-1/2}\rightarrow 0.$$%
$$K_{1,2}(t)\leq||G^2\varphi|| t^{-3/2}\int^t_0
\langle\rho,\varphi G_sG\varphi\rangle ds\leq ||G^2\varphi||^2
\langle\rho,\varphi\rangle t^{-1/2}\rightarrow 0.$$%
$$
K_{1,3}(t)\leq||G\varphi|| t^{-3/2}\int^t_0
\langle\rho,\varphi G_sG^2\varphi\rangle ds\leq||G\varphi||\,\,
||G^3\varphi||\langle\rho,\varphi\rangle t^{-1/2}\rightarrow 0.$$
\begin{eqnarray*}
\kern2cm K_{1,4}(t)&=&t^{-3/2}\int^t_0\int^s_0\int^\infty_0\langle \rho, 
\varphi GT_u (T_u T_rG_{s-r}\varphi )^2\rangle dudrds\\
&\sim& const.\;t^{-1/2}\int^t_0\int^\infty_0 \langle \rho, \varphi GT_u 
(T_{u+r}G_{t-r}\varphi)^2\rangle dudr
\mkern 120mu
\hbox{\rm (by l'H\^opital)}\\
&\sim&const.\;t^{1/2}\int^t_0\int^\infty_0\langle \rho, \varphi GT_u 
(T_{u+r}G_{t-r}\varphi \cdot T_{u+r}T_{t-r}\varphi )\rangle dudr
\quad\quad\,
\hbox{\rm (by l'H\^opital)}\\
&\leq&const.\;t^{1/2-\delta}\int^t_0 \int^\infty_0\langle \rho, \varphi 
GT_u (T_{u+r}G_{t-r}\varphi \cdot (t+u)^\delta T_{t+u}\varphi)\rangle 
dudr\\
&\leq& const. \;t^{3/2 -\delta}||G^3\varphi ||\langle \rho,
\varphi\rangle
\mkern 334mu
\hbox{\rm (by (2.2.1))}\\
&&\longrightarrow 0.
\end{eqnarray*}
\begin{eqnarray*}
\kern2cm
K_{1,5}(t)&=&t^{-3/2}\int^t_0 \int^s_0 \int^\infty_0 \langle \rho, 
T_{s-r}G_r \varphi \cdot T_u (T_u \varphi \cdot T_u T_r G_{s-r}\varphi 
)\rangle dudrds\\
&\leq& t^{-3/2}\int^t_0 \int^s_0 \int^\infty_0 \langle \rho, \varphi G^2 
T_u (T_u \varphi \cdot T_{u+r}G_{s-r}\varphi )\rangle dudrds\\
&\sim&const.\;t^{-1/2}\int^t_0\int^\infty_0 \langle \rho, \varphi G^2 T_u 
(T_u \varphi \cdot T_{u+r}G_{t-r}\varphi )\rangle dudr\qquad\qquad
\hbox{\rm (by l'H\^opital)}\\
&\sim& const.\;t^{1/2}\int^t_0 \int^\infty_0\langle \rho, \varphi G^2 T_u 
(T_u \varphi \cdot T_{u+r}T_{t-r}\varphi )\rangle dudr
\qquad\qquad\quad\hbox{\rm (by l'H\^opital)}\\
&\leq& const.\;t^{1/2-\delta}\int^t_0 \int^\infty_0 \langle \rho, \varphi 
G^2 T_u (T_u \varphi \cdot (t+u)^\delta T_{t+u}\varphi )\rangle dudr\\
&\leq& const.\;t^{3/2 -\delta}||G^3 \varphi || \langle \rho, \varphi 
\rangle\qquad
\mkern 280mu
 \hbox{\rm (by (2.2.1))}\\
&&\longrightarrow 0.
\end{eqnarray*}
$K_{2,1}(t)\rightarrow 0$, similarly to $K_{1,1}(t)\rightarrow 0.$

\noindent
$K_{2,2}(t)\rightarrow 0$, similarly to $K_{1,3}(t)\rightarrow 0$.
\begin{eqnarray*}
K_{3,1}(t)&\leq& ||G\varphi|| t^{-3/2}
\displaystyle\int^t_0\displaystyle\int^s_0
\langle\rho,\varphi rT_rG\varphi\rangle drds\\
&\leq& const. ||G\varphi||\,\,||G^3\varphi||\langle\rho,\varphi\rangle 
t^{-1/2}\rightarrow 0.
\end{eqnarray*}
\begin{eqnarray*}
\kern3cm
K_{3,2}(t)&=&t^{-2/3}\int^t_0 \int^s_0 \int^r_0 \langle \rho, \varphi GT_u 
(T_{r-u}G_{s-r}\varphi )^2 \rangle dudrds\\
&\sim& const.\;t^{-1/2}\int^t_0 \int^r_0 \langle \rho, \varphi G T_u 
(T_{r-u}G_{t-r}\varphi )^2\rangle dudr\qquad\quad\quad
\hbox{\rm (by l'H\^opital)}\\
&\sim&const.\;t^{1/2}\int^t_0\int^r_0\langle \rho, \varphi G T_u 
(T_{r-u}G_{t-r}\varphi \cdot T_{r-u}\varphi )\rangle dudr
\quad\,\,\,
\hbox{\rm (by l'H\^opital)}\\
&=&const.\;t^{1/2}\int^t_0\int^r_0\langle \rho, \varphi GT_{r-u}(T_u G_{  
t-r}\varphi \cdot T_u T_{t-r}\varphi )\rangle dudr\\
&=&const.\;t^{1/2}\int^r_0\int^{t-r}_0 \langle \rho, \varphi 
GT_{t-r-u}(T_uG_r \varphi \cdot T_{u+r}\varphi )\rangle dudr\\
&\sim& const.\;\;(M_1 H)+M_2(t)),
\end{eqnarray*}
where (by l'H\^opital),
$$
M_1(t)= t^{3/2}\int^t_0 \langle \rho, \varphi G(T_{t-r}G_r \varphi \cdot 
T_t\varphi )\rangle dr
$$
and, since $\displaystyle{d\over dt}GT_t\varphi =-T_t$,
$$
M_2 (t)=-t^{3/2}\int^t_0 \int^{t-r}_0\langle \rho, \varphi T_{t-r-u}(T_u 
G_r \varphi \cdot T_{u+r}\varphi )\rangle dudr.
$$
Now,
\begin{eqnarray*}
\kern3cm
M_1(t)&=&t^{3/2-\delta }\int^t_0 \langle \rho, \varphi G(T_{t-r}G_r \varphi 
\cdot t^\delta T_t \varphi \rangle dr\\
&\leq &const.\;t^{5/2 -\delta}||G^2\varphi ||\langle \rho, \varphi\rangle 
\qquad
\mkern 250mu
\hbox{\rm (by (2.2.1))}\\
&& \longrightarrow 0,
\end{eqnarray*}
\begin{eqnarray*}
\kern3cm
|M_2 (t)|&\leq & t^{3/2}\int^t_0\int^r_0 \langle \rho, \varphi T_{r-u}(T_u 
G_{t-r}\varphi \cdot T_{t-(r-u)}\varphi \rangle dudr\\
&=& t^{3/2}\int^t_0 \int^r_0 \langle \rho, \varphi T_u 
(T_{r-u}G_{t-r}\varphi \cdot T_{t-u}\varphi )\rangle dudr\\
&=& t^{3/2}\int^t_0 \left\langle \rho, \varphi T_u \left(\int^t_u 
T_{r-u}G_{t-r}
\varphi dr \cdot T_{t-u}\varphi\right)\right\rangle du\\
&\leq & ||G^2 \varphi ||t^{5/2}\langle \rho, \varphi T_t \varphi\rangle \\
&=& ||G^2 \varphi ||t^{5/2 -\delta}\langle \rho, \varphi t^\delta T_t 
\varphi \rangle\\
&\leq & const. \;||G^2\varphi ||\langle \rho, \varphi\rangle t^{5/2 
-\delta}
\mkern 280mu
\hbox{\rm (by (2.2.1))}\\
&&\longrightarrow 0.
\end{eqnarray*}
\vglue.2cm
Finally, the weak convergence of $Y_t$ as $t\rightarrow \infty $ follows
 as in  the $1$--level case. \hfill $\Box$
\vglue.5cm
\noindent
{\bf 5.3. Examples}
\vglue.5cm
\noindent
{\it Proof of Lemma 3.1.1:}
\vglue.4cm
\noindent
$\alpha$--stable process:

All the proofs for the $\alpha$--stable case can be done by using the 
self--similarly of the transition probability $p_t$. We will prove only (c) 
and (d) to exemplify.
\vglue.3cm
\noindent
$\mbox{\rm (c)}\hfill
G_t\varphi(x)=\displaystyle\int^t_0
\displaystyle\int_{\errita^d}p_s(x-y)\varphi(y)dyds
=\displaystyle\int^t_0s^{-d/\alpha}
\displaystyle\int_{\errita^d}p_1(s^{-1/\alpha}(x-y))\varphi(y)dyds\hfill
$%
$$
=t^{-d/\alpha+1}\int^1_0r^{-d/\alpha}
\int_{\errita^d}p_1(t^{-1/\alpha}r^{-1/\alpha}(x-y))\varphi(y)dydr,
$$
hence
$$\lim_{t\rightarrow\infty}t^{1-d/\alpha}G_t\varphi(x)=
{1\over 1-d/\alpha}p_1(0)
\int_{\errita^d}\varphi(y)dy,$$
and (by Fourier transform)
$$p_1(0)={1\over (2\pi)^d}\int_{\errita^d}e^{-|z|^\alpha}dz.$$
\vglue.3cm
\noindent
$\mbox{\rm (d)} \hfill G'_t\varphi(x)=\displaystyle
\int_{\errita^d}p_t(x-y)\varphi(y)dy
=t^{-1}\int_{\errita^d}p_1(t^{-1/\alpha}(x-y))\varphi(y)dy,\hfill$
\vglue.3cm
\noindent
hence
$$\lim_{t\rightarrow\infty}tG'_t\varphi(x)=p_1(0)
\int_{\errita^d}\varphi(y)dy,$$
with
$$p_1(0)={1\over (2\pi)^d}\int_{\errita^d}e^{-|z|^d}dz,$$
and the result follows by l'H\^opital's rule.
\vglue.5cm
\noindent
Hierarchical random walk:

Recall that $a^{\gamma+1}=\displaystyle{1\over N}$, by (4.2.3),
and  the definitions of $\theta$ and the functions $h$ in Subsection 4.2 
(see (4.2.4)--(4.2.9)).

\vglue.5cm
\noindent
(a) By (3.1.3), up to a summand which converges to 0 exponentially fast as 
$t\rightarrow\infty,p_t(0,x)$ equals
\begin{eqnarray*}
(N-1)\sum^\infty_{j=1}{1\over N^j} e^{-\theta a^jt}
&=&(N-1)(\theta t)^{-(\gamma+1)}\sum^\infty_{j=1}e^{-\theta a^jt}(\theta 
a^j t)^{\gamma+1}\\
&=&q_\gamma t^{-(\gamma+1)}
\biggl(h_t^{(1,\gamma+1)}-\sum_{j\leq 0}
e^{-\theta a^jt}(\theta a^jt)^{\gamma+1}
\biggr).
\end{eqnarray*}
To conclude the proof of (a) if suffices to observe that
$\sum_{j\leq 0}e^{-\theta a^jt}(\theta a^jt)^{\gamma+1}\rightarrow 0$ as 
$t\rightarrow\infty$. Indeed for all $t$ and all negative integers $j$ we 
have
$$(\theta a^jt)^{\gamma+1} e^{-\theta a^jt}\leq const. {1\over \theta 
a^jt},$$
hence
$\sum\limits_{j\leq 0}e^{-\theta a^jt} (\theta a^jt)^{\gamma+1}$
is majorized by $\displaystyle{1\over t}$ times a convergent geometric 
series.
\vglue.5cm
\noindent
(b) Because of (a) it suffices to compute
\begin{eqnarray*}
\int^\infty_t s^{-(\gamma+1)}h_s^{(1,\gamma+1)}ds
&=&\int^\infty_t\sum^\infty_{j=-\infty}(\theta a^j)^{\gamma+1}e^{-\theta 
a^js}ds\\
&=&\sum^\infty_{j=-\infty}(\theta a^j)^\gamma e^{-\theta 
a^jt}=t^{-\gamma}h_t^{(1,\gamma)}.
\end{eqnarray*}
(c) Since $G_t(0,x)=\int^t_0p_s(0,x)ds$, we infer from (a) that
\begin{eqnarray*}
G_t(0,x)&\sim&\int^t_0s^{-(\gamma+1)}h_s^{(1,\gamma+1)}ds\\
&=&q_\gamma\int^t_0\sum^\infty_{j=-\infty}(\theta a^j)^{\gamma+1}e^{-\theta 
a^js}ds\\
&=&q_\gamma\sum^\infty_{j-\infty}(\theta a^j)^\gamma(1-e^{-\theta a^jt})\\
&=&q_\gamma t^{-\gamma}h_t^{(2,\gamma)}.
\end{eqnarray*}
(d) Since $a=\displaystyle{1\over N}$ for $\gamma=0$, up to a summand which 
is uniformly bounded in $t$, $G_t(0,x)$ equals
$$
{N-1\over \theta}\sum^\infty_{j=1}(1-e^{-\theta a^j t}).
$$
Since $a^{j+1}\leq a^y\leq a^i$ for $y\in [j,j+1]$, then
\vglue.3cm
\noindent
$\hfill
\displaystyle\int^\infty_1 (1-e^{-a^yt})dy
\leq
\displaystyle\sum\limits^\infty_{j=1}(1-e^{-a^jt})\leq
\displaystyle\int^\infty_1 (1-e^{-a^y a^{-1}t})dy. \hfill (5.3.1)
$\vglue.3cm
\noindent
Now,
$$
\int^\infty _1 (1-e^{-a^yt})dy ={1\over -\log a}\int^a_0 {1-e^{-zt}\over z}
dz={1\over -\log a}\int^{at}_0 {1-e^{-r}\over r}dr,
$$
and by l'H\^opital's rule,
$$
\int^{at}_0 {1-e^{-r}\over r}dr\sim \log t.
$$
Hence it follows from (5.3.1) that
$$
\sum^\infty_{j=1}(1-e^{-a^jt})\sim {\log t\over -\log a},
$$
which implies the result.
\vglue.3cm
\noindent
(e) Since $G^2_t(0,x)=\int^t_0 \int^t_0 p_{u+v}(0,x)dudv$, we infer from 
(a) that
\begin{eqnarray*}
G^2_t (0,x)&\sim& q_\gamma \sum^\infty_{j=-\infty}(\theta a^j)^{\gamma 
+1}\int^t_0 \int^t_0 e^{-\theta a^j(u+v)}dudv\\
&=& q_\gamma \sum^\infty_{j=-\infty}(\theta a^j)^{\gamma -1}(1-e^{-\theta 
a^j t})^2\\
&=&q_\gamma t^{1-\gamma}h^{(3, \gamma -1)}_t.
\end{eqnarray*}
(f) Since $a=\displaystyle{1\over N^{1/2}}$ for $\gamma =1$, up to a 
summand which is uniformly bounded in
$t$, $G^2_t(0,x)$ equals
$${N-1\over \theta ^2}\sum^\infty_{j=1}(1-e^{-\theta a^jt})^2.$$
 The same argument used for (d) shows that
$$
\sum^\infty_{j=1}(1-e^{-a^jt})^2\sim {-\log t\over \log a},
$$
and the result follows.

\noindent
(g) Since $G^2_t G(0,x)=\int^\infty_0 \int^t_0 \int^t_0
p_{s+u+v}(0,x)dudvds$, we infer from (a) that
\begin{eqnarray*}
G^2_tG(0,x)&\sim& q_\gamma \sum^\infty_{j=-\infty}(\theta a^j)^{\gamma 
+1}\int^\infty_0\int^t_0\int^t_0 e^{-\theta a^j (s+u+v)} dudvds\\
&=&q_\gamma \sum^\infty_{j=-\infty}(\theta a^j)^{\gamma -2} (1-e^{-\theta 
a^j t})^2\\
&=&q_\gamma t^{2-\gamma}h^{(3, \gamma -2)}_t.
\end{eqnarray*}
(h) Since $a=\displaystyle{1\over N^{1/3}}$ for $\gamma =2$, up to a 
summand which is uniformly bounded in $t$, $G^2_t G(0,x)$ equals
$${N-1\over \theta^3}\sum^\infty_{j=1}(1-e^{-\theta a^jt})^2,$$
and the proof is the same as for (f).
\hfill $\Box$
\vglue.5cm
\noindent
{\it Proof of Lemma 3.1.2:}
\vglue.5cm
By  the observation before Lemma 2.4.1,
$||G^j\varphi||<\infty$
for ${\cal C}^+_c(S)$
if and only if $j<\gamma+1$, i.e. $\alpha j<d$ for the $\alpha$--stable 
case, and $c>N^{(j-1)/j}$ for the $c$--hierarchical case.

The expression for $G^j$ can be obtained by formula (2.1.1). For the 
$c$--hierarchical case the form of the semigroup $T_t$ is explicit
from
 the transition probability $p_t(0,x)$ given in (3.1.3). For the 
$\alpha$--stable case the transition probability $p_t(0,x)$ is not known 
explicitly for general $\alpha$, but its Fourier transform,
$$\int_{\errita^d}e^{-ix\cdot z}p_t(0,x)dx=e^{-t|z|^\alpha},$$
can be used for the proof. $\hfill\Box$
\vglue.5cm
\noindent
{\it Proof of Lemma 3.1.3:}
\vglue.5cm
We sketch the main idea of the proof. By Lemma 3.1.1(a), $T_t\sim 
t^{-(\gamma+1)}$. By Lemma  2.4.1, $\gamma>k+1$. Hence 
$T_t=o(t^{-(k+2)}).\hfill\Box$
\vglue.5cm
\noindent
{\it Proof of Proposition 3.1.1:}
\vglue.5cm
Let $h_t$ denote any of the functions defined in (4.2.5) -- (4.2.7).
We write $h_t$
 as $h_t=h^{(-)}_t+h^{(+)}_t$, where $h^{(-)}_t$ and $h^{(+)}_t$ stand for 
the sums $\sum\limits_j$ with $j<0$ and $j\geq 0$, respectively. Since 
$\lim_{t\rightarrow\infty}h^{(-)}_t=0$ (as in the proof of Lemma 3.1.1 for 
the hierarchical case),
and $h_t$ is periodic in a logarithmic scale,  in order to prove that
$$L_1\leq \inf_th_t\leq \sup_t h_t\leq L_2,$$
for some positive constants $L_1$ and $L_2$,
it suffices to show that
\vskip.3cm
\noindent
$\hfill
L_1\leq\liminf\limits_{t\rightarrow\infty}
h^{(+)}_t\leq \limsup\limits_{t\rightarrow\infty}
h^{(+)}_t\leq L_2.\hfill (5.3.2)
$\vskip.5cm
We will prove (5.3.2) for $h^{(1,\zeta)}_t$. Using the formula
$$q^{-j}={q\log q\over q-1}\int^{j+1}_j q^{-y}dy, q>0,\,\,\,\, q\neq 1,$$
and $a^{j+1}\leq a^y\leq a^j$ for $j\leq y\leq j+1$ (since $0<a<1$), we 
obtain
$$
{a^{-\zeta}\log a^{-\zeta}\over a^{-\zeta}-1}
\int^\infty_0
a^{y\zeta}e^{-a^{-y}a^{-1}t}dy\leq
\sum^\infty_{j=0}a^{j\zeta}e^{-a^jt}
\leq
{a^{-\zeta}\log a^{-\zeta}\over a^{-\zeta}-1}
\int^\infty_0a^{y\zeta }e^{-a^yt}dy.
$$
We have
$$\int^\infty_0
a^{y\zeta}e^{-a^yt}dy=
{1\over -\log a}\int^1_0z^{\zeta-1}e^{-zt}dz=
{t^{-\zeta}\over -\log a}\int^t_0 r^{\zeta-1}e^{-r}dr,$$
and since
$\int^t_0 r^{\zeta-1}e^{-r}dr\rightarrow\Gamma(\zeta)$ as 
$t\rightarrow\infty$, putting these results together we obtain
$$
{\Gamma(\zeta+1)\over a^{-\zeta}-1}\,\,\leq\,\, 
\liminf_{t\rightarrow\infty}
t^\zeta\sum^\infty_{j=0}a^{j\zeta}e^{-a^jt}
\leq\limsup_{t\rightarrow\infty}t^\zeta\sum^\infty_{j=0}a^{j\zeta}e^{-a^  
jt}\leq
{a^{-\zeta}\Gamma(\zeta+1)\over a^{-\zeta}-1},
$$
which finishes the proof.

This method can be used for the other functions
$h_t$
as well, with slightly more elaborate calculations. For 
$\widetilde{h}^{(2,\zeta)}_t$ and
$\widetilde{h}^{(3,\zeta)}_t$ we need to
use the fact that the functions $x\mapsto e^{-x}-1+x$ and $x\mapsto 
2e^{-x}-{1\over 2}e^{-2x}+x-{3\over 2}$, respectively, are increasing.
However, we can obtain  bounds for $\widetilde{h}^{(2,\zeta)}_t$ and
$\widetilde{h}^{(3,\zeta)}_t$
from the bounds for $h^{(2,\zeta)}_t$ and $h_t^{(3,\zeta)}$, simply by 
dividing
them
by $1-\zeta$. This is clear from the form of $\widetilde{h}_t$ in Lemma 
2.4.3.
$\hfill\Box$
\vglue.5cm
\noindent
{\it Proof of Theorem 3.2.1:}
\vglue.5cm
The proof is a direct application of Theorems 2.2.1, 2.2.1, 2.2.3, Lemmas 
3.1.1, 3.1.2, 3.1.3, and Corollary 3.1.1. $\hfill\Box$
\vglue1cm
\noindent
{\bf 5.4. Conditions for the results on infinite variance branching 
results}
\vglue.5cm
\noindent
{\it Proof of Porposition 3.2.1:}
We have to show that
\vskip.3cm
\noindent
$\hfill d>\alpha\left(1+\displaystyle{1\over\beta}\right)\hfill (5.4.1)
$
\vskip.3cm
\noindent
is necessary and sufficient for condition (2.5.1) of Theorem 2.5.1.

That (5.4.1) implies (2.5.1) is proved by Iscoe$^{(27)}$.
For the converse, note that
$$G 1_B(x)\geq k|x|^{-(d-\alpha)}\,\,\mbox{\rm for }\,\,|x|\geq 2,$$
where $B$ denotes the unit ball centered at the origin and $k$ is some 
positive constant. Hence, if $d\leq\alpha(1+1/\beta), (G 1_B)^{1+\beta}$ is 
not $\lambda$--integrable. $\hfill\Box$

\vglue.5cm
For the proof of Proposition 3.2.2 we need some preliminary results.

Let $R_{\infty}^{1}$ and $R_{\infty}$ be the canonical measure of the
equilibrium of the  particle system (started
off in the Poisson system
$\Pi_{\lambda}^{}$ with intensity $\lambda$) and that of the
superprocess (started off in $\lambda$), respectively (Appendix,
Subsection A.1).

\medskip

\noindent
\textbf{Lemma 5.4.1.} Assume that $\beta_{2}<\beta_{1}$ and let
$\varphi:\R^{d}\to [0,\infty]$. Then

\noindent
\textrm{(a)}
\[
\int_{{\cal M}_\tau(S)}^{}\langle \nu, 
\varphi\rangle^{1+\beta_{2}}R_{\infty}
(d\nu)\leq
\int_{{\cal M}_\tau(S)}^{}\langle \nu, \varphi\rangle^{1+\beta_{2}}R^{1}  
_{\infty}
(d\nu).
\]

\noindent
\textrm{b)}  If $\varphi$ is $\lambda$-integrable, then
	\[
	\int_{{\cal M}_\tau(S)}^{}\langle \nu, 
\varphi\rangle^{1+\beta_{2}}R^{1}_{\infty}
(d\nu) \leq C < \infty ,
	\]
	where the constant
 $C$ depends only on $d,\alpha,\beta_{1}, \beta_{2}$ and
$\langle \lambda, \varphi\rangle$.

\bigskip

\noindent
{\it Proof:}
\textrm{(a)}
We have,
by Jensen's inequality and (A.1.10),
\begin{eqnarray*}
	\int_{}^{} \langle \nu, \varphi\rangle^{1+\beta_{2}}
R_{\infty}\left(d\nu\right)& = &
	\int_{}^{}\left(\int_{}^{}\langle \mu,
	\varphi\rangle\Pi_{\nu}\left(d\mu\right)\right)^{1+\beta_{2}}
R_{\infty}\left(d\nu\right)\\
	 & \leq & \int_{}^{}\int_{}^{}\langle \mu,
\varphi\rangle^{1+\beta_{2}}
	 \Pi_{\nu}\left(d\mu\right) R_{\infty}\left(d\nu\right) \\
	 & = & \int_{}^{} \langle \mu, \varphi\rangle^{1+\beta_{2}}
	 R_{\infty}^{1}\left(d\nu\right).
\end{eqnarray*}


\noindent
\textrm{(b)}  By the Palm formula (A.2.1) it
suffices to show for some $\delta>\beta_{2}$,
\[
\int_{}^{}\langle \nu, \varphi\rangle^{\delta} (R_{\infty}^{1})_{x}
(d\nu) < C<\infty,
\]
\noindent
where the constant $C$  does not depend on $x$.
We can choose $\delta \in (\beta_{2},\beta_{1})$ such that
$d\beta_{2}/\alpha> 1$.

We will use the tree representation of
$(R_{\infty}^{1})_{x}^{red}$ given in (A.3.2), and we denote
$Z_{t,i}=\langle X_{t}^{W^x_{t},i}, \varphi\rangle$. Since
$\left(\sum_{j=1}^{n}a_{j}\right)^{\delta} \leq 
\sum_{j=1}^{n}a_{j}^{\delta}$,
for all nonnegative sequences $(a_{j})$ and $0<\delta<1$, we have

\begin{eqnarray*}
 	\int\limits_{}^{}\langle \nu,
\varphi\rangle^{\delta}\left(R_{\infty}^{1}\right)^{red}_{x}
 	\left(d\nu\right) & = & E
\biggl(\biggl(\int_{0}^{\infty}\biggl(\sum_{i=1}^{N_{t}}Z_{t,i}\biggr)
 	\pi(dt)\biggr)^{\delta}\biggr) \\
 	 & = & E
\biggl(E\biggl[\biggl(\int_{0}^{\infty}\biggl(\sum_{i=1}^{N_{t}}Z_{t,i}\  
biggr)
 	 \pi(dt)\biggr)^{\delta}\biggl| \, \pi, W^x\biggr]\biggr)
\\
 	 & \leq & E
\biggl(E\biggl[\biggl(\int_{0}^{\infty}\biggl(\sum_{i=1}^{N_{t}}Z_{t,i}\  
biggr)^
{\delta}
 	 \pi(dt)\biggr)\bigg| \, \pi, W^x\biggr]\biggr)\\
 	 & = & E \biggl(\int_{0}^{\infty}E
 	 \biggl[\biggl(\sum_{i=1}^{N_{t}}Z_{t,i}\biggr)^{\delta} \bigg| \,
W^x\biggr] \pi(dt)\biggl).
\end{eqnarray*}

\noindent
Now, by the self--similarity of the $\alpha$--stable process,
 \begin{eqnarray*}
 	E\left[Z_{t,i}|W^x\right] & = & \int_{}^{} \varphi\left(y\right)
p_{t}\left(W^x_{t},
 	y\right) \lambda \left(dy\right)\\
 	 & = & \int_{}^{}\varphi \left(y-W^x_{t}\right) p_{t}\left(0,y\right)
\lambda\left(dy\right) \\
 	 & = & \frac{1}{t^{d/\alpha}}
\int_{}^{}\varphi\left(y-W^x_{t}\right)p_{1}\left(0, y
 	 t^{-1/\alpha}\right)dy \\
 	 & \leq & K \frac{1}{t^{d/\alpha}} \langle \lambda, \varphi\rangle ,
 \end{eqnarray*}
 for some real constant $K$ not depending on $W^x$ and  $\varphi$.
Hence, by H\"older's inequality,
\begin{eqnarray*}
	E \biggl[\biggl(\sum_{i=1}^{N_{t}}Z_{t,i}\biggr)^{\delta} \bigg| W^x
\biggr] & = &
	\sum^\infty_{n=0}{}\, E \biggl[ \biggl(
\sum_{i=1}^{n}Z_{t,i}\biggr)^{\delta}\bigg| W^x\biggr] P
[N_{t}=n]\\
	 & = & \sum^\infty_{n=0}n^{\delta} \  E
\biggl[\biggl(\frac{1}{n}\sum_{i=1}^{n}Z_{t,i}\biggr)^{\delta}\bigg|W^x\  
biggr]P[N_{t}=n] \\
	 & \leq & \sum^\infty_{n=0} n^{\delta}\left(E \left[
	 \frac{1}{n}\sum_{i=1}^{n}Z_{t,i}\bigg| W^x\right]\right)^{\delta}
P\left[N_{t}=n\right] \\
	 & \leq  & \left(K\langle \lambda, \varphi\rangle\right)^{\delta}
t^{-d\delta/\alpha}
	 \sum^\infty_{n=0}n^{\delta}  P\left[N_{t}=n\right].
\end{eqnarray*}
\noindent For $(1+\beta_{1})$--branching and $\delta<\beta_{1}$,
\[
\sum^\infty_{n=0}n^{\delta} P\left[N_{t}=n\right]<\infty.
\]

\noindent On the other hand, since $d\delta/\alpha>1$, we have
\[
E\left(\int_{1}^{\infty}t^{-d\delta/\alpha}\pi \left(dt\right)
\right)< \infty.
\]
Putting these results together finishes the proof.
$\hfill \Box$

\bigskip

\noindent
\textbf{Corollary 5.4.1.}
In the setting of Lemma 5.4.1(b),
\[
\int_{}^{} \langle \nu, \varphi\rangle^{\beta_{2}} (R_{\infty})_{x}(d
\nu) < \infty, \qquad x\in \R^{d}.
\]

\medskip

\noindent
{\it Proof:}
Let $\psi$ be strictly positive and $\lambda$-integrable. Then, by the Palm 
formula (A.2.1),
\begin{eqnarray*}
	\int\!\int \langle \nu, \varphi\rangle^{\beta_{2}}
(R_{\infty})_{x}(d
\nu) \psi(x) \lambda(dx) & = & \int_{}^{}  \langle \nu,
\varphi\rangle^{\beta_{2}}
\langle \nu, \psi\rangle R_{\infty}(d\nu)\\
	 & \leq  & \int_{}^{} \langle \nu, \varphi+\psi
	 \rangle^{1+\beta_{2}} R_{\infty}(d\nu) < \infty.
\end{eqnarray*}
This shows that the assertion of the corollary holds for $\lambda$--almost
all $x$. The shift--invariance of the system implies that it is in fact
true for all $x\in \R^{d}$.
$\hfill \Box$

\bigskip

\noindent
\textbf{Lemma 5.4.2.}
\[
(R_{\infty}^{1})_{x}= \delta_{\delta_{x}} \ast
\int_{{\cal M}_\tau(S)}^{}\Pi_{\nu}(\cdot) (R_{\infty})_{x}(d\nu)
\,\,\,\,\mbox{\rm for}\,\, \lambda-
\mbox{almost all }x.
\]

\noindent
{\it Proof:}
Since $R_{\infty}^{1}$ has intensity measure  $\lambda$,
by (A.1.11),  the
assertion is obtained from the following chain of equalities,
where we use the Palm formula (A.2.1) for
$\Pi_\nu$ and for $R_\infty$,
the fact that $(\Pi_\nu)_x=\delta_{\delta_x} *\Pi_\nu$,
and (A.1.10),
\begin{eqnarray*}
	\int\!\int\!\int\!\!\!\!&\!\!\! f(x) &\!\!\!\!\! F(\mu)
(\delta_{\delta_{x}} \ast
	\Pi_{\nu}) (d\mu) (R_{\infty})_{x}(d\nu)\lambda (dx)  \\
	& = & \int\!\int\!\int f(x) F(\mu)
(\delta_{\delta_{x}} \ast
	\Pi_{\nu}) (d\mu) \nu(dx) R_{\infty}(d\nu)  \\
	 & = & \int\!\int_{}^{} F(\mu)\int_{}^{}f(x) \mu (dx)
	 \Pi_{\nu}(d\mu) R_{\infty}(d\nu) \\
	 & = & \int_{}^{}F(\mu) \int_{}^{}f(x) \mu(dx) R_{\infty}^{1}(d\mu).
\end{eqnarray*}
\vskip-1cm
$\hfill \Box$


\bigskip

\noindent
{\it Proof of Proposition 3.2.2:} We have to prove that
\vskip.3cm
\noindent
$\hfill
\beta_2<\beta_{1}<1 \hfill  (5.4.2)
$\vskip.3cm
\noindent
and
\vskip.3cm
\noindent
$\hfill
d>\alpha\left(1+\displaystyle\frac{1}{\beta_{2}}
\left(
1+\displaystyle\frac{1}{\beta_{1}}\right)\right).\hfill (5.4.3)
$
\vskip.3cm
\noindent
are necessary and sufficient for condition (2.5.2) of Theorem 2.5.2.

\noindent
1. Assume that (5.4.2) and (5.4.3) hold.
Because of Lemma 3.1.2 and (4.1.2), we have to show that
\[
\int_{}^{}\left(\
\int\!\!\int
\frac{\varphi(y)}{|x-y|^{d-\alpha}}\ dy \  \mu (dx) \right)^{1+ \beta_{2}}
R_{\infty}^{1} (d\mu)< \infty,\quad \varphi\in{\cal C}^+_c(\erre^d).
\]
First we observe that (Iscoe$^{(27)}$, Lemma 5.3)
\vskip.3cm
\noindent
$\hfill
\displaystyle\int_{}^{}
\displaystyle\frac{\varphi(y)}{|x-y|^{d-\alpha}}\  dy \leq  const.
(1\vee |x|)^{-d+\alpha}.\hfill  (5.4.4)
$\vskip.3cm
\noindent
Hence from Lemma 5.4.1(b) (denoting by $B_{r}$  the ball
with radius $r$ centered at $0$) we obtain
\[
\int_{}^{}\left( \int_{B_{1}}^{} \int_{\R^{d}}^{}
\frac{\varphi(y)}{|x-y|^{d-\alpha}}\  dy \ \mu(dx)
\right)^{1+\beta_{2}}\ R_{\infty}^{1}(d\mu) < \infty ,
\]
and it suffices to show that
\[
A:= \int_{}^{}\left( \int_{\R^{d}}^{}1_{B_{1}^{c}}(x) \
|x|^{-d+\alpha} \ \mu (dx) \right)^{1+\beta_{2}} \ R_{\infty}^{1}
(d\mu)<\infty.
\]

Using the Palm formula (A.2.1) and Lemma 5.4.2, we have
\begin{eqnarray*}
	A & = & \int1_{B_{1}^{c}}(x) \,
|x|^{-d+\alpha} \,
\int_{}^{}\left(\int_{}^{} 1_{B_{1}^{c}}(z) \,|z|^{-d+\alpha}
\,\mu(dz)\right)^{\beta_{2}} \,(R_{\infty}^{1})_{x} (d\mu)
\lambda(dx)\\
	 & = & \int 1_{B_{1}^{c}}(x) \,|x|^{-d+\alpha}
	 \int\int\left(1_{B_{1}^{c}}(x) \,|x|^{-d+\alpha}
	 + \int 1_{B_{1}^{c}}(z) \,|z|^{-d+\alpha}
\,\mu(dz)\right)^{\beta_{2}} \\
	 & & \hspace{10cm}  \cdot\Pi_{\nu}(d\mu) (R_{\infty})_{x}(d\nu)
\lambda(dx) \\
	 & \leq & \int_{}^{}1_{B_{1}^{c}}(x) \,|x|^{-d+\alpha}
	 (|x|^{-d+\alpha})^{\beta_{2}}\lambda(dx) \\
	 &&+ \int 1_{B_{1}^{c}}(x) \,|x|^{-d+\alpha}
	 \int \int \left( \int 1_{B_{1}^{c}}(z)
	 \,|z|^{-d+\alpha}  \mu(dz) \right)^{\beta_{2}}
         \Pi_{\nu}(d\mu) (R_{\infty})_{x}(d\nu) \lambda(dx).
\end{eqnarray*}

Since $\alpha+ \beta_{2}(-d +\alpha)<0$ by (5.4.3), the first
term  on the r.h.s. is finite.

Using H\"older's inequality, the second term can be bounded by
\[
B:=
\int_{B_{1}^{c}}^{}|x|^{-d+\alpha}\int\limits_{}^{}
\biggl(\int\limits
_{}^{}
|z|^{-d+\alpha}\nu(dz)\biggr)^{\beta_{2}} (R_{\infty})_{x}(d\nu) 
\lambda(dx),
\]
and by shift--invariance,
\[
B=
\int_{B_{1}^{c}}^{}|x|^{-d+\alpha}\int\limits_{}^{}
\biggl(\int\limits
_{}^{}|z-x|^{-d+\alpha}
\nu(dz)\biggr)^{\beta_{2}}(R_{\infty})_{0}(d\nu) \lambda (dx).
\]

The scaling property of $(R_{\infty})_{0}$ (Dawson and Perkins$^{(11)}$,
 Theorem 6.7)
yields
\begin{eqnarray*}
	B & = & \int_{B_{1}^{c}}^{}|x|^{-d+\alpha} |x|^{\alpha /
	\beta_{1}}\int\limits_{}^{}
        \biggl( \int\limits_{}^{}
	|\,z|x|\, -\, x \,|^{-d+\alpha}\nu(dz)\biggr)^{\beta_{2}}
	(R_{\infty})_{0} (d\nu) \lambda(dx) \\
	 & = & \int_{B_{1}^{c}}^{}|x|
	 ^{-d+\alpha+\alpha/\beta_{1}+\beta_{2}(-d+\alpha)}
	 \int\limits_{}^{}\biggl( \int\limits_{}^{}
	\left|z-\frac{x}{|x|}\right|^{-d+\alpha}\nu(dz)\biggr)^{\beta_{2
}}
	(R_{\infty})_{0} (d\nu) \lambda(dx).
\end{eqnarray*}

By isotropy of $(R_{\infty})_{0}$, the
integral w.r. to
$(R_{\infty})_{0}$
does not depend on $x/|x|$. For an arbitrary fixed
$x_{0}\in \R^{d}$ we put $e=x_{0}/|x_{0}|$ and we obtain
\[
B=\int_{B_{1}^{c}}^{}|x|
^{-d+\alpha+\alpha/\beta_1+\beta_{2}(-d+\alpha)}\lambda(dx)
 \int\limits_{}^{}\biggl( \int\limits_{}^{}
	|z-e|^{-d+\alpha}\nu(dz)\biggr)^{\beta_{2}}
	(R_{\infty})_{0} (d\nu) .
\]
Since
$\alpha+\alpha/\beta_{1}+\beta_{2}(-d+\alpha)<0$ by
 (5.4.3), the integral w.r. to $\lambda$ is finite.
Hence we will be done  if we can show that the
integral w.r. to $(R_\infty)_0$
is finite as well.

We will show that
\[
C_{1}  :=  \int_{}^{}\biggl(\int_{B_{2}}^{} |z-e|^{-d+\alpha}
\nu(dz)\biggr)^{\beta_{2}} (R_{\infty})_{0} (d\mu)
\]
and
\[
C_{2}  :=  \int_{}^{}\biggl(\int_{B_{2}^{c}}^{} |z-e|^{-d+\alpha}
\nu(dz)\biggr)^{\beta_{2}} (R_{\infty})_{0} (d\mu)
\]
are finite.

Let
\[
g(z)= 1_{B_{2}}(z) |z-e|^{-d+\alpha} .
\]
By Corollary 5.4.1
we have
\[
C_{1}  =  \int_{}^{}\left(\int g(z)
\nu(dz)\right)^{\beta_{2}} (R_{\infty})_{0} (d\mu)<\infty,
\]
since $\langle \lambda,g\rangle < \infty$.
On the other hand, for $z\in B_{2}^{c}$,
$
|z-e|\geq \frac{1}{2}|z|.
$
Therefore
\[
C_{2}\leq const. \int_{}^{}
\biggl(\int_{B_{2}^{c}}^{}
|z|^{-d+\alpha}
\nu(dz)\biggr)^{\beta_{2}} (R_{\infty})_{0} (d\nu).
\]

To show the finiteness of the latter integral we will decompose the
random variable $\int_{B_{2}^{c}}^{} |z|^{-d+\alpha}
\nu(dz)$ into a sum of terms whose $L^{\beta_{2}}$--norms add up to 
something
finite.
Put $D_{k}:=B_{2^{k+1}} \backslash B_{2^{k}}$, $k=0,1,\ldots$
Then, again by the scaling property of $(R_{\infty})_{0}$  we have
\begin{eqnarray*}
	 &&\int_{}^{} \left(\int_{D_{k}}^{} |z|^{-d+\alpha}
     \nu(dz)\right)^{\beta_{2}} (R_{\infty})_{0} (d\nu) \\
	&&= (2^{k})^{\alpha/\beta_{1}}\int_{}^{}
	 \left(\int_{D_{0}}^{}(2^{k}|z|)^{-d+\alpha}
	 \nu(dz)\right)^{\beta_{2}} (R_{\infty})_{0} (d\nu) \\
	 & &\leq  2^{k\left(\alpha/\beta_1+
	 \beta_{2}(-d+\alpha)\right)}
	  \int_{}^{}(\nu (D_{0}))^{\beta_{2}} (R_{\infty})_{0} (d\nu).
\end{eqnarray*}

The $L^{\beta_{2}}$--norm of the $k$--th summand is thus bounded by
$const.
2^{k\left(\alpha/\beta_1\beta_2-d+\alpha\right)}$,
which is summable since $\alpha/\beta_{1}\beta_{2}-d+\alpha< 0$
due to  (5.4.3).

\bigskip

\noindent
\textrm{2.} To prove the converse  we assume that
 (2.5.2) holds.

\medskip

\noindent \textrm{(a)}
Assume that $\beta_{2}\geq \beta_{1}$, $\beta_{1}<1$.
We will show that, for $\varphi \in C_{c}(S)$, $\varphi \geq 0$,
$\varphi \neq0$,
\[
\int_{}^{}\langle\mu, G \varphi\rangle^{1+\beta_{2}}
R_{\infty}^{1}(d\mu)=\infty .
\]

First note that by the Palm formula (A.2.1),
\vskip.3cm
\noindent
$\hfill
\displaystyle\int_{}^{}\langle\mu, G \varphi\rangle^{1+\beta_{2}}
     R_{\infty}^{1}(d\mu)
     = \displaystyle\int_{}^{} G \varphi(x)
	 \displaystyle\int_{}^{}
	 \langle\mu, G 
\varphi\rangle^{\beta_{2}}(R_{\infty}^{1})_{x}(d\mu)\lambda(dx). \hfill 
(5.4.5)
$\vskip.3cm
\noindent
Choose $\psi:\R^{d} \to \R_{+}$ such that $G \varphi \geq \psi(
\cdot - x)$ provided that $|x|\leq 1$. Then the r.h.s. of (5.4.5) is
bounded  above by
\[
\int_{}^{}1_{\{|x|\leq 1\}}G\varphi(x)
\int_{}^{}\langle\mu, 
\psi\rangle^{\beta_{2}}(R_{\infty}^{1})_{0}(d\mu)\lambda(dx).
\]
By the tree representation  (A.3.2) of $(R_{\infty}^{1})_{0}$ we have
\vskip.3cm
\noindent
$\hfill
	\displaystyle\int_{}^{}\langle\mu,
\psi\rangle^{\beta_{2}}(R_{\infty}^{1})_{0}(d\mu)\geq
E \left(\left(\sum_{i=1}^{N}\langle X_{\sigma}^{W^0_{\sigma},i},
	 \psi\rangle\right)^{\beta_{2}}\right),\hfill (5.4.6)
$\vskip.3cm
\noindent
where $\sigma$ is exponentially distributed with parameter $V_{1}$,
and $N$ is distributed like any of the $N_{t}$. Since $E
N^{\beta_{2}}= \infty$, and since, conditioned on $\sigma$ and $W^0$,
the random variables $\langle X_{\sigma}^{W^0_{\sigma},i}, \psi\rangle$ are
i.i.d. with positive expectation, it follows from the law of large
numbers that the r.h.s. of (5.4.6) is infinite.

\medskip

\noindent
\textrm{(b)}
Now assume that $1\geq\beta_{1}>\beta_{2}$, and suppose that $\varphi(x)
\geq 1$ for $|x|\leq1$.
Then, for some $k>0$,
\[
G \varphi(x) \geq k |x|^{-(d-\alpha)} \qquad \mbox{ if }|x| \geq 2.
\]
Assume that $\int_{}^{} \langle\mu, G \varphi\rangle^{1+\beta_{2}}
R_{\infty}^{1}(d\mu)< \infty$.
Then, by Lemma 5.4.1(a),
\[
\int_{}^{} \langle\nu, G \varphi\rangle^{1+\beta_{2}}
R_{\infty}^{1}(d\nu)<\infty.
\]
Therefore, by the Palm formula (A.2.1) we have
\begin{eqnarray*}
	\infty & > & \int_{}^{} \langle\nu, G \varphi\rangle^{1+\beta_{2}}
     R_{\infty}^{1}(d\nu) \\
	 & \geq & k^{^{1+\beta_{2}}} \int_{}^{} 1_{\{|x|\geq
	 1\} }|x|^{-d+\alpha} \int_{}^{}\left(\int_{}^{}
	 1_{\{|z|\geq 1\}} |z|^{-d+\alpha}
	 \nu(dz)\right)^{\beta_{2}}(R_{\infty})_{x}(d\nu)\lambda(dx).
\end{eqnarray*}
By shift--invariance this equals
\[
k^{^{1+\beta_{2}}} \int_{}^{}
 1_{\{|x|\geq
	 1\} }|x|^{-d+\alpha} \int_{}^{}\left(\int_{}^{}
	 1_{\{|z-x|\geq 1\}} |z-x|^{-d+\alpha}
	 \nu(dz)\right)^{\beta_{2}}(R_{\infty})_{0}(d\nu)\lambda(dx).
\]
The scaling property of $(R_{\infty})_{0}$  permits to
rewrite the inner integral (with  $e(x)=
x/|x|$)
as
\[
|x|^{\alpha/\beta_{1}} \int_{}^{}\left(\int_{}^{}
	 1_{\{|\, z \, |x|\, - \, x\, |\, \geq \, 1\}}
	 |\, z \, |x|\, -\,  x\, |^{-d+\alpha}
	 \nu(dz)\right)^{\beta_{2}}(R_{\infty})_{0}(d\nu)
\]

\[
=|x|^{\alpha/\beta_1+\beta_{2}(-d+\alpha)}
\int_{}^{}\left(\int_{}^{}
	 1_{\{|x| |z-e(x)|\geq 1\}} |z-e(x)|^{-d+\alpha}
	 \nu(dz)\right)^{\beta_{2}}(R_{\infty})_{0}(d\nu).
\]
For $|x|\geq1$, the latter integral is bounded below by
$$\int\left(\int 1_{\{|z|\geq 
1\}}|2z|^{-d+\alpha}\nu(dz)\right)^{\beta_2}(R_\infty)_0(d\nu)>0.$$
Now,
\[
\int_{}^{}|x|^{-d+\alpha+\alpha/\beta_{1}+\beta_{2}(-d+\alpha)}
1_{\{|x|\geq 1\}}\lambda(dx) <\infty
\]
implies that
$\alpha+\alpha/\beta_{1}+\beta_{2}(-d+\alpha)<0$, or
equivalently, (5.4.3) holds. $\hfill\Box$

\vglue1cm
\noindent
{\bf APPENDIX}
\vglue.5cm
\noindent
{\bf A.1. Background on 1- and 2-level branching systems}
\vglue.5cm
We consider particle systems in a locally compact Abelian group $S$ with
Haar measure $\rho$. Recall that $T_t$ denotes the semigroup of the 
particle motion and $G$ the corresponding Green operator.

Let ${\cal C}(S)$ denote the space of bounded continuous functions on
$S,\,\,{\cal C}_0(S)$ the subspace of functions vanishing at infinity, and 
${\cal C}_c(S)$ that of functions with compact support. For a strictly 
positive function $\tau\in{\cal C}_0(S)$, let
$${\cal C}_\tau(S)=\{\varphi\in{\cal C}(S):\varphi \tau^{-1}\in{\cal 
C}_0(S)\}$$
with the norm
$||\varphi||_\tau=||\varphi \tau^{-1}||$. We assume that $\tau$ is such 
that
$t\mapsto T_t\varphi$ is a continuous curve in
$({\cal C}_\tau(S),||\cdot||_\tau)$ for each $\varphi\in{\cal C}_\tau(S)$. 
For example, in the case of the $\alpha$-stable motion in $\erre^d$ we may 
take
$\tau(x)=(1+|x|^2)^{-q}$ with $d/2<q<(d+\alpha)/2$
(Dawson and Gorostiza$^{(6)}$). The subspaces of non--negative
functions are indicated with the
superscript `+', e.g. ${\cal C}^+_\tau(S)$. Let ${\cal M}_\tau(S)$ denote 
the space of non--negative Radon measures $\mu$ on $S$ such that 
$\langle\mu,\tau\rangle<\infty$, endowed with the smallest topology which 
makes the maps
$\mu\mapsto\langle\mu,\varphi\rangle$ continuous for all
$\varphi\in{\cal C}^+_c(S)\cup \{\tau\}$.
We assume that $\rho\in{\cal M}_\tau(S)$. The subspace of ${\cal 
M}_\tau(S)$ of integer--valued measures is designated by ${\cal 
N}_\tau(S)$.

The Laplace
functional of the occupation time of the 1--level branching particle
system  $X_t$
with $(1+\beta)$--branching at rate $V$ and started off from a
Poisson system  with intensity $\rho$ is given by
\vskip.2cm
\noindent
$\hfill
E\mbox{\rm exp}\left\{-
\left\langle
\displaystyle\int^t_0X_sds,\varphi\right\rangle\right\}
=\,\mbox{\rm exp}\{-\langle\rho,u_\varphi(t)\rangle\},\,\,
\varphi\in{\cal C}^+_\tau(S),\hfill (A.1.1)$
\vskip.2cm
\noindent
where $u_\varphi(x,t)$ with values in $[0,1]$
 is the unique solution of the non--linear evolution equation
\vskip.2cm
\noindent
$\hfill
u_\varphi(t)=-\displaystyle{V\over
1+\beta}\displaystyle\int^t_0
T_{t-s}(u_\varphi(s)^{1+\beta})ds+\displaystyle\int^t
_0T_{t-s}(\varphi(1-u_\varphi(s)))ds.\hfill (A.1.2)$
\vskip.2cm
\noindent
This is shown by the same argument of Theorem 5 in
Gorostiza and L\'opez--Mimbela$^{(18)}$ (formulas
(4.8) and (4.9)). It follows that
\vskip.2cm
\noindent
$\hfill
u_\varphi(t)\leq\displaystyle\int^t_0
T_{t-s}
(\varphi(1-u_\varphi(s)))ds\leq G_t\varphi.
\hfill (A.1.3)$
\vskip.3cm
Let $U_t$ denote the semigroup of the 1--level branching particle system.
We have
\vskip.2cm
\noindent
$\hfill
U_t(\langle\cdot,\varphi\rangle)(\mu)=\langle\mu,T_t
\varphi\rangle,\quad
\varphi\in {\cal C}^+_\tau(S),\mu\in{\cal N}_\tau(S),\hfill (A.1.4)$
\vskip.2cm
\noindent
and, if $\beta=1$,
\begin{eqnarray*}
\qquad\qquad\qquad
U_t(\langle\cdot,\varphi\rangle\langle\cdot,\psi\rangle)(\mu)
&=&\langle\mu,T_t\varphi\rangle\langle\mu,T_t\psi\rangle
+\langle\mu,T_t(\varphi\psi)-T_t\varphi\cdot T_t\psi\rangle\\
&+&V\displaystyle\int^t_0\langle\mu,T_s(T_{t-s}\varphi\cdot T_{t-s}\psi)
\rangle ds,\quad \varphi,\psi\in{\cal C}^+_\tau(S),\,\, \mu\in{\cal N}_\tau
\qquad (A.1.5)
\end{eqnarray*}
The formulas (A.1.4) and (A.1.5)
can be derived by martingale methods from the Markov
property of the system (see  e.g. Gorostiza and Rodrigues$^{(20)}$ for
 explicit calculations of this type). In particular,
\vskip.2cm
\noindent
$\hfill
U_t(\langle\cdot,\varphi\rangle)(\delta_x)=T_t\varphi(x),\hfill (A.1.6)$%
$$\hfill
\qquad\qquad\qquad\qquad
U_t(\langle\cdot,\varphi\rangle\langle\cdot,\psi\rangle)(\delta_x)=
T_t(\varphi\psi)(x)
+V\displaystyle\int^t_0T_s(T_{t-s}\varphi\cdot T_{t-s}\psi)(x)ds.
\hfill\qquad\qquad\qquad\quad (A.1.7)$$
\vskip.2cm
\noindent
We have from (A.1.4)
$$
U_t(\langle\cdot,\varphi\rangle\langle\cdot,\psi\rangle)(\mu)
\leq \langle\mu,T_t\varphi\rangle\langle\mu,T_t\psi\rangle
+\langle\mu,T_t(\varphi\psi)\rangle+V
\displaystyle\int^t_0\langle\mu,T_s(T_{t-s}\varphi\cdot T_{t-s}\psi)
\rangle ds,
$$%
$
\qquad\qquad\qquad\qquad\qquad\qquad\qquad\qquad\qquad\qquad\qquad\qquad
\qquad\qquad\qquad\qquad\qquad
\varphi,\psi\in{\cal C}^+_\tau(S).\hfill (A.1.8)$
\vglue.4cm
If the 1--level branching system is persistent, it
 has a ``Poisson type'' equilibrium (in the
sense of Liemant et al$^{(31)}$, section 2.3), which is an infinitely 
divisible random element
of ${\cal M}_\tau(S)$. Its canonical measure,   which
is a measure on ${\cal M}_\tau(S)$, is denoted by $R^1_\infty$.
A sufficient condition for persistence is
$$\int^\infty_0\langle\rho,(T_t\varphi)^{1+\beta}\rangle dt<\infty,
\,\, \varphi\in{\cal C}^+_c(S)$$
(Gorostiza and Wakolbinger$^{(22)}$, Theorem 2.1).

For each $t>0$, the random measure $X_t$ is infinitely divisible and its
canonical measure $R_t$ has the form
\vskip.3cm
\noindent
$\hfill R^1_t=\displaystyle\int_SP[X^x_t\in(\cdot)]\rho(dx),\hfill (A.1.9)
$
\vskip.3cm
\noindent
where $X^x_t$ corresponds to the branching system starting with
a single
ancestor in $x$ at time $0$ (Gorostiza and Wakolbinger$^{(21)}$,
formula (3.1),
Liemant et al$^{(31)}$).

The measure $R^1_\infty$ is the ``Poissonization'' of the equilibrium 
canonical measure
$R_\infty$ of the superprocess counterpart of
the particle system, i.e.,
\vskip.3cm
\noindent
$\hfill
R^1_\infty=\displaystyle\int_{{\cal M}_\tau(S)}
\Pi_\nu (\cdot) R_\infty(d\nu),\hfill (A.1.10)
$
\vskip.3cm
\noindent
where $\Pi_\nu$ is
the distribution of
a Poisson random measure  on $S$ with intensity measure
$\nu$.
Indeed, in Gorostiza et al$^{(19)}$ it is shown that the distibution
$L_{t}^{1}$ of the branching particle system $X_{t}$ is a Cox process, 
i.e.,
\[
L_{t}^{1}=\int_{{\cal M}_\tau(S)}^{}\Pi_{\nu}(\cdot) L_{t}(d\nu),
\]
where $L_{t}$ is the distribution of the superprocess counterpart of
$X_{t}$. By continuity and the assumed persistence, this relation
carries over to $t= \infty$:
\[
L_{\infty}^{1}=\int_{{\cal M}_\tau(S)}^{}\Pi_{\nu}(\cdot) L_{\infty}(d\nu).
\]
Together with the L\'evy-Khinchin formula (Kallenberg$^{(28)}$),
this implies  the following chain
of equalities for each $\varphi \in C_{c}(S)$:
\begin{eqnarray*}
	\mbox{exp}\left\{ -\int_{}^{}R_{\infty}^{1}(d\mu)
(1-e^{-\langle \mu,\varphi
	\rangle})\right\} & = & \int_{}^{} e^{-\langle
\mu,\varphi\rangle}
	L_{\infty}^{1} \left(d\mu\right) \\
	& = & \int\!\!\int_{}^{}e^{-\langle \mu,\varphi
	\rangle} \Pi_{\nu}\left(d\mu\right) L_{\infty}\left(d\nu\right) \\
	 & = & \int_{}^{}e^{-\langle \nu,1-e^{-\varphi} \rangle}
	 L_{\infty}\left(d\nu\right) \\
	 & = &  \mbox{exp}\left\{ -\int_{}^{}R_{\infty}\left(d\nu \right)
	 (1-e^{-\langle \nu,1-e^{-\varphi} \rangle})\right\}  \\
	 & = & \mbox{exp}\left\{  -\int_{}^{}R_{\infty}\left(d\nu\right)
	 \biggl(1- \int_{}^{} e^{-\langle \mu,\varphi
\rangle}\Pi_{\nu}(d\mu)\biggr)\right\}\\
	 & = & \mbox{exp}\left\{ -\int_{}^{}R_{\infty}\left(d\nu\right)
	 \int_{}^{}\Pi_{\nu}\left(d\mu\right) (1-e^{-\langle
\mu,\varphi \rangle})\right\},
\end{eqnarray*}
which yields (A.1.10).

Let ${\cal M}^2_\tau(S)$ denote the space of Radon measures
$\underline{\mu}$ on ${\cal M}_\tau(S)$ such that
$\langle\!\langle\underline{\mu},\langle\cdot,\tau\rangle
\rangle\!\rangle<\infty$, where
$$\langle\!\langle
\underline{\mu},F\rangle\!\rangle
=\int_{{\cal M}_\tau(S)}F(\nu)\underline{\mu}(d\nu)$$
(sometimes we use the notation on the r.h.s in order to avoid confusion).
We have that $R^1_\infty\in{\cal M}^2_\tau(S)$, $R^1_\infty$ is invariant
(but not reversible)
 for the 1--level dynamics (Liemant et al$^{(31)}$, Chapter 2,
Dawson and Perkins$^{(11)}$), and it has intensity $\rho$ in the sense that
\vskip.2cm
\noindent
$\hfill
\langle\!\langle
 R^1_\infty,\langle\cdot,\varphi\rangle\rangle\!\rangle
=\langle\rho,\varphi\rangle,\quad\varphi\in{\cal C}^+_\tau(S).\hfill 
(A.1.11)$
\vskip.2cm
\noindent
If $\beta =1$, then $ R^1_\infty$ has finite moments of all orders and the 
second
 and third moments are given by
\setcounter{equation}{11}
\def\theequation{A.1.\arabic{equation}}
\begin{eqnarray}
\langle\!\langle
 R^1_\infty,
\langle\cdot,
\varphi\rangle\langle\cdot,\psi\rangle\rangle\!\rangle
&=&
\langle\rho,\varphi\psi\rangle+
{V\over 2}\langle\rho,\varphi G\psi\rangle,
\quad \varphi,\psi\in{\cal C}^+_c(S), \\
&&\mkern -260mu
\langle\!\langle R^1_\infty,\langle\cdot,\varphi\rangle\langle\cdot,\psi
\rangle
\langle\cdot,\zeta\rangle\rangle\!\rangle =
\langle\rho,\varphi\psi\zeta\rangle+{V\over 2}\langle\rho,
\varphi\psi G\zeta+\varphi\zeta G\psi+\psi\zeta G\varphi\rangle\nonumber\\
&&\mkern -260mu
+{ V^2\over 2}\left\langle\rho,\int^\infty_0
[\varphi GT_t(T_t\psi\cdot T_t\zeta)+\psi GT_t(T_t\varphi\cdot T_t\zeta)+
\zeta G T_t(T_t\varphi\cdot T_t\psi)]dt
\right\rangle,\quad\varphi,\psi,\zeta\in{\cal C}^+_c(S).
\def\teheequation{A.1.13}
\end{eqnarray}
See Subsection A.4 for a proof.

Note also, from (A.1.9) and $T_t$--invariance of $\rho$, that for each 
$t>0$,
\vskip.3cm
\noindent
$\hfill
\langle\!\langle R^1_t, 
\langle\cdot,\varphi\rangle\rangle\!\rangle=\langle\rho,\varphi\rangle,\  
quad\varphi\in{\cal C}^+_\tau(S).\hfill (A.1.14)
$\vskip.3cm

We pass now to the 2--level branching system.
 A ``2--level particle'' is an element $\mu$ of ${\cal
N}_\tau(S)$ of the form $\mu=\sum^n_{i=1}\delta_{x_i}$.
A ``clan'' is the progeny under the 1--level dynamics (i.e., individual 
particles undergoing $(1+\beta)$--branching at rate $V_1$) of a family of 
particles which constitute an initial 2--level particle. Clans undergo 
$(1+\beta_2)$--branching at rate $V_2$.
Assuming persistence of the 1--level system,
the 2--level system starts off from a
Poisson system of ``2--level particles'' with intensity measure
$R^1_\infty$.
The empirical measures
of the 2--level system take values in ${\cal M}^2_\tau(S)$. Restricting to
test functions on ${\cal M}_\tau(S)$ of
the form $\mu\mapsto\langle\mu,\varphi\rangle$,  $\varphi\in{\cal
C}^+_c(S)$, amounts to considering the {\it aggregated}
system,
i.e.,
we consider the empirical measure of all the point masses disregarding
which 2--level particles they belong to. Note that the moments of
$R^1_\infty$ in (A.1.11), (A.1.12), (A.1.13) correspond to the aggregated 
Poisson
system. The empirical measures $X_t$ of the aggregation of the
2--level system take values in ${\cal N}_\tau(S)$.

The same argument used to obtain the Laplace functional of the 1--level
system can be used for the 2--level system.
Hence, analogously to (A.1.1),
 the Laplace functional of the occupation time of the 2--level system is 
given by
\vskip.2cm
\noindent
$\hfill
E\mbox{\rm exp}\left\{-
\left\langle\displaystyle\int^t_0
X_sds,\varphi\right\rangle\right\}=\mbox{\rm exp}\{-
\langle\!\langle R^1_\infty,{\mathbf u}_\varphi(t)\rangle\!\rangle\},
\quad\varphi\in{\cal C}^+_\tau(S),\hfill (A.1.15)$
\vskip.2cm
\noindent
where ${\mathbf u}_\varphi$ with values in $ [0,1]$
 is the unique solution of the non-linear evolution equation
\vskip.2cm
\noindent
$\hfill
{\mathbf u}_\varphi(t)=-
\displaystyle{V_2\over 1+\beta_2}
\displaystyle\int^t_0U_{t-s}({\mathbf u}_\varphi(s)^{1+\beta_2})ds+
\displaystyle\int^t_0U_{t-s}(\langle\cdot,\varphi\rangle(1-
{\mathbf u}_\varphi(s)))ds.\hfill (A.1.16)$
\vskip.2cm
\noindent
If follows from (A.1.4) and (A.1.16) that
\vskip.2cm
\noindent
$\hfill
{\mathbf u}_\varphi(t)(\mu)\leq
\displaystyle\int^t_0U_s(\langle\cdot,\varphi\rangle)
(\mu)ds=\displaystyle\int^t_0\langle\mu,T_s\varphi\rangle
ds=\langle\mu,G_t\varphi\rangle ,\quad \mu\in{\cal M}_\tau(S).
\hfill (A.1.17)$
\vskip.5cm
\noindent
{\bf A.2. The Palm formula}
\vglue.5cm
Let $M$ be a measure on ${\cM}_{\tau}(S)$ whose intensity measure
$$\Lambda_{M}:= \int_{{\cal M}_\tau(S)}\nu(\cdot) M(d\nu)$$
is locally finite, i.e.
$\langle\Lambda_{H}, \varphi\rangle < \infty$ for all $\varphi \in
\cC_c(S)$.
The Palm measures of $M$ are a family $(M_{x})_{x\in S}$ of
probability measures on $\cM_{\tau}(S)$ which satisfy
\vglue.3cm
\noindent
$\hfill
\displaystyle\int_{{\cal M}_\tau(S)}
\langle \nu, \varphi\rangle F(\nu) M(d\nu) =
\displaystyle\int_S\varphi(x)
\displaystyle\int_{{\cal M}_\tau(S)} F(\nu) \, M_{x}(d\nu) \Lambda_{M}(dx)
\hfill (A.2.1)
$\vskip.3cm
\noindent
for all measurable $F:\cM_{\tau}(S) \to \R_{+}$ and $\varphi:S\to
\R_{+}$.
If $M$ is supported by $\{\mu\in {\cal M}_{\tau}(S)\,|\,\mu \mbox{ is }
\{0,1,2,\ldots,\infty\}$--valued$\}$, then
$$M_{x}(\{\mu |\, \mu(x)\geq 1\})=1$$
for $\Lambda_{M}$--almost all $x$, and in this case the reduced Palm 
measures
$M^{red}_{x}$, $x\in S$, are defined by
\vglue.3cm
\noindent
$\hfill
M^{red}_{x}= M_{x}(\{\mu-\delta_{x}\in (\cdot)\}).\hfill
(A.2.2)
$
\vskip.5cm
\noindent
{\bf A.3.  Tree representations of the Palm measures
of $R^1_t$ and $R^1_\infty$}
\vglue.5cm
The Palm measures of $R_{t}^{1}$ and of the equilibrium canonical
measure $R_{\infty}^{1}$ described in Subsection A.1  have a
representation in terms of a backward tree which we recall here.

\vspace{1cm}

\noindent
\textbf{Lemma A.3.1.} (Gorostiza and Wakolbinger$^{(21)}$, Theorem 2.3).
Let $W^{x}$ be a random path of the motion process starting in $x\in
S$, let $\pi$ be a random Poisson configuration on $\R_{+}$ with
intensity $V$, and for each $y\in S$,  $r>0$ and $i=1,2,\ldots$,
let $X_{r}^{y,i}$
be a branching particle system arising from one
ancestor at site $y$ and developing over time $r$. Let
$N_{r}$ be an integer--valued random variable with generating
function $1-(1+\beta) q (1-s)^{\beta}$
(see Subsection 2.5). Assume
all these objects are independent.
Then the particle systems
\vskip.3cm
\noindent
$\hfill
\Phi_{x}^{t}:=\displaystyle\int_{0}^{t}
\sum_{i=1}^{Nr}X_{r}^{W_{r}^{x},i}(\cdot)\, \pi(dr)\hfill (A.3.1)
$\vskip.3cm
\noindent
and
\vskip.3cm
\noindent
$\hfill
\Phi_{x}^{\infty}:=\displaystyle\int_{0}^{\infty}
\displaystyle\sum_{i=1}^{Nr}X_{r}^{W_{r}^{x},i}(\cdot)\,
\pi(dr) \hfill (A.3.2)
$\vskip.3cm
\noindent
have distributions $(R_{t}^{1})^{\mbox{{\tiny \it red}}}_{x}$ and
$(R_{\infty}^{1})^{\mbox{{\tiny \it red}}}_{x}$, respectively,
for $\rho$--almost all $x\in S$.

\vspace{1cm}

It follows immediately from (A.3.1), (A.3.2),
(A.1.11), (A.1.14) and the Palm formula
(A.2.1) that all the moments of $R_{t}^{1}$ increase
 to those of $R_{\infty}^{1}$  as $t\to\infty$.
\vskip.5cm
\noindent
{\bf A.4. Second and third moments of $R^1_\infty$}
\vglue.5cm
\noindent
{\it Proof of (A.1.12) and (A.1.13):}
\vglue.5cm
The proof of  can be carried out directly by using
the explicit form of the Laplace transform of $R_{\infty}^{1}$ (as
given, e.g., in Gorostiza and Wakolbinger$^{(21)}$, Theorem 3.3).
Here we give an argument which uses the
structure of $R_{t}^{1}$  in (A.1.9) and the monotone
convergence  of the moments of $R^1_t$ mentioned above.

Let us introduce the  notation
\vskip.3cm
\noindent
$\hfill
A_{t,x}(\varphi, \psi)= \displaystyle\int_{0}^{t}T_{s}(T_{t-s}\varphi\cdot
T_{t-s}\psi)(x)ds, \hfill (A.4.1)
$\vskip.3cm
\noindent
$\hfill
B_{t,x}(\varphi,\psi,\zeta)=
\displaystyle\int_{0}^{t}T_{s}\left(T_{t-s}\varphi
\left(\displaystyle\int_{0}^{t-s}
T_{t-s-r}\left(T_{r}\psi\cdot T_{r}\zeta\right)
dr
\right)\right)(x)ds.\hfill (A.4.2)
$\vskip.5cm
The following formulae, which can be obtained either by
differentiating the Laplace functional of $X_{t}^{x}$ or by means of a
tree representation of the Palm measures of the distribution of
$X_{t}^{x}$ (similar to (A.3.1)), are well
known (e.g. Klenke$^{(30)}$, Lemma 3.1),
\setcounter{equation}{2}
\def\theequation{A.4.\arabic{equation}}
\begin{eqnarray}
& & E[\langle X_{t}^{x}, \varphi\rangle \langle X_{t}^{x},
\psi\rangle]=T_{t}(\varphi\psi)(x) + VA_{t,x}(\varphi,\psi), \\
\nonumber\\[-.5cm]
& & E[\langle X_{t}^{x}, \varphi\rangle \langle X_{t}^{x}, \psi\rangle
\langle X_{t}^{x}, \zeta\rangle]\nonumber \\
 && = T_{t}(\varphi\psi\zeta)(x)+V(A_{t,x}(\varphi,
\psi\zeta)+ A_{t,x}(\psi,
\varphi\zeta)+ A_{t,x}(\zeta, \varphi\psi))\nonumber\\
&& +V^{2}(B_{t,x}(\varphi,\psi,\zeta)+B_{t,x}
(\psi,\varphi,\zeta)+B_{t,x}(\zeta,\varphi,\psi)).
\end{eqnarray}
Note that
(A.4.3) is just a rewriting of (A.1.7). We need the following lemma.

\vspace{1cm}

\noindent
{\bf Lemma A.4.1.}

	\textrm{(a)} $\displaystyle\int_S
A_{t,x}(\varphi,\psi)\rho(dx)
	\mathop{\longrightarrow}
        \displaystyle\frac{1}{2}\langle \rho, \varphi G\psi\rangle$ as
        $t\rightarrow\infty$.
\vglue.5cm
	\textrm{(b)}  $\displaystyle\int_S B_{t,x}(\varphi, \psi,\zeta)\rho(dx)
	\mathop{\longrightarrow}
        \displaystyle\frac{1}{2}
        \left\langle \rho, \varphi
         \int_{0}^{\infty}GT_{r}(T_{r}\psi\cdot T_{r}\zeta)dr\right\rangle$
         as $t\rightarrow\infty$.

\vspace{1cm}

\noindent
{\it Proof:}
\vglue.5cm
\noindent
\textrm{(a)}\qquad\quad  $\displaystyle\int_S 
A_{t,x}(\varphi,\psi)\rho(dx)=
\left\langle \rho, \displaystyle\int_{0}^{t} (T_{s} \varphi\cdot
	T_{s}\psi)ds\right\rangle=
\left\langle\rho, \varphi \displaystyle\int_{0}^{t} \,
	T_{2s}(\psi)ds\right\rangle \mathop{\longrightarrow}
        \displaystyle\frac{1}{2}\langle \rho, \varphi G \psi\rangle$
        as $t\rightarrow\infty.\qquad\qquad$

\vglue.5cm
\noindent
(b)
\vglue-1.5cm
\begin{eqnarray*}
\displaystyle\int_S B_{t,x}(\varphi,\psi,\zeta)\rho(dx)
 & = & \left\langle
		\rho, \int_{0}^{t} \  T_{s} \varphi
\left(\int_{0}^{s}
		 \
		T_{s-r}(T_{r}\psi\cdot T_{r}\zeta)dr\right)ds\right\rangle\\
		 & = & \left\langle \rho , \varphi\int_{0}^{t}
\int_{0}^{s} \
		 T_{2s-r}(T_{r}\psi\cdot T_{r}\zeta)drds\right\rangle\\
		 & = & \left\langle \rho, \varphi \int_{0}^{t} \
\int_{r}^{t}
		 T_{2s-r}(T_{r}\psi\cdot T_{r}\zeta)dsdr\right\rangle\\
		 & = & \frac{1}{2}\left\langle \rho, \varphi
\int_{0}^{t} \
		 \int_{r}^{2t-r} \   T_{u}(T_{r}\psi
\cdot T_{r}\zeta)dudr\right\rangle\\
		 & = & \frac{1}{2}\left\langle \rho, \varphi
\int_{0}^{t} \
		 \int_{0}^{2(t-r)} \   T_{v}T_{r}(T_{r}\psi\cdot
T_{r}\zeta)dvdr\right\rangle \\
		 & 	\mathop{\longrightarrow}
                 & \frac{1}{2}\left\langle \rho,
\varphi\int_{0}^{\infty} \   G
	     T_{r}(T_{r}\psi\cdot T_{r}\zeta)dr\right\rangle\quad\mbox{\rm 
as}\quad t\rightarrow\infty.
\end{eqnarray*}
\vskip-1cm
$\hfill\Box$
\vglue.5cm
Combining Lemma A.4.1 with (A.1.9), (A.4.3), (A.4.4), and
using the above mentioned monotonicity of the moments of $R^1_t$,
the proof of (A.1.12) and (A.1.13) is complete.
$\hfill \Box$
\vskip.5cm
\noindent
{\bf ACKNOWLEDGEMENTS}
\vskip.5cm
This research was carried out mainly during visits to The Fields Institute 
(Toronto, Canada), the Johann Wolfgang Goethe University (Frankfurt, 
Germany), and the Center for Mathematical Research (CIMAT, Guanajuato, 
Mexico). We express our gratitude for the hospitality of these 
institutions. We thank Dr. Mike Porter, who helped us  by his computer 
analysis
to unravel the oscillations of the hierarchical random walk.
=====
\vskip.5cm
\noindent
{\bf REFERENCES}

\begin{enumerate}
\item Barlow, M. T. and  Perkins, E. A. (1988).
 Brownian motion on the Sierpi\'nski gasket,
{\it Probab. Th. Rel. Fields} {\bf 79} 543--623.

\item Cartwright, D. I. (1988).  Random walks on direct sums of discrete 
groups, {\it J. Theor. Probab.}  {\bf 1}  341--356.

\item Collet, P. and Eckmann, J-P. (1978).  A renormalization
group analysis of the hierarchical model in statistical mechanics,
{\it  Lecture Notes in Physics} {\bf  74}, Springer-Verlag, Berlin-New 
York.

\item  Cox, J. T. and Griffeath, D. (1984).  Large deviations for
Poisson
systems
of independent random walks, {\it Probab. Th. Rel. Fields} {\bf 66} 
543-558.

\item  Cox, J. T. and   Griffeath, D. (1985).  Occupation times for
critical
branching Brownian motions, {\it Ann. Probab.}
{\bf 13} 1108-1132.

\item  Dawson D. A. and  Gorostiza, L. G. (1990).
 Generalized solutions of a class of nuclear--space--valued stochastic
evolution equations, {\it Appl. Math. Optim.} {\bf 22} 241--263.

\item  Dawson, D. A. and  Greven, A.  Multiple space--time
scale
analysis for interacting branching models, {\it Electron. J.  Probab.}
{\bf 1}, no. 14, 84 pp.

\item Dawson, D. A. and  Ivanoff, G. (1978). Branching diffusions
and random
measures, in {\it Branching Processes},  Joffe, A.  Ney P. (editors), M.
Dekker,
New York,  61--103.

\item  Dawson,D. A. and  Hochberg, K. J. (1991).  A multilevel
branching
model, {\it Adv. Appl. Prob.} {\bf 23} 701--715.

\item Dawson,D. A., Hochberg, K. J. and Vinogradov, V. (1996). 
High--density limits of hierarchically structured branching--diffusing 
populations, {\it Stoch. Proc. Appl.} {\bf 62} 191--222.

\item Dawson, D. A. and Perkins, E. (1991).
Historical Processes. {\it Memoirs of the AMS} {\bf 454}, Providence, R. I.

\item Dawson, D. A. and Perkins, E. (1999). Measure--valued processes and 
renormalization of branching particle systems, in {\it Stochastic Partial 
Differential Equations: Six Perspectives}, Carmona, R. A. and Rozovskii, B. 
(editors), {\it Math. Surveys and Monographs} {\bf 64} 45--106, AMS.

\item Deuschel, J. D. and Rosen, J. (1998). Occupation time large 
deviations for critical branching Brownian motion, super--Brownian motion 
and related processes, {\it Ann. Probab.} {\bf 26} 602--643.

\item Deuschel, J. D. and Wang, K. (1994).  Large deviations for the
occupation
time functional of a Poisson system of independent Brownian
particles, {\it Stoch.}
{\it Proc. Appl.}  {\bf 52}  183-209.

\item  Fleischmann, K. and Greven, A. (1994).  Diffusive clustering
in an
infinite system of hierarchically interacting diffusions, {\it Probab.
Th. Rel.
Fields} {\bf 98} 517--566.

\item Gorostiza, L.G. (1996).  Asymptotic fluctuations and
critical
dimension for a two-level branching system, {\it Bernoulli} {\bf 2}
109-132.

\item  Gorostiza, L. G. Hochberg, K. and  Wakolbinger, A. (1995).
Persistence of a
critical super-2 process. {\it J. Appl. Prob.} {\bf 32}  534-540.

\item  Gorostiza, L. G.  and  L\'opez--Mimbela, J. A. (1994).
An occupation time approach for convergence of measure--valued processes, 
and the death process of a branching system. {\it Stat. Prob. Lett.} {\bf 
21} 59--67.

\item  Gorostiza, L. G., Roelly--Coppoletta, S. and Wakolbinger, A. 
 (1990).
Sur la persistence du processus de Dawson--Watanabe stable, in {\it 
S\'eminaire de Probabilit\'es XXIV}, Az\'ema, J, Meyer, P. A. and Yor, M. 
(editors).
{\it Lect. Notes Math.} {\bf 1426}, Springer--Verlag, Berlin, 275--281.

\item Gorostiza, L. G. and Rodrigues, E.R. (1999). A stochastic model for 
transport of particulate matter in air: an asymptotic analysis, {\it Acta 
Applicandae Mathematicae} (in press).

\item Gorostiza, L. G. and  Wakolbinger, A. (1991).  Persistence
criteria
for a class of critical branching particle systems in continuous
time, {\it Ann.
Probab.} {\bf 19}  266--288.

\item Gorostiza, L. G. and   Wakolbinger, A. (1994).
 Long time behavior of
critical particle systems and applications, in {\it Measure--Valued 
Processes, Stochastic Partial Differential Equations, and Interacting 
Systems}, Dawson, D. A. (editor), {\it CRM Proc. and
Lecture Notes}
{\bf 5}, AMS, 119--137.

\item Grabner, P. J. and  Woess, W. (1997).
 Functional iterations and periodic oscillations for simple random walk on 
the Sierpi\'nski gasket, {\it Stoch. Proc. Appl.} {\bf 69}  127--138.

\item Greven, A. and Hochberg, K. J. (1999). New behavioral patterns for 
two--level branching systems, in {\it Stochastic Models}, Gorostiza, L. G. 
and Ivanoff, B. G. (editors), {\it Conference Procedings Series}, CMS (in 
press).

\item Hochberg, K. J. (1995). Hierarchically structured branching 
populations with spatial motion. {\it Rocky Mountain J. Math.} {\bf 25} 
269--283.

\item  Hochberg, K. J. and  Wakolbinger, A. (1995).
Non-persistence of two-level
branching particle systems in low dimensions, in {\it Stochastic Partial 
Differential Equations}, Etheridge, A. (editor),
 {\it London Mathem.
Soc. Lecture
Note Series} {\bf 216},  Cambridge Univ. Press, 126--140.

\item Iscoe, I. (1986).   A weighted occupation time for a class
of
measure-valued branching processes, {\it Probab. Th. Rel. Fields}
{\bf 71} 85 -- 116.

\item Kallenberg, O. (1983). Random Measures, 3rd. Edition, 
Akademie--Verlag, Berlin, Academic Press, New York.

\item Kesten, H. and Spitzer, F. (1965). {\it Random walks on countably 
infinite Abelian groups}, Acta Math. {\bf 114} 237--265.

\item Klenke, A. (1997). Multiple scale analysis of clusters in spacial 
branching models. {\it Ann. Probab.} {\bf 25}, 1670--1711.

\item  Liemant, A.,  Matthes, K. and  Wakolbinger, A. (1998).
Equilibrium
Distributions of Branching Processes, {\it Akademie-Verlag}, Berlin and
Kluwer,
Dordrecht.

\item  M\'el\'eard, S. and  Roelly, S. (1992).
 An ergodic result for critical spatial branching systems,
{\it Stochastic Analysis and Related Topics,
Progress in Probability} {\bf 31}, Birkhauser, Boston,  333--341.

\item Port, S. C. and   Stone, C. J. (1971).
Infinitely divisible processes and their
potential theory (First Part), {\it Ann. Inst. Fourier} {\bf 21} (2)
157--275.

\item Sato, K. (1996). Criteria of weak and strong transience
for L\'evy
processes, in {\it Probability Theory and Mathematical Statistics,
Proceedings
of the Seventh Japan-Russia Symposium}, World Scientific, Singapore,
438--449.

\item  Sawyer, S. and  Felsenstein, J. (1983).  Isolation by
distance in a
hierarchically clustered population, {\it J. Appl. Prob.} {\bf 20}
1--10.

\item Sinai, Ya.  G. (1982).  Theory of Phase Transitions: Rigorous
results. Pergamon Press.

\item Spitzer, F. (1964).  Principles of Random Walk.
Van Nostrand, Princeton.

\item  St\"ockl, A. and   Wakolbinger, A. (1994).  On clan-recurrence
and
-transience in time stationary branching Brownian particle systems,
in {\it Measure--Valued Processes, Stochastic Partial Differential 
Equations, and Interacting Systems}, Dawson, D. A. (editor),
{\it CRM Proc. and Lecture Notes} {\bf 5},  AMS, 213--219.

\item Wilson, K. (1976).  The renormalization group and block
spins, in {\it Proceedings of the International Conference on
Statistical Physics},  P\'{a}l, L.  and  Sz\'{e}pfalusy, P. (editors), 
North
Holland, Amsterdam.

\item Wu, Y. (1994). Asymptotic behavior of two level measure branching 
processes, {\it Ann. Probab.} {\bf 22} 854--874.

\end{enumerate}

\end{document}